\documentclass[11pt]{article}
\usepackage[all]{xy}
\usepackage{amsmath}
\usepackage{amsfonts}
\usepackage{amssymb}
\usepackage{amscd}
\usepackage{amsthm}
\usepackage{latexsym}
\usepackage{amsbsy}
\newtheorem{lem}{Lemma}[section]
\newtheorem{prop}{Proposition}[section]
\newtheorem{theorem}{Theorem}[section]

\newtheorem{definition}{Definition}[section]
\newtheorem{example}{Example}[section]
\newtheorem{examples}{Examples}[section]

\newtheorem{remark}{Remark}

\newtheorem{exercise}{Exercise}[section]

\newcommand{\ot}{\otimes}
\newcommand{\A}{\mathcal{A}}

\begin{document}
\title{Lectures on Noncommutative Geometry}
\author {Masoud Khalkhali}
\date{}
\maketitle
\tableofcontents
\begin{abstract}

 \end{abstract}This text is an introduction to a few selected areas of Alain Connes' noncommutative geometry
 written for the volume of the school/conference ``Noncommutative Geometry 2005" held at IPM
 Tehran. It is an expanded version of my lectures  which was directed at
 graduate students and novices in the subject.

\section{Introduction}
David Hilbert once said  ``No one will drive us from the paradise
which Cantor created for us"  \cite{hilbert}. He was of course
referring to set theory and the vast new areas of  mathematics
which were made possible
  through it. This was  particularly relevant for geometry where our ultimate
  geometric intuition of {\it space}  so far
   has been a set endowed  with some extra structure (e.g.  a topology, a measure, a smooth structure,
   a sheaf,  etc.).
  With Alain Connes'  noncommutative geometry \cite{ac80, co85, cob}, however, we are now gradually moving into a
  new  `paradise', a `paradise' which contains the `Hilbertian paradise' as  one of its old
  small neighborhoods.
  Interestingly enough though,  and I hasten to say this,  methods of functional analysis,
  operator algebras, and spectral theory,  pioneered by Hilbert and his disciples,
  play a big role in Connes'  noncommutative geometry.

  The following text is a greatly  expanded version of talks I gave
  during the  conference on noncommutative geometry  at IPM, Tehran. The
  talks were directed  at graduate students,
  mathematicians, and physicists  with no background in noncommutative geometry. It inevitably covers only 
  certain selected parts
   of the subject   and many important topics  such as: metric and spectral aspects of noncommutative geometry, the local
  index formula, connections with number theory, and interactions with
  physics, are left out. For an insightful and comprehensive  introduction to the current state of  the art
    I refer to Connes and Marcolli's  article `A Walk
  in the Noncommutative Garden' \cite{coma2} in this volume. For a deeper
  plunge into  the subject one   can't do better than directly  going to Connes' book \cite{cob} and original articles.

  I would like to thank  
      Professors Alain Connes and  Matilde Marcolli for their support
      and encouragement over a long period of time.   
   Without their kind support  and advice this project
   would have taken much longer to come to its  conclusion, if ever. My
   sincere thanks go also to 
    Arthur Greenspoon   who took a keen interest in the
   text and  kindly and carefully  edited the entire manuscript. Arthur's superb skills resulted in
     substantial improvments  in  the original text.

\section{From $C^*$-algebras to noncommutative spaces}

Our working definition of a  {\it noncommutative space} is a
noncommutative algebra, possibly endowed with some extra
structure. Operator algebras, i.e. algebras of bounded linear operators
on a Hilbert space,   provided the first really deep insights into
this noncommutative realm. It is generally agreed that  the
classic series of papers of Murray and von Neumann starting with \cite{muvo}, and Gelfand and
Naimark \cite{gena} are the foundations upon which the theory of operator
algebras is  built. The first is the birthplace of
  von Neumann algebras as the  noncommutative counterpart of measure
  theory,  while in
  the second   $C^*$-algebras were shown to be the noncommutative analogues of locally compact
  spaces. For lack of space we shall say nothing about von Neumann
  algebras and their  place in noncommutative geometry (cf. \cite{cob} for the general theory as well as links with
 noncommutative geometry).
  We start with the definition of $C^*$-algebras and results of Gelfand
  and Naimark. References include \cite{bl2, cob, da, fi}.

\subsection{Gelfand-Naimark theorems}

By an {\it algebra} in these notes  we shall mean an {\it associative algebra} over the field of complex
numbers $\mathbb{C}$.  Algebras are not assumed to be commutative or unital, unless explicitly
specified so.  An {\it involution} on
an algebra $A$ is a conjugate linear map $a \mapsto a^*$
satisfying
$$(ab)^\ast = b^\ast a ^\ast \quad \text{and} \quad
(a^\ast )^\ast = a $$
 for all $a$ and $b$ in $A$.
By a {\it normed algebra} we mean an algebra
$A$ such that $A$ is a normed vector space and
\begin{eqnarray*}
 \| ab \| \le \| a \|\| b \|,
\end{eqnarray*}
for all $a, b$ in $A$.
If $A$ is unital, we assume that $\| 1 \| = 1$.
A {\it Banach algebra} is a  normed algebra which is {\it complete} as a metric space. One of the main consequences of completeness is that
norm convergent series are convergent; in particular if $\| a \|<1$ then  the geometric series $\sum_{n=1}^{\infty}a^n $ is
convergent.  From this it easily follows that the group
of invertible elements of a unital Banach algebra is  open in the norm topology.

\begin{definition}
A {\it $C^\ast$-algebra} is  an involutive Banach algebra $A$ such that for all $a
\in A$  the {\it $C^*$-identity}
\begin{equation}\label{Cid}
\| aa^\ast \| = \| a \|^2
\end{equation}
holds.
\end{definition}
A {\it morphism} of  $C^*$-algebras is an algebra homomorphism
$ f: A \longrightarrow B $
which   preserves the $*$ structure, i.e.
  $$ f(a^*)=f(a)^*, \quad \text{for all} \; \, a\in A.$$

 The $C^*$-identity \eqref{Cid}  puts     $C^*$-algebras
in a  unique place among all Banach algebras, comparable  to the
unique position enjoyed by Hilbert spaces among all Banach spaces.
Many facts  which are true for $C^*$-algebras are not necessarily
true for  an arbitrary  Banach algebra. For example, one can show,
using the spectral radius formula $\rho (a)= \text{Lim}\,\, \|a^n
\|^{\frac{1}{n}}$ coupled with the $C^*$-identity, that the norm
of a $C^*$-algebra is unique.  In fact it can be shown that a
morphism of $C^*$-algebras is  automatically {\it contractive}
  in the sense that for all $a \in A$,
  $\| f(a)\| \le \| a\|.$
  In particular they are always  continuous. It follows that  if $(A, \, \|\,  \|_1)$ and
   $(A, \, \|\, \|_2)$ are both $C^*$-algebras then
$$\|a \|_1 =\|a \|_2,$$
for all $a\in A$. Note also that a morphism  of $C^*$-algebras
is an {\it isomorphism} if and only if it is one to one and onto.
Isomorphisms of $C^*$-algebras are necessarily {\it isometric}.

\begin{example} \rm{
let $X$ be a locally compact Hausdorff space and let $A=C_0(X)$
denote the algebra of continuous complex valued functions on $X$
 vanishing at infinity. Equipped with the sup norm and the
involution  defined by complex conjugation,  $A$ is easily seen to
be a commutative $C^*$-algebra. It is unital if and only if $X$ is
compact, in which case $A$ will be denoted by $C(X)$. By a fundamental theorem of Gelfand and Naimark, to be recalled below,
any commutative
$C^*$-algebra is of the form  $C_0(X)$ for a canonically defined locally compact Hausdorff space $X$.  }
\end{example}

\begin{example} \rm{ The algebra  $A= \mathcal{L} (H)$ of all bounded linear
operators on a complex Hilbert space $H$ endowed with the operator norm and the usual adjoint operation is a $C^*$-algebra.
The crucial $C^*$-identity $\|a a^*\|=\|a\|^2$ is easily checked. When $H$
is finite dimensional of dimension $n$  we obtain the matrix algebra
$A=M_n (\mathbb{C})$. A direct sum of matrix algebras
$$A=M_{n_1}(\mathbb{C})\oplus M_{n_2}(\mathbb{C})\oplus \cdots
\oplus M_{n_k}(\mathbb{C})$$
is a $C^*$-algebra as well. It can be shown that any finite dimensional
$C^*$-algebra is unital and is a direct sum of matrix algebras \cite{bl2, da}. In other words,
finite dimensional $C^*$-algebras are semisimple.

Any norm closed subalgebra of $\mathcal{L} (H)$ which is also
closed under the adjoint map is clearly a $C^*$-algebra. A nice
example is the algebra $\mathcal{K} (H)$ of {\it compact operators} on
$H$. By definition, a bounded
 operator $T : H \to H$ is called compact if it is the norm limit of a sequence of finite rank operators.
 By the second fundamental theorem of Gelfand
 and Naimark, to be discussed below, any $C^*$-algebra is realized as a subalgebra of
 $\mathcal{L} (H)$ for some Hilbert space $H$. }
\end{example}

Let $A$ be an algebra. A {\it character} of $A$ is a nonzero
multiplicative linear map
$$\chi: A \to \mathbb{C}.$$
 If $A$ is
unital then necessarily $\chi (1)=1$. Let $\widehat{A}$ denote the
set of characters of $A$. It is also known as the {\it maximal spectrum} of $A$.
 If $A$ is a Banach algebra it can be
shown that any character of $A$ is automatically continuous and
has norm one. We can thus endow $\widehat{A}$ with the weak*
topology inherited from $A^*$, the continuous dual of $A$. The unit ball of
$A^*$ is compact in the weak* topology and one can deduce from
this fact that $\widehat{A}$ is a locally compact Hausdorff space.
It is compact if and only if $A$ is unital.

When $A$ is a unital Banach algebra
  there is a one to one correspondence between characters
of $A$ and the set of maximal ideals of $A$: to a character $\chi$
we associate its kernel, which is a maximal ideal,  and to a maximal
ideal $I$ one associates the character $\chi: A \to A/I \simeq
\mathbb{C}$.
   If $A$ is a $C^*$-algebra a
character  is necessarily  a $C^*$ morphism.

 An arbitrary $C^*$-algebra may well
have no characters at all. This happens, for example,  for simple
$C^*$-algebras, i.e. $C^*$-algebras with no non-trivial closed two-sided ideal; a  simple  example of these are  matrix algebras
$M_n(\mathbb{C})$. A commutative $C^*$-algebra, however, has plenty
of characters to the extent that characters separate points of $A$
and in fact completely characterize it, as we shall see.
\begin{example} \rm{Let $X$ be a locally compact Hausdorff space.
For any $x\in X$ we have the {\it evaluation character}
$$\chi = \text{ev}_x: C_0(X) \longrightarrow  \mathbb{C}, \quad  \text{ev}_x (f)=f(x).$$
It is easy to see that all characters of $C_0(X)$ are of this form
and that the map
$$ X \to \widehat{C_0(X)}, \quad \quad x  \mapsto \text{ev}_x$$
is a homeomorphism.}
\end{example}

For any commutative Banach algebra $A$,   the {\it Gelfand transform}
$$ \Gamma : A \to C_0 (\widehat{A})$$
is defined by $\Gamma (a)=  \hat{a}$, where
$$\hat{a}(\chi)= \chi (a).$$
It is a  norm contractive algebra homomorphism,  as can be easily seen. In general $\Gamma$
need not be injective or surjective, though its image separates
the points of the spectrum. The kernel of $\Gamma$ is the {\it nilradical} of $A$ consisting of nilpotent
elements of $A$.

The  paper of Gelfand and
Naimark \cite{gena} is the birthplace of the theory of
$C^*$-algebras. Together with Murray-von Neumann's series of papers on von Neumann algebras \cite{muvo}, they form
 the foundation stone  of operator algebras.
  The following two fundamental results on the structure
 of $C^*$-algebras are proved in this paper. The first result is the
 foundation for the belief that noncommutative $C^*$-algebras  can  be regarded as noncommutative
 locally compact Hausdorff spaces.

\begin{theorem}\label{gn} (Gelfand-Naimark \cite{gena}) a) For any  commutative
$C^*$-algebra $A$  with spectrum $\widehat{A}$  the Gelfand transform
$$ A \longrightarrow C_0(\widehat{A}), \quad a\mapsto \hat{a},$$
defines  an isomorphism of $C^*$-algebras. \\
b) Any $C^\ast$-algebra is  isomorphic to a $C^\ast$-subalgebra of the algebra $\mathcal{L}(H)$
of bounded  operators on a Hilbert
space $H$.
\end{theorem}

Using part a), it is easy to see that the functors
$$ X \leadsto C_0(X),  \quad \quad A \leadsto \widehat{A}$$
define an equivalence between the category of locally
compact Hausdorff spaces and proper continuous maps and the opposite of
the category of commutative $C^*$-algebras and
proper  $C^*$-morphisms. Under this correspondence compact
Hausdorff spaces correspond to unital commutative $C^*$-algebras.
  We can therefore think of the opposite of the category
of  $C^*$-algebras as the  category of locally
compact noncommutative spaces.

\begin{exercise} Let $X$ be a compact Hausdorff space and  $x_0 \in X $. To test your understanding of
the Gelfand-Naimark theorem, give
a completely $C^*$-algebraic definition of the fundamental group $\pi_1
(X, x_0)$.
\end{exercise}

\subsection{GNS, KMS, and the flow of time }

We briefly indicate the proof of part b) of Theorem \eqref{gn} in
the unital case. It is based on the notion of  {\it state} of a
$C^*$-algebra and the accompanying `left regular
representation',   called the GNS (Gelfand-Naimark-Segal)
construction. We then look at the KMS (Kubo-Martin-Schwinger)
condition characterizing the  equilibrium states in quantum
statistical  mechanics and the time evolution  defined by a
state. In the context of von Neumann algebras, a fundamental result of
Connes states  that the time evolution is unique up to inner
automorphisms.

The concept  of state  is the  noncommutative analogue of Borel
probability measure.    A state of a unital $C^*$-algebra $A$ is a
positive normalized linear functional $\varphi : A \mapsto
\mathbb{C}$:
$$ \varphi (a^*a) \geq 0 \quad \forall a \in A,  \quad  \text{and} \quad \varphi(1)=1.$$
 The
{\it expectation value} of an element (an `observable') $a \in A$, when the system is
in the state $\varphi$, is defined by $\varphi (a)$.
  This terminology is motivated by states in statistical mechanics where  one abandons the idea of describing the state of a
system by  a point in the phase space. Instead, the only
reasonable question to ask is  the probability of  finding  the
system
 within a  certain region in the phase space. This probability is of course  given
by a probability measure $\mu$. Then the expected value of an
observable $f: M\to \mathbb{R}$, if the system is in the state
$\mu$,  is  $\int fd\mu$.

Similarly, in
quantum statistical mechanics the idea of describing the quantum
states of a system by a vector (or ray) in a Hilbert space is
abandoned and one instead uses a {\it density matrix}, i.e. a  trace class positive operator $p$
 with Tr$(p)=1$. The expectation
value of an observable $a$, if the system is in the  state $p$,
 is given by Tr$(ap)$.

A state is called {\it pure} if it is not a non-trivial convex
combination of two other states. This corresponds to point masses
in the classical case and to vector states in the quantum case.

\begin{example} \rm{  By  the Riesz representation theorem,  there is a one to one
correspondence  between states  on  $A=C(X)$  and Borel
probability measures on $X$,  given by
$$ \varphi (f)=\int_X fd\mu.$$
$\varphi$ is pure if and only if $\mu$ is a Dirac mass at some
point $x\in X$.}
\end{example}

\begin{example} \rm{
For
$A=M_n(\mathbb{C})$, there is a one to one correspondence between states
$\varphi$
of $A$ and positive matrices $p$ with Tr$(p)=1$ ($p$ is called a {\it density matrix}). It is  given by
$$\varphi (a)=\text{Tr} \,(ap).$$
$\varphi$ is pure if and only if $p$ is of rank 1.}
\end{example}

\begin{example} \rm{
For a different example,  assume $\pi :A\to
\mathcal{L}(H)$ is a {\it representation} of $A$ on a Hilbert
space $H$. This simply means $\pi$ is a morphism of $C^*$-algebras.
Given any unit vector $v\in H$, we can define a state on $A$  by
$$\varphi (a)= \langle \pi (a) v, \,  v \rangle.$$
Such states are called {\it vector states}. As a consequence of the GNS
construction, to be described next,  one knows that any state is a vector
state in the corresponding GNS representation.}
\end{example}

Let $\varphi$ be a positive linear functional on $A$. Then
$$\langle a, \, b \rangle:=\varphi (b^*a)$$
is a positive semi-definite bilinear form  $A$ and hence satisfies the {\it Cauchy- Schwarz}
inequality. That is,   for all $a, b$ in $A$ we have
$$ |\varphi(b^*a)|^2\leq \varphi(a^*a)\varphi(b^*b).$$
Let
$$N=\{ a\in A; \quad \varphi (a^*a)=0\}.$$
It is easy to see, using the above Cauchy-Schwarz inequality,
that $N$ is a closed left ideal  of $A$ and
$$\langle a+N, \, \,b+N \rangle:= \langle a, \, b \rangle,$$
is a  positive definite inner product on the quotient space $A/N$.
Let $H_{\varphi}$ denote the Hilbert space completion of $A/N$.
The {\it GNS representation}
$$\pi_{\varphi} : A \longrightarrow \mathcal{L}(H_{\varphi})$$
is, by definition,  the unique extension of the left regular
representation  $A \times A/N \to A/N, \, \, (a,\, \, b+N)\mapsto
ab+N$. Let $v= 1 + N$. Notice that we have
$$\varphi (a)= \langle (\pi_{\varphi}a)(v),\, \,v \rangle,$$
for all $a\in A$. This shows that the state
$\varphi $ can be recovered from the  GNS representation
 as a vector state.

\begin{example} \rm{For $A= C(X)$ and $\varphi (f)=\int_X fd\mu$, where
$\mu$ is a Borel probability measure on $X$, we have  $
H_{\varphi}= L^2 (X, \, \mu)$. The GNS representation is the
representation of $C(X)$ by multiplication operators on $ L^2 (X,
\, \mu)$. In particular, when $\mu$ is the Dirac mass at a point $x
\in X$, we have  $H_{\varphi} \simeq \mathbb{C}$ and
$\pi_{\varphi} (f)= f(x)$.}
\end{example}

The GNS representation $(\pi_{\varphi}, H_{\varphi})$ may fail to be
faithful. It can be shown that it is irreducible if and only if
$\varphi$ is a pure state \cite{bl2, da}. To construct a faithful
representation, and hence an embedding of $A$ into the algebra of
bounded operators on  a Hilbert space, one first shows that there
are enough pure states on $A$. The proof of the following result
is based on the Hahn-Banach and Krein-Milman theorems.

\begin{lem} For any positive element $a$ of $A$, there exists a pure
state $\varphi$ on $A$ such that $\varphi (a)=\|a\|.$
\end{lem}

Using the GNS representation associated to $\varphi$, we can then
construct, for any $a \in A$,  an irreducible representation $\pi$
of $A$ such that $\|\pi (a) \|=|\varphi (a)|=\|a\|.$

We can now prove the second theorem of Gelfand and Naimark.
\begin{theorem} Every $C^*$-algebra is isomorphic  to a
$C^*$-subalgebra of the algebra of bounded operators on a Hilbert space.
\end{theorem}
\begin{proof}
Let $\pi=\sum_{\varphi \in \mathcal{S}(A)} \pi_{\varphi}$ denote
the direct sum of all GNS representations for  all  states of $A$.
By the above remark $\pi$ is faithful.
\end{proof}

\begin{exercise} Identify the GNS representation in the case where $A=
M_n(\mathbb{C})$ and $\varphi $ is the normalized trace on $A$.
For which state is the standard representation on $\mathbb{C}^n$ 
obtained?
\end{exercise}

In the remainder of this section we look at the
Kubo-Martin-Schwinger (KMS) equilibrium condition for states  and
some of its consequences. KMS states  replace the Gibbs
equilibrium states for interacting systems with infinite degrees
of freedom.  See \cite{brro} for an introduction to quantum statistical mechanics; see also
\cite{coma1} and the forthcoming book of
Connes and Marcolli \cite{coma3} for relations between quantum statistical mechanics,
number theory and noncommutative geometry. For relations with
Tomita-Takesaki theory and Connes' classification of factors the best
reference is Connes' book \cite{cob}.

  A $C^*$-{\it dynamical system} is a triple $(A, \,G, \, \sigma)$ consisting of a $C^*$-algebra $A$, a
locally compact group $G$ and a continuous action
$$\sigma : G \longrightarrow  \text{Aut} \,(A)$$
of $G$ on $A$, where $\text{Aut}\, (A)$ denotes the
group of $C^*$-automorphisms of $A$.   The correct  continuity
condition for $\sigma$ is {\it strong continuity} in the sense
that for all $a \in A$ the map $g \mapsto \sigma_g  (a)$ from $G
\to A$ should be continuous. Of particular interest is the case
$G= \mathbb{R}$ representing a quantum mechanical system evolving in time. For
example, by  Stone's  theorem one knows that one-
parameter groups of  automorphisms of $A =\mathcal{L}
(\mathcal{H})$ are of the form
$$\sigma_t (a)=e^{itH} a e^{-itH},$$
where $H$, the {\it Hamiltonian} of the system,  is a self-
adjoint, in general  unbounded,  operator on $\mathcal{H}$. Assuming the operator $e^{-\beta H}$ is trace class,  the corresponding
  {\it Gibbs equilibrium state} at inverse temperature $\beta =\frac{1}{k T} > 0$ is the state
  \begin{equation} \label{gist}
   \varphi (a) =\frac{1}{Z (\beta)} \, \text{Tr} \, (a e^{-\beta H}),
   \end{equation}
  where the {\it partition function} $Z$ is defined  by
  $$ Z (\beta) = \,\text{Tr}\, (e^{-\beta H}).$$

According to Feynman \cite{fe},  formula \eqref{gist}  for the
Gibbs equilibrium state (and its classical analogue) is the apex
of statistical mechanics. It should however be added that
\eqref{gist} is not powerful enough to deal with interacting
systems with an infinite number of degrees of freedom (cf. the first chapter  of
Connes' book \cite{cob} for an example), and in general  should be
replaced by the KMS equilibrium condition.

 Let   $(A, \, \, \sigma_t)$ be a $C^*$-dynamical system evolving in time. A state $\varphi:  A \to \mathbb{C}$
 is called a {\it KMS state at
 inverse temperature} $  \beta >0 $  if for all $a,  b  \in  A$  there
 exists a function $F_{a,b}(z)$ which is continuous and bounded on the
 closed
 strip $ 0 \leq \text{Im}\, z \leq \beta$ in the complex plane and holomorphic in the interior such
 that for all $t \in \mathbb{R}$
 $$ F_{a, b}(t)= \varphi (a \sigma_t (b)) \quad \quad \text{and} \quad
 \quad F_{a, b}(t+ i\beta)= \varphi (\sigma_t (b) a).$$

Let $\mathcal{A} \subset A$ denote the set of {\it analytic
vectors} of $\sigma_t$ consisting of those elements $a\in A$ such
that $t\mapsto \sigma_t (a)$ extends to a holomorphic function on
$\mathbb{C}$. One shows that $\mathcal{A}$ is a dense
$*$-subalgebra of $A$. Now the KMS condition is equivalent to a {\it
twisted trace property} for $\varphi$: for all  analytic vectors
$a, b \in \mathcal{A}$ we have
$$ \varphi (b a) = \varphi (a \sigma_{i \beta } (b)).$$
Notice that the automorphism $\sigma_{i \beta }$ obtained by analytically continuing $\sigma_t$ to imaginary time (in fact imaginary
temperature!) is only densely defined.

\begin{example} \rm{ Any Gibbs state is a KMS state as can be easily
checked. }
\end{example}

\begin{example} (Hecke algebras,  Bost-Connes  and
Connes-Marcolli systems \cite{boco, coma1}) \rm{ A subgroup $\Gamma_0$ of a group $\Gamma$ is called
{\it almost normal} if every left coset $\gamma \Gamma_0$ is a
{\it finite} union of right cosets. In this  case we say
$(\Gamma, \, \Gamma_0)$ is  a {\it Hecke pair}. Let $L (\gamma)$ denote
the number of distinct right cosets $\Gamma_0 \gamma_i$  in the
decomposition
$$ \gamma \Gamma_0 =\cup_i \Gamma_0 \gamma_i,$$
and let $R (\gamma) = L(\gamma^{-1})$.

The rational {\it Hecke
algebra}  $\mathcal{A}_{\mathbb{Q}} = \mathcal{H}_{\mathbb{Q}} (\Gamma, \, \Gamma_0)$  of a
 Hecke pair $(\Gamma, \, \Gamma_0)$ consists of
 functions with finite support
 $$ f: \Gamma_0 \setminus \Gamma \to
\mathbb{Q} $$
 which are right $\Gamma_0$-invariant, i.e. $f (\gamma
\gamma_0) =f (\gamma)$ for all $\gamma \in \Gamma$ and $\gamma_0 \in
\Gamma_0$.
 Under the convolution product
$$(f_1 * f_2)(\gamma) :=\sum_{\Gamma_0 \setminus \Gamma} f_1 (\gamma \gamma_1^{-1})
f_2 (\gamma_1), $$
$\mathcal{H}_{\mathbb{Q}} (\Gamma, \, \Gamma_0)$ is an associative
unital algebra. Its complexification
$$ \mathcal{A}_{\mathbb{C}}=\mathcal{A}_{\mathbb{Q}}
\otimes_{\mathbb{Q}} \mathbb{C}$$
is a $*$-algebra with an involution given by
 $$f^*(\gamma):=  \overline{f (\gamma^{-1}}).$$
  Notice that if $\Gamma_0$ is normal in $\Gamma$ then one
obtains the group algebra of the quotient group $\Gamma /\Gamma_0$. We
refer to \cite{boco, coma1} for the $C^*$-completion of
$\mathcal{A}_{\mathbb{C}}$,  which is similar to the $C^*$-completion of
group algebras.
There is a one-parameter group of automorphisms of this Hecke algebra
(and its $C^*$-completion) defined by
$$ (\sigma_t f) (\gamma) = \left( \begin{array}{c} L (\gamma)\\
\overline{R(\gamma)}\end{array} \right)^{-it} f( \gamma).$$
 Let $P^+$ denote the subgroup of the ``$ax+b$'' group with $a >0$. The corresponding $C^*$-algebra
 for the Hecke pair $(\Gamma_0, \, \Gamma)$ where $\Gamma =P_{\mathbb{Q}}^+$ and
$\Gamma_0 =P_{\mathbb{\mathbb{Z}}}^+$ is the Bost-Connes
$C^*$-algebra. We refer to \cite{coma1} for a description of the
Connes-Marcolli system. One feature of these systems is that
their partition functions are expressible in terms of zeta  and
$L$-functions of number fields.}
  \end{example}

Given a state $\varphi$ on a $C^*$-algebra $A$ one may ask if there is a
one-parameter group of automorphisms  of $A$ for which $\varphi$ is a
KMS state at  inverse temperature $\beta =1$. Thanks to Tomita's
theory (cf. \cite{cob, bl2}) one knows that the answer is positive if  $A$ is a von
Neumann algebra which we  will  denote by $ M $ now.  The
corresponding automorphism group $\sigma_t^{\varphi}$, called the
{\it modular automorphism group}, is uniquely defined subject to the
condition $\varphi \sigma_t^{\varphi}= \varphi$ for all $t \in
\mathbb{R}$.

A von Neumann algebra  typically carries many states. One of the first achievements of Connes,  which set his
grand classification program of von Neumann algebras   in motion, was his
proof that   the modular automorphism  group is unique up to inner
automorphisms. More precisely   he showed that for any other state
$\psi$ on $M$ there is a continuous map $u$ from $\mathbb{R}$ to
the group of unitaries   of $M$ such that
$$\sigma_t^{\varphi} (x)=
u_t \sigma_t^{\psi} (x) u_t^{-1} \quad \text{and} \quad  u_{t+s} =u_t
\sigma_s^{\varphi}u_s.$$
It follows that  the modular
automorphism group is independent, up to inner automorphisms, of
the state (or weight) and if  $\text{Out} \,  (M)$ denotes the
quotient of the group of automorphisms of $M$ by inner
automorphisms, any von Neumann algebra  has  a god-given dynamical system
$$ \sigma : \mathbb{R} \to \text{Out} \,  (M)$$
 attached to it. This is a purely non-abelian  phenomenon as
the modular automorphism group is trivial for abelian von Neumann
algebras as well for type II factors.  For  type III factors it
turns out the the modular automorphism group possesses a  complete set of
invariants for the isomorphism type of the algebra  in the
injective case. This is the beginning  of Connes'  grand
classification theorems for von Neumann algebras,  for which we
refer the reader to  his book \cite{cob} and references therein.

\subsection{From groups to noncommutative spaces}
Many interesting $C^*$-algebras are defined as group
$C^*$-algebras or as crossed product $C^*$-algebras. Group
$C^*$-algebras are completions of group algebras with respect to
certain pre $C^*$-norms. To illustrate some of the general ideas
of noncommutative geometry and noncommutative index theory, we
shall sketch  Connes' proof  of `connectedness' of the group
$C^*$-algebra of free groups.

 \begin{example}(group $C^{\ast}-$algebras) \rm{
  To any locally compact topological group $G$ one  can associate two $C^*$-algebras,  the
 {\it full} and the {\it reduced} group $C^*$-algebras of $G$, denoted by $C^*
 (G)$ and $C^*_r (G)$, respectively. There is a 1-1
 correspondence between unitary representations of $G$ and the
 representations of $C^* (G)$ and a 1-1 correspondence between
 unitary representations of $G$ which are equivalent to a sub-representation of its left regular
 representation  and representations of   the reduced group
 $C^*$-algebra $C^*_r (G)$.
  Both algebras are completions of the convolution algebra of $G$ (under different norms).
 There is always a surjective $C^*$-morphism
 $C^* (G)\to C^*_r (G)$ which  is injective if and only if the group $G$ is
 {\it amenable}. It is known that abelian, solvable,  as well as compact groups are amenable while, for example,  non-abelian
 free groups are non-amenable. Here we consider only  discrete groups.

 Let $\Gamma$ be  a discrete group and
let $H=\ell^2 (\Gamma )$ denote the Hilbert space of square summable functions on
$\Gamma$.  It has an orthonormal basis consisting of delta functions $\{\delta_g\}$, $g\in \Gamma$.
The {\it left regular representation} of $\Gamma$ is the unitary
representation
$$\pi: \Gamma \longrightarrow \mathcal{L} (\ell^2 (\Gamma ))$$
 defined by
$$(\pi g)f (h)=f(g^{-1}h).$$
Let $ \mathbb{C} \Gamma $ be the group algebra of $\Gamma$.
There is  a unique linear extension of $\pi$ to an (injective) $*$-algebra homomorphism
$$\pi :\mathbb{C} \Gamma \longrightarrow \mathcal{L}(H), \quad \quad \pi (\sum a_g g)=\sum a_g \pi (g). $$
The {\it reduced group $C^*$-algebra} of $\Gamma$, denoted by $C^*_r \Gamma$,
is the norm closure of $\pi (\mathbb{C} \Gamma)$ in $\mathcal{L}(H)$. It
is obviously a unital $C^*$-algebra.
The {\it canonical trace} $\tau$ on $C^*_r \Gamma$ is defined by
$$\tau (a)= \langle a \delta_e, \, \, \delta_e \rangle, $$
 for all $a \in  C^*_r \Gamma$. Notice that on $\mathbb{C} \Gamma$ we have
 $\tau (\sum a_g g) = a_e.$
  It is easily seen that $\tau$ is {\it positive} and {\it faithful} in the sense that
 $\tau (a^* a ) \geq 0$ for all $a$ with equality holding only for
 $a=0$.

 The {\it full}  group $C^*$-algebra of $\Gamma$ is the norm completion
of  $\mathbb{C}\Gamma$ under the norm
$$\|f\|=\mbox{sup}\; \{\|\pi (f)\|; \pi \, \mbox{ is a $\ast$-representation of $\mathbb{C}\Gamma$}\},$$
where by a $\ast$-representation we mean a  $\ast$-representation on a
Hilbert space. Note that $\|f\|$ is finite since for $f=\sum_{g \in
\Gamma}a_g g$ (finite sum) and any $*$-representation $\pi$ we have
$$\| \pi (f)\|\leq \sum \| \pi (a_gg)\| \leq  \sum |a_g| \|\pi (g)\| \leq \sum |a_g|.$$

By its very definition it is clear that there is a 1-1
correspondence between unitary representations of $\Gamma$ and $C^*$
representations of $C^* \Gamma$.
Since the identity map $\text{id}: (\mathbb{C}\Gamma, \| \,\|) \to (\mathbb{C}\Gamma, \| \,\|_r)$ is continuous,
we obtain  a surjective $C^*$-algebra homomorphism
$$C^* \Gamma \longrightarrow C^*_r \Gamma.$$
It is known that this map is an isomorphism if and only if $\Gamma$ is
an amenable group \cite{bl2, da, fi}. }
\end{example}

\begin{example} \rm{  By Fourier transform, or the Gelfand-Naimark theorem, we have an algebra isomorphism
$$C^*_r \mathbb{Z}^n \simeq C
(\mathbb{T}^n).$$ Under this isomorphism, the canonical trace
$\tau$ is identified with the Haar measure on the torus
$\mathbb{T}^n$. More generally, for any abelian group  $\Gamma$ let
$\widehat{\Gamma}=\text{Hom}\,(\Gamma, \, \mathbb{T})$ be the
group of unitary  characters  of $\Gamma$. It is a  compact
group which is  in fact   homeomorphic to  the space of
characters  of the commutative $C^*$-algebra  $C^*_r \Gamma$. Thus
the Gelfand transform defines an algebra isomorphism
\begin{equation}\label{dg}
 C^*_r\Gamma \simeq C (\widehat{\Gamma})
 \end{equation}
 Again the canonical trace $\tau$ on the left hand side is identified with the Haar measure on  $C (\widehat{\Gamma})$. }
\end{example}

In general one should think of the  group $C^*$-algebra of a group
$\Gamma$  as the ``algebra of coordinates" on the noncommutative
space representing the unitary dual of $\Gamma$. Note that by the
above example this is fully justified
  in the commutative case.
 In the noncommutative case, the unitary dual is a badly
behaved space in general but the noncommutative dual is a perfectly legitimate noncommutative space (see
 the unitary dual of the infinite dihedral group in \cite{cob, coma2} and its noncommutative replacement).

\begin{example} \rm{
For a finite group  $\Gamma$  the group $C^*$-algebra coincides with the group
algebra of $\Gamma$. From basic representation theory we know that
the group algebra $\mathbb{C}\Gamma$ decomposes as a sum of matrix algebras
$$ C^* \Gamma \simeq \mathbb{C} \Gamma\simeq \oplus M_{n_i}(\mathbb{C}),$$
where the summation is over the set of conjugacy classes of $\Gamma$.
}
\end{example}

\begin{example} (A noncommutative `connected' space)  \rm{ A {\it projection} in an $*$-algebra is a selfadjoint
idempotent, i.e. an element $e$ satisfying
$$e^2=e=e^*.$$
It is clear that a compact  space $X$ is connected if and only if
$C(X)$ has no non-trivial projections. Let us agree to call a
noncommutative space represented by a $C^*$-algebra $A$
 `connected' if $A$ has no non-trivial projections. The
{\it Kadison conjecture} states that the reduced group
$C^*$-algebra of a torsion-free discrete group is connected. This
conjecture is still open although it has now been verified for
various classes of groups \cite{val}. The validity of the conjecture is
known to follow from the surjectivity of the {\it Baum-Connes assembly
map}, which is an equivariant index map
$$ \mu: K_*^{\Gamma} (\underline{ E \Gamma}) \longrightarrow K_* (C^*_r
(\Gamma).$$ This means that if there are enough `elliptic operators' on
the  classifying space for proper actions of $\Gamma$ then one can
prove an integrality theorem for values of $\tau (e)$, which then
immediately implies the conjecture ($e$ is a projection and $\tau$
is the canonical trace on $C^*_r (\Gamma)$). This principle is
best described in the example below, due to Connes,  where the
validity of the conjecture  for  free groups $F_n$  is established
\cite{co85}.
 Notice that the
 conjecture is obviously true for finitely
 generated  torsion-free abelian groups $\mathbb{Z}^n$
 since  by Pontryagin duality $C^* (\mathbb{Z}^n)\simeq C(\mathbb{T}^n)
 $ and $\mathbb{T}^n$ is clearly connected.

Let $\tau: C^*_r(F_2) \to \mathbb{C}$ be the canonical trace.
Since $\tau $ is positive and faithful, if we can show that, for a
projection $e$,  $\tau (e) $  is an integer then we   can deduce
that $e=0$ or $e=1$. The proof of the {\it integrality}
 of $\tau (e)$ is remarkably similar to integrality theorems for
 characteristic numbers in topology
proved   through index theory. In fact we will show that there is
a Fredholm operator $F_e^+$ with the property that
$$ \tau (e) =\, \text{index}\,  (F_e^+),$$
which clearly implies the integrality of $\tau (e)$. For $p \in
[1, , \infty)$, let $\mathcal{L}^p (H) \subset \mathcal{L} (H)$
denote the {\it Schatten ideal} of $p$-summable compact operators
on $H$.

The proper context for noncommutative index theory is the following \cite{co85}:
\begin{definition} \label{fm} A   $p$-summable Fredholm module over an algebra
$\mathcal{A} $ is a pair $(H, \, F)$ where \\
1. $H=H^ +\oplus H^{-}$ is a $\mathbb{Z}/2$-graded Hilbert space with grading operator $\varepsilon$, \\
2. $H$ is a  left even $A$-module,\\
3. $F \in \mathcal{L}(H)$ is an odd operator with $F^2=I$ and for all $a
\in \mathcal{A}$ one has
\begin{equation} \label{psabl}
[F, \,  a]=Fa -aF \in \mathcal{L}^p(H).
\end{equation}
We say that $(H, \, F)$ is a Fredholm module over $\mathcal{A}$ if
instead of \eqref{psabl}  we have:
\begin{equation} \label{cres1}
[F, \, \, a] \in \mathcal{K} (H)
\end{equation}
for all $a \in \mathcal{A}$.
\end{definition}

The $p$-summability condition \eqref{psabl} singles out `smooth
subalgebras' of  $\mathcal{A}$: the higher the summability order $p$ is the smoother
 $a \in \mathcal{A}$ is. This principle is easily corroborated  in the commutative case.
 The smoother a function is the more rapidly
 decreasing its Fourier coefficients are. That is why if
 $\mathcal{A}$ is a $C^*$-algebra, the natural condition to consider
 is the `compact resolvent' condition \eqref{cres1}  instead of \eqref{psabl}. In general any
 Fredholm module over a $C^*$-algebra $A$ defines  a series of subalgebras of `smooth functions' in
 $A$ which are
 closed under holomorphic functional calculus and have the same $K$-theory as $A$
 (cf. Section 3.1 for an explanation of these terms).

 For
  $A =C ^{\infty}(M)$, $M$ a closed $n$-dimensional  smooth manifold,
   $p$-summable Fredholm  modules for $p > n$ are defined using
  elliptic operators $D$ acting  between sections of vector bundles $E^+$ and $E^{-}$ on  $M$.
  Then one lets $ F =\frac{D}{|D|}$ be the phase of $D$ (assuming
  $D$ is injective), and $H^{+}$ and $H^{-}$  the Hilbert spaces of
  square integrable sections of $E^+$ and $E^{-}$, respectively.  The algebra of continuous (respec. smooth)
  functions on $M$ acts by multiplication on sections and the resulting pair $ (H, \, F)$ is a Fredholm module
  (resp. $p$-summable Fredholm module for $p >n$). }
 \end{example}

Now let $(H, \, F)$ be a 1-summable Fredholm module over an algebra
$\mathcal{A}$. Its {\it character} is the linear functional $\text{Ch} \,(H, F):
\mathcal{A} \to \mathbb{C}$ defined by
$$ \text{Ch} \,(H, F) (a)=\frac{1}{2} \text{Tr}\, (\varepsilon F [F, \, a]),$$
which is finite for all $a \in \mathcal{A}$ thanks to the
1-summability condition on $(H, F)$. As a good exercise the reader
should check that $ \text{Ch} \,(H, F)$ is in fact a trace on $\A$.

The second ingredient
that we need is  an index formula. Let $e\in \mathcal{A}$ be an
idempotent.  Let $F= \left[ \begin{array}{cc} 0 & Q\\P& 0
\end{array} \right]$,  $H_1 =eH^+, \, \, H_2 =eH^{-1}$ and $P':
H_1 \to H_2$ and let  $Q': H_2 \to H_1$ be the restrictions of $P$ and
$Q$ to $H_1$ and $H_2$, respectively. The trace  condition $[F, \,
e] \in \mathcal{L}^1 (H)$  is equivalent to $P'Q'- 1_{H_2} \in
\mathcal{L}^1 (H_2)$ and  $Q'P'- 1_{H_1} \in \mathcal{L}^1 (H_1), $
which of course imply that the operator
\begin{equation} \label{fop}
F_e^+: =P'
\end{equation}
 is Fredholm
and its Fredholm index is given by
\begin{eqnarray*}
 \text{Index} \, F_e^+ &=& \text{Tr} (I-Q'P')- \text{Tr}
(I-P'Q')=\frac{1}{2} \text{Tr} (\varepsilon F [F, e])\\
&=& \text{Ch} \,(H, F) (e).
\end{eqnarray*}

We see that we are done provided we can construct a 1-summable
Fredholm module over a  dense and closed under holomorphic
functional calculus subalgebra $\A$ of  $A = C^*_r (F_2)$ such
that for any projection $e \in \A$,
$$\tau (e)= \text{Ch} (H, F) (e).$$
It is
known that a group is free if and only if it  acts freely  on a
tree. Let then $T$ be a tree with a free action of $F_2$ and
let $T^0$ and $T^1$ denote the set of vertices and edges of $T$
respectively. Let $H ^+ =\ell^2 (T^0)$, $H^-= \ell^2 (T^1) \oplus
\mathbb{C}$ with orthonormal basis denoted by $\varepsilon_q$. The
action of $F_2$ on $T$ induce an action of $C^*_r (F_2)$ on $H^+$
and $H^-$ respectively. Fixing  a vertex $p \in T^0$ we can define
a one to one correspondence $\varphi: T^0 -{p} \to T^1$  by
sending $q \in   T^0 -{p}$ to the edge containing $q$ and lying
between $p$ and $q$. This defines a unitary operator $ P: H^+ \to
H^-$ by
$$ P(\varepsilon_q) =
\varepsilon_{\varphi(q)} \quad \text{and} \quad P(\varepsilon_p)=(0, 1).  $$
The compact resolvent condition \eqref{cres} is a consequence of {\it
almost equivariance} of $\varphi$ in the sense that for any $q \in
T^0 -{p}$, $\varphi (gq)=g \varphi (q)$ for all but a finite
number of $g$,  which is not difficult to prove.

\begin{exercise} Prove the last statement. To see what is going on start
with $\Gamma =\mathbb{Z}$ and go through all the steps in the proof.
\end{exercise}

Let $\A \subset A$ denote the subalgebra of all $a\in A$ such that
$$ [F, \, a] \in \mathcal{L}^1 (H).$$
Clearly $(H, \, F)$ is a 1-summable Fredholm module over $\A$.
It can be shown that $\A$ is stable  under holomorphic functional calculus in $A$ (see Section 3.1).
  It still remains to be checked that
for all  $a \in \mathcal{A}$  we have
$ \text{Ch} \, (H, F) (a) =\tau (a),$
 which we leave as an exercise.

\begin{example}(crossed product algebras) \rm{
Let $\alpha : G \to \text{Aut} \, (A)$ be an action of a group $G$ by
algebra automorphisms on  an algebra $A$.  The action of $g$ on
$a$ will be denoted by $g(a)$. The (algebraic) {\it crossed product algebra} $A\rtimes G$ is
the algebra generated by the two subalgebras $A$ and $\mathbb{C} G$ subject to
relations
$$ gag^{-1}=g(a),$$
for all $g$ in $G$ and $a \in A$.
Formally we define $A\rtimes G=A\otimes  \mathbb{C} G$ with
product given  by
$$(a\otimes g)(b\otimes h)=ag(b)\otimes gh.$$
One checks that this is an associative product and that the two definitions are in fact  the same.

In many cases $A=C(X)$ or $A=C^{\infty}(X)$ is an  algebra of
 functions on a space $X$ and $G$ acts on $X$
by homeomorphisms or diffeomorphisms. Then there is an induced
action of $G$ on $A$ defined by $(gf)(x)=f(g^{-1}(x))$. One of the
key ideas  of noncommutative geometry is Connes' dictum that in
such situations the crossed product algebra $ C(X) \rtimes G$
replaces the algebra of functions on the quotient space $X  /G$
(see Section 4.1 for more on this).}
\end{example}

\begin{example} \rm{ Let $G = \mathbb{Z}_n$ be the cyclic group of order $n$ acting by translations  on  $X
=\{1, 2, \cdots, n\}$. The crossed product algebra $C(X) \rtimes G$ is the algebra generated by elements $U$ and $V$ subject
to  the relations
$$ U^n=1,  \quad V^n=1, \quad UVU^{-1} =\lambda V,$$
where $\lambda = e^{\frac{2\pi i}{n}}$.
Here $U$ is a generator of $\mathbb{Z}_n$, $V$ is a generator of
$\hat{\mathbb{Z}}_n$ and we  have used the isomorphism $ C(X)\simeq
\mathbb{C} \hat{\mathbb{Z}}_n$.

We have an isomorphism
\begin{equation} \label{nct1}
 C(X) \rtimes \mathbb{Z}_n \simeq M_n (\mathbb{C}).
\end{equation}
To  see this consider the $n \times n$ matrices

\[ u= \left( \begin{array}{ccccc}
1 & 0 &\cdots  & \ &0\\
0 & \bar{\lambda} & 0 &\ &0 \\
0 & 0 & \bar{\lambda}^2 &\  &0\\
\\
0 &  \cdots & \cdots &0 &\bar{\lambda}^{n-1} \end{array} \right), \quad v=
\left( \begin{array}{ccccc}
0& 1 &  &   &0\\
 & 0 & 1&\cdots &0 \\
 &  & 0 & 1 & \\
 &  &  &  &1\\
1&   &  &  & 0 \end{array} \right)\] They clearly satisfy the
relations
$$ u^n=1,  \quad v^n=1, \quad uvu^{-1} =\lambda v.$$
Moreover one checks that the matrices $u^p v^q$ for $ 1\leq p, q \leq n$
are linearly independent, which shows that the algebra
generated by
 $u$ and $v$ in $M_n(\mathbb{C})$ is all of  $M_n(\mathbb{C})$.
 The isomorphism \eqref{nct1} is now  defined by
 $$ U \mapsto u, \quad V \mapsto v.$$}
 \end{example}

This isomorphism  can be easily generalized:
\begin{exercise} Let $\mathbb{Z}_n$ act on  an algebra $A$. Show that the dual group ${\hat{\mathbb{Z}}}_n$ acts on
$A \rtimes  \mathbb{Z}_n$ and that we have an isomorphism
$$ A \rtimes \mathbb{Z}_n \rtimes {\hat{\mathbb{Z}}}_n \simeq M_n (A).$$
Extend this to actions of finite abelian groups.
\end{exercise}
The above isomorphism is  a discrete analogue of Takai duality for
$C^*$-crossed products of actions of locally compact abelian
groups (cf. \cite{bl2})
$$ A \rtimes G \rtimes \hat{G} \simeq A \otimes \mathcal{K} (L^2 (G)),$$
where $\mathcal{K}$ is the algebra of compact operators.

The above purely algebraic theory can be extended to a
$C^*$-algebraic context.  For any locally compact topological
group $G$ acting by $C^*$-automorphisms on a $C^*$-algebra $A$,
one defines the reduced $ A \rtimes_r G$ and the full $A \rtimes
G$ crossed product $C^*$-algebras. Assuming $G=\Gamma$ is a discrete
group, to define the crossed product let $\pi : A \to \mathcal{L}
(H)$ be a faithful representation of $A$ (the definition turns out
to be independent of the choice of $\pi$). Consider the Hilbert
space $\ell^2 (\Gamma, \, H)$ of square summable functions $\xi :
\Gamma \to H$ and define a representation $\rho : A \rtimes \Gamma \to \mathcal{L}
(H)$ by
$$ (\rho (x) \xi) (g) =\sum_h \pi (g^{-1} (x(h)))\xi (h^{-1}g)$$
for all $x \in A \rtimes \Gamma$, $\xi \in H$, and $ g \in \Gamma$. It
is an injective $*$-algebra homomorphism and   the reduced
$C^*$-crossed product algebra $ A \rtimes_r \Gamma$ is defined to be the norm completion of the image of $\rho$.

\begin{example} \rm{The crossed product  $C^*$-algebra
$$ A_{\theta} = C(S^1) \rtimes \mathbb{Z}$$
is known as the {\it noncommutative torus}. Here $\mathbb{Z}$ acts on
$S^1$ by rotation through the angle $ 2 \pi \theta$.}
\end{example}

\subsection{Continuous fields of $C^*$-algebras}
 Let $X$ be a
locally compact Hausdorff
topological space and for each $x\in X$ let a $C^*$-algebra $A_x$ be
given. Let also
 $\Gamma \subset \, \coprod A_x$  be  a $*$-subalgebra.
We say this data defines    a {\it continuous field } of
$C^*$-algebras over $X$
if \\
1) for all $x \in X$, the set $\{ s(x); \, s \in \Gamma \}$ is dense in $A_x$,\\
2) for all $s\in \Gamma$, $x \to  \| s(x)\|$ is continuous on $X$,\\
3) $\Gamma $ is locally uniformly closed, i.e.  for any
section $s \in \coprod A_x$,  if for any $x \in X$  we can approximate
$s$  around $x$ arbitrarily closely by  elements of $\Gamma$, then $s\in \Gamma$.  \\

To any continuous field of $C^*$-algebras over a locally compact
space $X$ we associate the   $C^*$ algebra of its {\it continuous
sections} vanishing at infinity (a section $s: X \to   \coprod
A_x$ is called continuous if $x \to \| s(x)\|$ is continuous).

\begin{examples} \rm{ i) A rather trivial example is the $C^*$-algebra
$$A= \Gamma (X,\,  \text{End}\,  (E))$$
of continuous  sections of the endomorphism bundle  of
a vector bundle $E$ over a compact space $X$. A nice example of this is
   the noncommutative torus
 $A_{\theta}$ for $\theta = \frac{p}{q}$  a rational number.
 As we shall see
 $$A_{\frac{p}{q}}\simeq \Gamma ( \mathbb{T}^2, \, \text{End}\, (E))$$
 is the
 algebra of continuous
 sections of the endomorphism bundle of a vector
 bundle of rank $q$ over the torus $\mathbb{T}^2$. Notice that in general $C(X) \subset Z (A)$ is a subalgebra
 of the center of $A$. As we shall see  $ Z(A_{\frac{p}{q}}) \simeq
 C(\mathbb{T}^2)$.}

\noindent ii) A
less trivial example is the field over the interval $[0, 1]$, where $A_x
=M_2(\mathbb{C})$ for $0<x<1$ and $A_0 =A_1 =\mathbb{C} \oplus
\mathbb{C}.$

\noindent iii) The noncommutative tori $ A_{\theta}$,  $ 0\leq \theta \leq 1$,  can be put together
to form a continuous field of $C^*$-algebras over the circle (see  \cite{bl2}).
\end{examples}

\begin{exercise} An {\it Azumaya algebra} is the algebra of continuous sections of a (locally trivial) bundle of finite
dimensional full matrix algebras over a  space $X$. Give an example of an
Azumaya algebra which is not of the type $\Gamma (\text{End}\,  (E))$ for
some vector bundle $E$.
\end{exercise}

Of particular interest  are  $C^*$-algebras of continuous sections
of a locally trivial bundle of algebras with fibers the algebra
$\mathcal{K} (H)$ of compact operators and structure group $\text{Aut} \,\mathcal{K} (H)$. When $H$ is infinite
dimensional such bundles can be completely classified in terms of
their {\it Dixmier-Douady invariant} $\delta (E) \in H^3 (X, \,
\mathbb{Z})$ \cite{dido}. (cf.  \cite{rawi} for a modern detailed
account). The main reason for this is that any automorphism of
$\mathcal{K} (H)$ is inner. In fact for any Hilbert space $H$
there is an exact sequence of topological groups
\begin{equation} \label{caut}
1 \longrightarrow U(1) \longrightarrow U (H) \overset{\text{Ad}}{\longrightarrow} \text{Aut} \,\mathcal{K} (H) \longrightarrow 1,
\end{equation}
where the unitary group  $ U (H)$ has its strong operator topology
and the automorphism group $\text{Aut} \,\mathcal{K} (H)$ is taken with
its norm topology. It is rather easy to see that if $H$ is
infinite dimensional then $ U(H)$ is contractible. (A much
harder theorem of Kuiper states that it is also contractible under
the norm topology, but this is not needed here.)

Let $X$ be a locally compact Hausdorff and paracompact  space.
Locally trivial bundles over $X$ with fibers $\mathcal{K} (H)$ and structure group $\text{Aut} (\mathcal{K} (H))$ are
classified by their classifying class in the \v{C}ech cohomology group\\
$H^1 (X, \, \text{Aut} (\mathcal{K} (H)).$ Now, since the middle
term in \eqref{caut} is contractible, it stands to reason to
expect that
$$ H^1 (X, \, \text{Aut} \,\mathcal{K} (H)) \simeq H^2 (X, \, U(1))
\simeq H^3 (X, \, \mathbb{Z}).$$ The first isomorphism is actually
not so  obvious because \eqref{caut} is an exact sequence of
nonabelian  groups, but it is true (cf. \cite{rawi} for a proof).  The resulting class
$$ \delta (E) \in H^3 (X, \, \mathbb{Z})$$
is a complete isomorphism invariant of such bundles.

Continuous fields of $C^*$-algebras are often constructed by
crossed product algebras. Given a crossed product algebra $ A
\rtimes_{\alpha} \mathbb{R}$, by rescaling the action $\alpha_t$ to
$\alpha_{st}$, we obtain a one-parameter family of crossed product
 algebras
$A \rtimes_{\alpha^s} \mathbb{R}$
over $[0, \infty)$ whose  fiber at 0 is $A \otimes C_0(\mathbb{R})$.
For example,  the translation action of $\mathbb{R}$ on itself will give a field of
$C^*$-algebras whose fiber at 0 is $C_0(\mathbb{R}) \otimes C_0
(\mathbb{R}) \simeq \, C_0(\mathbb{R}^2)$ and at $t >0$ is
$C_0 (\mathbb{R})\rtimes_{\alpha} \mathbb{R}\simeq \mathcal{K}(L^2 (\mathbb{R}))$.

\begin{example} \rm{Let $H$ be the discrete {\it Heisenberg group} of upper triangular matrices
$$ \left( \begin{array}{ccc} 1& a &c\\0 & 1 &b\\0& 0 &1
\end{array}\right),$$
with integer entries.
It is a solvable, and hence  amenable,  group which means that  the
reduced and full group $C^*$-algebras of $H$ coincide. The group
$C^*$-algebra $C^* (H)$ can also be defined as the $C^*$-algebra
generated by three unitaries $U, V, W$  subject to relations
\begin{equation} \label{hga}
 UV = W V U, \quad UW =W U, \quad V W=WV.
 \end{equation}
They correspond to elements
$$ u=\left( \begin{array}{ccc} 1& 1&0\\0 & 1 & 0\\0& 0 &1
\end{array}\right), \quad v=\left( \begin{array}{ccc} 1& 0 &0\\0 & 1 &1\\0& 0 &1
\end{array}\right), \quad w=\left( \begin{array}{ccc} 1& 0 &1\\0 & 1 &0\\0& 0 &1
\end{array}\right).$$
Notice that $W$ is central and in fact it generates the center
$ Z(A) \simeq \, C(S^1)$. For each  $ e^{2 \pi i \theta} \in S^1$, we have a surjective $C^*$-morphism
$$  f_{\theta} : C^* (H) \longrightarrow A_{\theta}$$
by sending $ U\mapsto U$, $V \mapsto V$ and $W \mapsto e^{2 \pi i
\theta} 1$. This is clear from \eqref{hga}. Thus, roughly
speaking,  the Heisenberg group $C^*$-algebra $C^* (H)$ can be
viewed as a `noncommutative  bundle'  over the circle with
noncommutative tori $A_{\theta}$ as fibers.   Notice, however,
that this bundle is not a locally trivial bundle of algebras
 as different fibers are not
isomorphic to each other. One can show  in fact  $C^* (H)$ is
the $C^*$-algebra of continuous sections of a field of
$C^*$-algebras over the circle whose fiber at
  $e^{2\pi i \theta}$ is the noncommutative torus
$A_{\theta}$  (cf. \cite{bl2, fi} and references therein). }
\end{example}

\subsection{ Noncommutative tori}
 These algebras can be defined in a variety of ways, e.g. as the $C^*$-algebra of the Kronecker
foliation of the two-torus by lines of constant slope $dy=\theta
dx$,  as the crossed product algebra $C(S^1)\rtimes \mathbb{Z}$
associated to the automorphism of the circle through  rotation by an
angle $2\pi \theta$, as strict deformation quantization, or by
generators and relations as we do  here.

Let   $\theta \in
\mathbb{R}$ and $\lambda =e^{2\pi i \theta}$. The {\it noncommutative torus}  $A_{\theta}$ is the
{\it universal} unital  $C^*$-algebra
generated by unitaries $U$ and $V$ subject to the relation
\begin{equation}\label{nct}
UV=\lambda VU
\end{equation}
The universality property here means  that  given any  unital $C^*$-algebra $B$ with two unitaries $u$ and $v$ satisfying
$uv=\lambda vu$, there exists a unique unital $C^*$ morphism
$A_{\theta} \to B$ sending $U \to u$ and $V \to v$.

Unlike the purely algebraic case where any set of generators and
relations automatically defines a universal algebra, this is not
the case for universal $C^*$-algebras.  Care must be applied in
defining a norm satisfying the $C^*$-identity, and in general the
universal problem does not have a solution (cf. \cite{bl2} for more on this). For the noncommutative
torus we proceed as follows.  Consider the unitary operators $U,
V: \ell^2 (\mathbb{Z}) \to  \ell^2 (\mathbb{Z})$ defined by
$$   (U f)(n)=e^{2 \pi i\theta n}
f(n), \quad \quad  (V f)(n)= f(n-1).$$
They satisfy $UV= \lambda VU$. Let $A_{\theta}$
be the unital $C^*$-subalgebra of $\mathcal{L} (\ell^2
(\mathbb{Z}))$ generated by $U $ and $V$.  To check the
universality condition, it suffices to check that the operators $U^m
V^n$,  $m, n \in \mathbb{Z}$, are linearly independent, which we
leave as an exercise.

Using  the Fourier isomorphism  $ \ell^2 (\mathbb{Z}) \simeq L^2 (S^1)$, we see that $A_{\theta}$ can also be
described as the $C^*$-subalgebra of $ \mathcal{L} (L^2 (S^1))$
generated by the operators $ (Uf) (z) = f (\lambda z)$ and $(V f)(z) =z f(z)$
for all $f \in L^2 (S^1)$. This presentation also makes it clear that we have an
isomorphism
$$ A_{\theta} \simeq C(S^1) \rtimes \mathbb{Z}$$
where $\mathbb{Z}$ acts by rotation though the angle $2\pi \theta$.

Since the Haar measure on $S^1$ is  invariant under rotations, the formula
$$ \tau (\sum f_i V^i) =\int_{S^1} f_0$$
defines a faithful
 positive
normalized trace
$$\tau : A_{\theta} \to \mathbb{C}.$$
This means that $\tau$ is a trace and
$$ \tau (aa^*) >0, \quad  \quad \tau (1)=1,$$
for all $a \neq 0$.
On the dense subalgebra of $A_{\theta}$ generated  by $U$ and $V$ we
have
\begin{equation} \label{trnct}
\tau (\sum a_{mn} U^m V^n) =a_{00}.
\end{equation}

\begin{exercise} Show that if $\theta$ is irrational then there is a unique
trace, given by \eqref{trnct}, on the subalgebra of $A_{\theta}$ generated  by $U$ and $V$. For
rational values of $\theta$ show that there are uncountably many traces.
\end{exercise}

The map $U \mapsto V$ and $ V \mapsto U$ defines an isomorphism $
A_{\theta} \simeq A_{1-\theta}$. This shows that we can restrict
to $\theta \in [0, \frac{1}{2})$. It is known that in this range $A_{\theta_1}$ is isomorphic to
$A_{\theta_2}$ if and only if $\theta_1 =\theta_2$.    Notice that $A_{\theta}$ is commutative if
and only if $\theta$ is an integer, and
  simple Fourier
theory shows that $A_{0}$ is isomorphic to the algebra
$C(\mathbb{T}^2)$ of continuous functions  on the 2-torus.   For
  irrational  $\theta$, $A_{\theta}$ is known
to be a simple $C^*$-algebra, i.e. it has no proper closed two-sided
ideal. In particular it has no finite dimensional representations \cite{bl2, da, fi}.

For  $\theta =\frac{p}{q}$  a rational number, $A_{\theta}$ has a
finite dimensional representation. To see this   consider the unitary $q\times q$ matrices

\[ u= \left( \begin{array}{ccccc}
1 & 0 &\cdots  & \cdots &0\\
0 & \lambda & 0 &\cdots &0 \\
0 & 0 & \lambda^2 &\cdots &0\\
\cdots\\
0 &  \cdots & \cdots &0 &\lambda^{q-1} \end{array} \right), \quad
v=  \left( \begin{array}{ccccc}
0 & 0 &\cdots  & \cdots &1\\
1 & 0 & 0 &\cdots &0 \\
0 & 1& 0 &\cdots &0\\
\cdots\\
0 &  \cdots & \cdots &1 & 0 \end{array} \right) \]
They satisfy the relations
\begin{equation}
\label{nctr}
 uv=\lambda vu, \quad \quad u^q=v^q=1.
 \end{equation}
It follows that there is a unique $C^*$ map
 from $A_{\theta} \to
 M_q(\mathbb{C})$ sending $U \to u$ and $V \to v$. This,   of course,  implies
 that $A_{\theta}$ is not simple for $\theta$ a rational
 number. In fact it is known that in this case $A_{\theta}$ is Morita
 equivalent to the commutative algebra $C(\mathbb{T}^2)$, as we shall see
 shortly.

 \begin{exercise} Show that the $C^*$-algebra generated  by relations
 \eqref{nctr} is isomorphic to $M_q(\mathbb{C}).$ (hint: it suffices to
 show that the  matrices $u^i v^j$,  $ 1\leq i, j \leq q,$ are linearly
 independent.)
 \end{exercise}

Assuming $\theta =\frac{p}{q}$ is rational, we have  $U^q V= VU^q$ and $
V^qU=UV^q$,   which show that
 $U^q$ and $V^q$ are in the center of $A_{\theta}$. In fact it can be shown that they generate the center and we have
$$ Z(A_{\theta})\simeq C(\mathbb{T}^2).$$
One shows that there is a (flat) vector bundle $E$ of rank $q$ over
$\mathbb{T}^2$  such that  $A_{\frac{p}{q}}$ is isomorphic to the
algebra of continuous sections of the endomorphism bundle of $E$,
   $$ A_{\frac{p}{q}} \simeq \Gamma ( \mathbb{T}^2, \, \text{End}\, (E)).$$

\begin{exercise} Show that if $\theta$ is irrational then $ Z
(A_{\theta})= \mathbb{C}1.$
\end{exercise}

There is a  dense $*$-subalgebra $\mathcal{A}_{\theta} \subset A_{\theta}$ that
deserves to be called the algebra of {\it `smooth functions'} on
the noncommutative torus (see Section 3.1 for more on the meaning of the word smooth in noncommutative geometry).
By definition $\A_{\theta}$ consists of elements of the form
$$ \sum_{(m, n)\in \mathbb{Z}^2} a_{mn} U^m V^n$$
where $(a_{mn}) \in  \mathcal{S} (\mathbb{Z}^2)$ is a rapid decay
 Schwartz class function, i.e. for all $ k\geq 1$
$$ \underset{m,n}{\text{Sup}} \, \,  |a_{mn}| (1+|m|+|n|)^k < \infty.$$
\begin{exercise} Show that  $\A_0 \simeq C^{\infty} (\mathbb{T}^2)$.
\end{exercise}

A {\it derivation} or  {\it infinitesimal automorphism}  of  an
algebra $\A$ is a linear map $ \delta : \A \to \A$ satisfying $
\delta (ab) =\delta (a) b +a \delta (b)$ for all $a, b \in \A$.
For the algebra $\A = C^{\infty} (M)$ of smooth functions on a
manifold $M$ one checks that there is a one to one correspondence
between derivations on $\A$ and vector fields on $M$.  For
 this reason derivations are usually considered as noncommutative analogues of vector fields. The
  fundamental   derivations
 of the noncommutative torus,
$\delta_1, \, \delta_2 :  \mathcal{A}_{\theta} \to
\mathcal{A}_{\theta}$ are the derivations  uniquely defined by their values on
generators of the algebra
$$ \delta_1 (U)=U, \quad \delta_1 (V)=0, \quad \text{and} \quad  \delta_2 (U)=0, \quad \delta_2 (V)=V,$$
or equivalently
$$ \delta_1 (\sum_{(m, n)\in \mathbb{Z}^2} a_{mn} U^m V^n) = \sum_{(m,
n)\in \mathbb{Z}^2} ma_{mn} U^m V^n,$$
and similarly for $\delta_2$.

The invariance property of the Haar measure for the torus has a
noncommutative counterpart. It is easy to see that the canonical
trace  $\tau: \A_{\theta} \to \mathbb{C} $  is invariant under
$\delta_1$ and $\delta_2$, i.e. $\tau (\delta_i (a))=0$ for  $i=1,
2$ and all $a \in \mathcal{A}_{\theta}$.

\begin{exercise} Find all derivations of $\A_{\theta}$.
\end{exercise}

\begin{remark} Although they are deformations of Hopf algebras (in fact
deformations of a group), noncommutative tori are not Hopf
algebras. They are however in some sense Hopf algebroids. Since
they are groupoid algebras, this is not surprising.
\end{remark}

The definition of noncommutative tori can be extended to higher dimensions.   Let  $\Theta$ be a real skew-symmetric
$n\times n$ matrix and let $\lambda_{j k}=e^{2\pi i \Theta_{j
k}}$. The {\it noncommutative torus} $A_{\Theta}$ is the universal
unital $C^*$-algebra generated by unitaries $U_1, U_2, \cdots,
U_n$  with
$$U_k U_j =\lambda_{j k} U_j U_k$$
for all $1\leq j, k\leq n.$ Alternatively, $A_{\Theta}$ is the
universal $C^*$-algebra generated by unitaries $U_l$, $l\in
\mathbb{Z}^n$ with
$$U_l U_m =e^{\pi i \langle l, \,\Theta m \rangle} U_{l+m},$$
for all $l, m \in \mathbb{Z}^n$. From this it follows hat
for any $B \in GL(n, \mathbb{Z})$ we have $A_{B\Theta B^t}\simeq
A_{\Theta}$. The definitions of the invariant trace, smooth subalgebras, and fundamental derivations are easily extended to higher
dimensions. The isomorphism types and Morita equivalence classes of higher dimensional noncommutative tori,
as well as finite projective modules over them has been extensively studied in recent years.

\section{Beyond $C^*$-algebras}
To plunge deeper into noncommutative geometry one must employ various  classes of noncommutative algebras that
are not $C^*$-algebras. They include
 dense subalgebras of $C^*$-algebras which  are {\it stable under
holomorphic functional calculus}, e.g.  algebras of smooth functions on
 noncommutative tori. Another class consists of {\it almost commutative algebras}. They  have a
 {\it Poisson algebra} as their
 {\it semiclassical limit} and   include algebras of
differential operators and  enveloping algebras.   They  appear  in questions of
deformation quantization and its applications.

\subsection{Algebras stable  under holomorphic functional calculus}

 We are going to describe a situation where many features of the
 embedding $C^{\infty} (M) \subset C(M)$  of the algebra of smooth functions on a closed manifold
 into the  algebra of continuous functions is captured and extended
 to  the  noncommutative world.

 Let $A$ be a unital Banach algebra and let $f$ be a  holomorphic function defined on a neighborhood
 of $ \text{sp}\,  (a)$, the spectrum of $a \in A$. Let
 \begin{equation} \label{hfc}
 f(a):= \frac{1}{2\pi i}\int_{\gamma} f(z) (z1-a)^{-1}dz,
 \end{equation}
 where the contour $\gamma$ goes around the spectrum of $a$ only once (counter clockwise). The
 integral is independent of the choice of the contour  and
 it can be shown that for a fixed $a$, the map $f \mapsto f(a)$ is a
 unital algebra map from the algebra of holomorphic functions on a
 neighborhood of $\text{sp}\,  (a)$ to $A$. It is called the {\it holomorphic functional calculus}.
 If $f$  happens to be  holomorphic in a disc containing $\text{sp} (a)$ with power series expansion
 $f(z) =\sum c_i z^i$,   then one shows, using the Cauchy
 integral formula, that $ f(a) =\sum c_i a^i$. If $A$ is a $C^*$-algebra
 and $a$ is a normal element then, thanks to the Gelfand-Naimark theorem,  we have the  much more powerful  {\it
 continuous functional  calculus} from $ C( \text{sp}\, (a)) \to A$. It extends the
 holomorphic functional calculus (see below).

 \begin{definition} Let $B\subset A$ be a unital subalgebra of a unital
 Banach algebra $A$. We say $B$ is stable  under holomorphic
 functional calculus if for all $a \in B$ and all holomorphic functions
 on  $\text{sp}\, (a)$, we have $ f(a) \in B.$
 \end{definition}

 \begin{example} \rm{ 1) The algebra $C^{\infty}(M)$ of smooth functions on a closed smooth manifold
  $M$ is stable under holomorphic functional calculus
 in $C(M)$.  The same can be said about the algebra $C^{k}(M)$ of $k$-times differentiable functions.
 The algebra $\mathbb{C} [X]$ of polynomial functions is not
 stable under holomorphic functional calculus in $C[0, 1]. $\\
 2) The smooth noncommutative torus  $\mathcal{A}_{\theta} \subset
 A_{\theta}$  is stable under holomorphic functional calculus.}
 \end{example}

 Let $\text{sp}_B \, (a)$
 denote the spectrum of $a \in B$ with respect to the subalgebra $B$.
 Clearly $\text{sp}_A \, (a)\subset \text{sp}_B \, (a)$ but the reverse inclusion holds if
 and only if invertibility in $A$ implies invertibility in $B$. A good
 example to keep in mind is $ \mathbb{C} [x] \subset C [0, 1]$.
  It is
  easy to see that if $B$ is stable  under holomorphic functional
  calculus, then we have the {\it spectral  permanence property}
 $$\text{sp}_A \, (a)\subset \text{sp}_B \, (a).$$
 Conversely, under some  conditions on the subalgebra $B$ the above  spectral  permanence property implies that
 $B$ is stable under holomorphic functional calculus.  In fact,  in this case for all $ z \in
 \text{sp}_A \, (a), \, \, (z 1 -a)^{-1} \in B$ and if there is a suitable
 topology in $B$, stronger than the topology induced from $A$,   in which $B$ is
 complete,
  one can then show that the integral \eqref{hfc} converges in $B$.  We give two instances where this technique works.

Let $(H, F)$ be a  Fredholm module over a Banach algebra $A$ and
assume that the action of $A$ on $H$ is continuous.

\begin{prop} (\cite{cob}) For each $p\in [1, \infty )$, the subalgebra
$$ \mathcal{A}= \{ a \in A ; \quad [F, \, a] \in \mathcal{L}^p (H) \}$$
 is  stable under holomorphic
functional calculus.
\end{prop}

 \begin{exercise} (smooth compact operators). Let $\mathcal{K}^{\infty}\subset
 \mathcal{K}$ be the algebra of infinite matrices $(a_{i j})$  with rapid decay
 coefficients. Show  that $\mathcal{K}^{\infty}$ is stable under holomorphic functional calculus  in the algebra of
 compact operators $\mathcal{K}$.
 \end{exercise}

 Another source of examples  are {\it smooth vectors}  of Lie group actions.
 Let $G$ be a Lie group acting continuously on a Banach
 algebra $A$.  An element $a \in A$ is called {\it smooth} if the map $G \to A$ sending $g \mapsto g(a)$ is smooth.
 It can be shown that smooth vectors  form a dense subalgebra of $A$
 which is stable under   holomorphic functional calculus.

 \begin{example} \rm{ The formulas
 $$ U \mapsto \lambda_1 U, \quad \quad V \mapsto \lambda_2 V,$$
 where $(\lambda_1, \lambda_2) \in \mathbb{T}^2$,  define an action of the two-torus $\mathbb{T}^2$ on the noncommutative torus
 $A_{\theta}$.  Its set  of smooth vectors  can be shown to coincide
 with the smooth noncommutative torus $\mathcal{A}_{\theta}$
 \cite{co78}}.
 \end{example}

 For applications  to $K$-theory and density theorems, the following result is crucial \cite{sc}.
 \begin{prop} If  $B$ is  a dense subalgebra of a Banach algebra $A$
 which is stable  under holomorphic functional calculus then so is $M_n(B)$ in $M_n(A)$ for all $n \geq 0$.
 \end{prop}
Now let $e \in A $ be an idempotent in $A$. For any $\epsilon >0$  there is an idempotent $e'\in B$ such that
 $ \|e-e'\|< \epsilon. $
 In fact,  since $B$ is dense in $A$ we can first approximate it by an element
 $g \in B$. Since $\text{sp}\, (e) \subset \{0, 1\}$,  $\text{sp}\, (g)$
 is concentrated around 0 and 1. Let $f$ be a holomorphic function which
 is identically equal to 1 around 0 and 1. Then
 $$e'=f(g)= \frac{1}{2\pi i}\int_{\gamma} f(z) (z1-e)^{-1}dz,$$
 is an idempotent in  $B$ which is close to $e$. In particular $[e]=[f(g)]$ in $K_0(A)$ (see Section 5.1). Thanks to the above
 Proposition, we
 can repeat this argument for  $M_n(B)\subset M_n(A)$ for all $n$. It
 follows that if $B$ is dense is $A$ and is stable  under holomorphic functional calculus, the natural embedding $B \to A$ induces an isomorphism
 $K_0(B) \simeq K_0(A) $ in
 $K$-theory (cf. also the article of J. B. Bost \cite{bo} where a more general density theorem along these lines is proved).  

\begin{example} (Toeplitz algebras) \rm{The original Toeplitz algebra $\mathcal{T}$ is defined as the universal
unital $C^*$-algebra generated by an {\it isometry}, i.e. an element $S$ with
$$ S^* S =I.$$
It can be concretely realized as the $C^*$-subalgebra of $ \mathcal{L}
(\ell^2 (\mathbb{N}))$ generated by the unilateral forward  shift operator $ S
(e_i)=e_{i+1}, \, i=0, 1, \cdots$. Since the algebra $C(S^1)$ of
continuous functions on the circle is the universal algebra defined by a unitary $u$, the map $S \mapsto u$ defines a $C^*$-algebra
surjection $$ \sigma: \mathcal{T} \longrightarrow C(S^1),$$
 called the {\it symbol map}. It is an example
of the symbol map for pseudodifferential operators of order zero over a closed manifold
(see below).

The rank one
projection $I-SS^*$ is in the kernel of $\sigma$. Since the closed ideal generated by $I-SS^*$ is the ideal $\mathcal{K}$ of
compact operators, we have $ \mathcal{K} \subset \text{Ker}\, \sigma$. With
some more work one shows that in fact $ \mathcal{K} = \text{Ker}\, \sigma$
and therefore we have a short exact sequence of $C^*$-algebras, called the  {\it Toeplitz extension}
(due to Coburn, cf. \cite{dou, fi})
\begin{equation} \label{text}
 0\longrightarrow \mathcal{K} \longrightarrow \mathcal{T}
\overset{\sigma}{\longrightarrow}
C(S^1)\longrightarrow 0.
\end{equation}

There is an alternative description of the Toeplitz algebra and
extension \eqref{text} that makes its relation with pseudodifferential
operators and index theory more transparent.
Let $H= L^2(S^1)^+$ denote the Hilbert
space of square integrable functions on the circle whose negative
Fourier coefficients vanish and let $P: L^2(S^1) \to L^2(S^1)^+$
denote the canonical projection. Any  continuous function
$f\in C(S^1)$  defines a {\it Toeplitz operator}
$$T_f: L^2(S^1)^+ \to
L^2(S^1)^+, \quad \quad T_f(g)=P(gf).$$
It can be shown that the $C^*$-algebra generated by the set of Toeplitz operators   $\{ T_f; \, \, f
\in C(S^1)\}$ is isomorphic to the Toeplitz algebra $\mathcal{T}$. The relation
$$ T_f T_g -T_{fg} \in \mathcal{K}(H)$$
shows that any Toeplitz operator $T$ can be written as
$$ T= T_f +K$$
where $K$ is a compact operator. In fact this decomposition is
unique and gives another definition of the symbol map $\sigma$ by
$\sigma (T_f +K)= f$. It is also clear from extension \eqref{text}
that a Toeplitz operator $T$ is Fredholm if and only if its symbol
$\sigma (T)$ is an invertible function on $S^1$.

The algebra generated by Toeplitz operators $T_f$ for $ f \in C^{\infty} (S^1)$ is called  the
{\it smooth Toeplitz algebra} $\mathcal{T}^{\infty} \subset \mathcal{T}$. Similar to \eqref{text}  we have  an extension
\begin{equation} \label{stext}
 0\longrightarrow \mathcal{K}^{\infty} \longrightarrow \mathcal{T}^{\infty}
\overset{\sigma}{\longrightarrow}
C^{\infty} (S^1)\longrightarrow 0.
\end{equation}

\begin{exercise} Show that $\mathcal{T}^{\infty}$ is stable under
holomorphic functional calculus in $\mathcal{T}$.
\end{exercise}

\begin{exercise} Show that
$$ \varphi (A,\, B) = \text{Tr} \, ([A, \, B])$$
defines a cyclic 1-cocycle on $\mathcal{T}^{\infty}$. If $f$ is a smooth
non-vanishing function on the circle, show that
$$\text{index} \, (T_f)= \varphi \, (T_f, \, T_{f^{-1}}).$$
\end{exercise}

The Toeplitz extension \eqref{text} has a grand generalization. On
any closed smooth manifold $M$, a (scalar) pseudodifferential
operator $D$ of order zero defines a bounded linear map $D : L^2
(M) \to L^2 (M)$ and its principal symbol $\sigma (D)$ is a
continuous function on $S^* (M)$,  the unit cosphere bundle of $M$
(cf. \cite{lami}).  Let $\Psi^0 (M) \subset \mathcal{L} (L^2 (M))$
 denote the $C^*$-algebra generated by all  pseudodifferential
operators of order zero on $M$. We  then have a short exact sequence of $C^*$-algebras
$$0\longrightarrow \mathcal{K}(L^2 (M)) \longrightarrow \Psi^0 (M)
\overset{\sigma}{\longrightarrow}
C (S^* M)\longrightarrow 0.$$
For $M=S^1$, the cosphere bundle splits as the disjoint union of two copies of $S^1$ and the above sequence is the
direct sum of two identical copies, each of which is isomorphic to the Toeplitz extension
\eqref{text}. }
\end{example}

 \subsection{Almost commutative and Poisson algebras}

Let $A$   be a unital complex algebra.  We say $A$ is a {\it filtered algebra}
if it has   an increasing filtration
     $F^i(A)\subset F^{i+1}(A)$, $i=0, 1, 2, \cdots$, with
     $F^0(A)=\mathbb{C} 1$,
$F^i (A) F^j(A)\subset F^{i+j}(A)$ for all $i, j$ and $\cup_i F^i(A)=A$.
Let $F^{-1}(A)=0$.
The {\it associated graded algebra} of a filtered algebra is the graded
algebra
$$\text{Gr}\, (A)=\bigoplus_{i\geq 0}\frac{F^i(A)}{F^{i-1}(A)}.$$

\begin{definition} An almost commutative algebra is  a filtered algebra
 whose  associated graded
algebra $\text{Gr}\, (A)$ is  commutative.
\end{definition}
\noindent Being almost commutative  is equivalent to the commutator condition
$$[F^i (A), \, \, F^j (A)]\subset F^{i+j-1} (A),$$
for all $i, j$. As we shall see Weyl algebras and more generally
algebras of differential operators on a smooth manifold, and
universal enveloping algebras are examples of almost commutative
algebras.

Let $A$ be an almost commutative algebra. The original Lie bracket $[x, y]=xy-yx$  on $A$ induces a Lie bracket
$\{\,\, \}$ on $Gr (A)$ via the
formula
$$\{x+F^i, \, \, y+F^j\}:= [x, \, \, y] + F^{i+j-2}.$$
 \noindent Notice that by the almost commutativity assumption, $[x, \, \, y]$ is in
$F^{i+j-1}(A)$. The induced Lie bracket on $\text{Gr}\,  (A)$ is compatible with its multiplication
in the sense that for all $a \in \text{Gr}\,  (A)$, the map $b \mapsto \{a, \, b\}$ is a derivation.
The algebra $\text{Gr}\,  (A)$ is called the {\it semiclassical limit} of the almost commutative algebra $A$.
It is an example of a Poisson algebra  as we recall next. The quotient
map
$$ \sigma: A \to \text{Gr}\,  (A)$$
 is called the {\it principal symbol map}.

 Any  splitting $
 q: \text{Gr}\,  (A) \to A$ of this map  can be regarded as a `naive quantization map'. Linear
 splittings always exist but they are  hardly interesting. One usually
 demands more. For example one wants  $q$ to be  a Lie algebra map in the sense that
 \begin{equation}\label{dqr}
 q \{a, \, \, b\}=[q (a), \, \, q (b)]
 \end{equation}
 for all $a, b $ in $\text{Gr}\,  (A)$. This is one form of {\it Dirac's  quantization
 rule} going back to \cite{dir}.
 {\it No-go theorems},   e.g. Groenvald-van Hove's (cf. \cite{gust} for
 a discussion and precise statements), state  that,
 under reasonable non-degeneracy
 conditions,   this is almost never possible.
 The remedy is to have  \eqref{dqr} satisfied only in an asymptotic sense. As we shall discuss later in this section, this can be done
 in different ways either  in the context of formal deformation
 quantization \cite{bffls1, ko} or through strict $C^*$-algebraic deformation quantization \cite{rie5}.

The notion of a Poisson algebra captures  the structure of
semiclassical limits.
\begin{definition}
Let $P$ be a commutative algebra. A {\it Poisson structure} on $P$ is a Lie algebra bracket
$(a, b)\mapsto \{a, b\}$ on $A$ such that for any $a \in A$, the map
$b\mapsto \{a, b\} : A \to A$ is
a derivation of $A$. That is,
for all $b, c$ in $A$ we have
$$\{a, bc\}= \{a, b\}c +b\{a,  c\}.$$
\end{definition}
In geometric examples (see below) the vector field defined by this
derivation is called the {\it Hamiltonian vector field} of the {\it Hamiltonian function} $a$.

\begin{definition}
A {\it Poisson algebra} is a pair $(P, \, \{, \, \})$ where $P$ is  a commutative algebra and
$ \{ \, \, , \, \}$ is
a Poisson structure on $P$.
\end{definition}

We saw that the semiclassical limit $P= \text{Gr}\,  (A)$ of  any almost
commutative algebra $A$  is a Poisson algebra. Conversely, given a
Poisson algebra $P$ one may ask if it is the semiclassical limit
of an almost commutative algebra. This is one form of the problem of
quantization of Poisson algebras,  the answer to which  for general Poisson algebras
is negative. We give a few concrete examples of
Poisson algebras.

\begin{example}  \rm{ A {\it Poisson manifold} is a manifold $M$ whose  algebra of smooth
 functions $A=C^{\infty}(M)$ is a Poisson algebra (we should also
 assume that the bracket $\{\, \, , \, \}$ is continuous in the Fr\'{e}chet
 topology of $A$).
 It is not difficult to see that all Poisson
 structures on $A$ are of the form
$$ \{f, \, g\}:= \langle df\wedge dg, \, \, \pi\rangle,$$
where $\pi \in C^{\infty}(\bigwedge^2 (TM))$
is  a smooth 2-vector field on $M$.
This bracket clearly satisfies the Leibniz rule in each variable
and  one  checks that it satisfies the Jacobi identity  if and
only if $[ \pi, \pi ]=0,$ where the {\it Schouten bracket} is defined in local
coordinates  by
$$  [ \pi, \pi ]= \sum_{l=1}^n(\pi_{lj}\frac{\partial \pi_{ik}}{\partial x_l}+ \pi_{li}\frac{\partial \pi_{kj}}{\partial x_l}
+\pi_{lk}\frac{\partial \pi_{ji}}{\partial x_1})=0.$$
 The Poisson bracket  in local coordinates is given by
$$\{f, \, g\}=\sum_{i j}\pi_{ij}\frac{\partial f}{\partial x_i}
\frac{\partial g}{\partial x_j}.$$
Symplectic manifolds are the simplest examples of Poisson manifolds. They correspond to non-degenerate Poisson
structures. Given a symplectic form $\omega$, the associated  Poisson
bracket is given  by
$$ \{ f, \, g \}= \omega (X_f, X_g),$$
where the vector field $X_f$ is  the symplectic dual of $df$.

Let $C^{\infty}_{\text{poly}}(T^* M)$  be the algebra of smooth functions
on $T^*M$ which are polynomial in the cotangent direction. It is a
Poisson algebra under the natural symplectic structure of $T^*M$.
This Poisson algebra  is the semiclassical limit
 of the algebra of differential operators on $M$, as we will see in the
 next example. }
  \end{example}

\begin{example}  (Differential operators on commutative algebras). \rm{ Let $A$
be a commutative unital algebra. Let $\mathcal{D}^0 (A)$ denote the set of {\it differential operators of order
zero} on $A$, i.e. $A$-linear maps from $A\to A$, and for $n\geq 1$, let $\mathcal{D}^n(A)$ be the set of all operators
$D$ in $\text{End}_{\mathbb{C}} (A)$ such that for any $a\in A$,
$[D,  \, \,a] \in \mathcal{D}^{n-1} (A).$
 The set
   $$ \mathcal{D} (A) =\bigcup_{n\geq 0} \mathcal{D}^n (A)$$
   is a subalgebra of $\text{End}_{\mathbb{C}} (A)$, called the {\it algebra of
   differential operators} on $A$. It is an almost commutative algebra under the  filtration given by $\mathcal{D}^n, \, n\geq0.$
 Elements of $\mathcal{D}^n(A)$ are called differential operators of order $n$. For example, a linear map $D: A \to A$ is a differential operator of
order one if and only if it is of the form $D=\delta +a$, where $\delta $ is a derivation on $A$ and $a\in A$.

 For general $A$, the semiclassical limit $\text{Gr}\,  (\mathcal{D} (A))$ and its Poisson structure are not easily identified except
 for coordinate rings of smooth affine varieties or algebras of smooth functions on a manifold.  In this case
  a differential operator $D$ of order $n$   is locally given  by
$$ D=\sum_{|I|\leq k} a_I(x)\partial^I$$
where $I=(i_1, \cdots, i_n)$ is a multi-index and $\partial^I =\partial_{i_1} \partial_{i_2} \cdots \partial_{i_n}$
is a mixed partial derivative.  This expression depends on the local
coordinates but its leading terms of total degree $n$ have
 an invariant meaning provided we replace $\partial_i \mapsto \xi_i \in T^* M$.  For $\xi \in T^*_x M$, let
 $$\sigma_p (D) (x, \xi): =\sum_{|I|= k} a_I(x)\xi^I$$
 Then the function $\sigma_p(D): T^*M\to \mathbb{C}$, called the {\it principal
 symbol} of $D$,  is invariantly defined and belongs to $C^{\infty}_{\text{poly}}(T^* M)$. The algebra
 $C^{\infty}_{\text{poly}}(T^* M)$ inherits a canonical Poisson structure as a subalgebra of the Poisson algebra
 $C^{\infty}(T^* M)$ and we have the following
\begin{prop} The principal symbol map induces an isomorphism of Poisson
algebras
$$ \sigma_p: Gr \, \mathcal{D}(C^{\infty}(M))\overset{\simeq}\to  C^{\infty}_{poly}(T^* M).$$
\end{prop}}
\end{example}
See \cite{chgi} for a proof of this or, even better, try to prove it yourself by proving it for Weyl algebras first.

\begin{example} (The Weyl algebra) \rm{ Let $ A_1: = \mathcal{D} \,  \mathbb{C}[X]$ be the {\it Weyl algebra}
  of differential operators on the line.
 Alternatively,  $A_1$  can be described as the
 unital complex algebra defined  by generators $x$ and $p$ with
$$px-xp=1.$$
The map $x\mapsto x, \, p \mapsto \frac{d}{dx}$ defines  the
isomorphism. Physicists   prefer to write the defining relation  as the {\it
canonical commutation relation} $pq-qp= \frac{h}{2\pi i}1$, where
$h$ is Planck's constant and $p$ and $q$ represent momentum and
position  operators. This is not without merit because we can then
let $h\to 0$  and obtain the commutative algebra of polynomials in
$p$ and $q$ as the  semiclassical limit. Also, $i$ is necessary if we want to consider $p$ and $q$ as
selfadjoint operators (why?).

Any element of $A_1$ has
a unique expression as a differential operator with polynomial
coefficients $ \sum a_i(x)\frac{d^i}{dx^i}$
 where the standard filtration is by degree of the differential
 operator. We have an algebra isomorphism
 $Gr (A_1)\simeq  \mathbb{C} [x, y]$ under which the (principal) symbol map is  given by
$$\sigma (\sum_{i=0}^n a_i(x)\frac{d^i}{dx^i})= a_n(x)y^n.$$
The induced Poisson bracket on $\mathbb{C} [x, y]$ is the classical Poisson
bracket
$$\{f, \, g\} =\frac{\partial f}{\partial x} \frac{\partial g}{\partial y}-\frac{\partial f}{\partial y}
\frac{\partial g}{\partial x}.$$

In general, the Weyl algebra $A_n$ is the algebra of differential operators
on $\mathbb{C}[x_1, \cdots, x_n]$. Alternatively it can be defined as
the universal algebra defined by $2n$ generators $x_1, \cdots x_n, p_1, \cdots,
p_n$  with
$$ [p_i, x_i]=\delta_{ij}, \quad \text{and} \quad [p_i, p_j]= [x_i,
x_j]=0$$
for all $i, j$. Notice  that $ A_n \simeq A_1 \otimes \cdots \otimes A_1$ ($n$ factors).   }
\end{example}

A lot is known about Weyl algebras and a lot remains to be known,
including the Dixmier conjecture about the automorphisms of $A_n$
and its relation with the Jacobian conjecture recently studied  by
Kontsevich in \cite{beko}.  The Hochschild and cyclic cohomology of $A_n$
is  computed in \cite{ft}.

\begin{exercise}. Show that $A_1$  is a simple algebra, i.e. it has no
non-trivial two-sided ideals; that
any derivation of
$A_1$ is inner; and that
  $[A_1, A_1]=A_1$, i.e. there  are no non-trivial traces on
$A_1$. Prove all of this for $A_n$. Any derivation of $A_1$ is inner; is it true that any automorphism of $A_1$ is inner?
\end{exercise}

\begin{exercise} Let $A=\mathbb{C}[x]/(x^2)$ be the algebra of dual
numbers. Describe its algebra of differential operators.
\end{exercise}

\begin{example} (Universal enveloping algebras) \rm{ Let $U (\mathfrak{g})$ denote the
{\it  enveloping algebra} of a Lie algebra $\mathfrak{g}$. By definition,
$U(\mathfrak{g})$ is the quotient of the tensor algebra $T(\mathfrak{g})$ by the two-sided
ideal generated by $x\otimes y-y\otimes x -[x, y]$ for all $x, y \in \mathfrak{g}$.
 For $p\geq\, 0$, let   $F^p (U (\mathfrak{g}))$ be the subspace generated by tensors of degree at most $p$.  This turns $U(\mathfrak{g})$ into
 a filtered algebra and the
Poincar\'{e}-Birkhoff-Witt theorem asserts that its  associated
graded algebra is canonically isomorphic to the symmetric
algebra $ S (\mathfrak{g})$. The algebra isomorphism is induced by
the {\it symmetrization map} $ s : S(\mathfrak{g}) \to \text{Gr}\, (U
(\mathfrak{g}))$, defined by
$$ s (X_1 X_2 \cdots X_p)= \frac{1}{p!} \sum_{\sigma \in S_p}
X_{\sigma (1)}\otimes \cdots \otimes X_{\sigma (p)}.$$
Note that $ S(\mathfrak{g})$ is the algebra
of polynomial functions on the dual space $\frak{g}^*$  which is a Poisson manifold under the bracket
$$\{f \, , \,g\} (X)=[ Df (X), \ \, Dg  (X)]$$
for all $f, g \in C^{\infty}(\mathfrak{g}^*)$ and $X\in
\mathfrak{g}^*$. Here we have used the canonical isomorphism
$\mathfrak{g}\simeq \mathfrak{g}^{**}$, to regard the differential
$Df (X) \in \mathfrak{g}^{**}$ as an element of $\mathfrak{g}$.
The induced Poisson structure on polynomial functions  coincides
with the Poisson structure in $\text{Gr}\,  (U (\mathfrak{g}))$. }
\end{example}

\begin{example} (algebra of formal pseudodifferential  operators on the
circle) \rm{ This algebra is obtained by formally inverting the differentiation
operator
 $\partial: =\frac{d}{dx}$ and then completing the resulting algebra.
 A formal  pseudodifferential operator on the circle is an expression
of the form $\sum_{-\infty}^n a_i(x)\partial^i$, where each $a_i(x)$ is a Laurent polynomial.
The multiplication in uniquely
defined by the rules $\partial x -x \partial =1$ and $ \partial
\partial^{-1} = \partial^{-1} \partial =1$.
  We denote the resulting algebra by $\Psi_1$.  The
 {\it Adler-Manin trace} on $\Psi_1$ \cite{ma79}, also called the
{\it noncommutative residue}, is defined by
$$\text{Tr}\,  (\sum_{-\infty}^n a_i(x)\partial^i)= \text{Res} \, (a_{-1} (x); \,\, 0)=\frac{1}{2\pi i}\int_{S^1}
a_{-1}(x).$$
This is a trace on $\Psi_1$. In fact one can show that
$\Psi_1 /[\Psi_1, \, \, \Psi_1]$ is one-dimensional which means
 that any  trace on $\Psi_1$ is a multiple of Tr. Notice that for the Weyl algebra $A_1$ we have $[A_1, \,
 A_1]= A_1.$

Another interesting difference between $\Psi_1$ and $A_1$ is that $\Psi_1$ 
 admits non-inner derivations (see exercise below). The algebra
$\Psi_1$ has a nice generalization to  algebras of
pseudodifferential operators in  higher dimensions. The
appropriate extension of the above trace   is the {\it noncommutative
residue} of    Guillemin and Wodzicki (cf. \cite{wod}. See also
\cite{cob} for relations with the Dixmier trace and its place
in noncommutative Riemannian geometry).
}
\end{example}

\begin{exercise} Unlike the algebra of  differential operators,  $\Psi_1$
admits non-inner derivations. Clearly $log \, \partial \notin \Psi_1$,
but show that
 for any $a\in \Psi_1$, we have $[log \,\partial, \,  \, a]\in \Psi_1$ and therefore the
map
$$ a \mapsto \delta (a): =[log \, \partial, \, \, a]$$
defines a non-inner derivation of $\Psi_1$ \cite{krkh}. The
corresponding Lie algebra 2-cocycle
$$ \varphi (a, b) = Tr (a [log \,
\partial, \, \, b])$$ is the Radul cocycle.
\end{exercise}

\subsection{Deformation theory}

In this section we shall  freely use results about  Hochschild cohomology from
Section 6.5.  In the last Section we saw one way to formalize
the idea of quantization   through the notion of an almost commutative algebra and its
semiclassical limit which is a Poisson algebra. A closely related notion is {\it formal deformation
quantization}, or {\it star  products} \cite{bffls1, ko}.  Let $ (A, \, \{\, , \, \})$ be a Poisson algebra and let
$A[[h]]$ be the algebra of formal power series over $A$. A formal
deformation of $A$ is an associative $\mathbb{C}[[h]]$-linear
multiplication
$$*_h: A[[h]] \otimes A[[h]] \to A[[h]]$$
such that $*_0$ is the original multiplication  and for all
$a, b$ in $A$,
$$\frac{a*_h b-b*_ha}{h} \to \{a, b\}$$
as $h\to 0$.
Writing
$$ a*_h b =B_0 (a, b) +h B_1 (a, b)+h^2 B_2(a, b)+\cdots   $$
where $B_i: A \otimes A \to A$ are Hochschild 2-cochains on $A$ with
values in $A$, we see that the initial conditions on $*_h$ are equivalent to
$$ B_0 (a, b) =ab, \quad \text{and} \quad B_1(a, b)-B_1(b, a)=\{a, \, b \}.$$
The associativity condition on $*_h$  is equivalent to an infinite  system of
 equations involving the cochains $B_i$. They are given by
$$ B_0 \circ B_n +B_1 \circ B_{n-1}+\cdots +B_n\circ B_0=0, \quad \text{for all} \, n\geq 0,$$
or equivalently
\begin{equation}\label{ae}
 \sum_{i=1}^{n-1} B_i \circ B_{n-i} =\delta B_n.
 \end{equation}
Here, the {\it Gerstenhaber $\circ$ product} of 2-cochains $f, \, g: A \otimes A \to
A$ is defined  as the 3-cochain
$$ f\circ g (a, b, c)= f(g(a, b), c) -f(a, g(b, c)).$$
Notice that a  2-cochain  $f$ defines an associative product if
and only if $f\circ f =0$. Also notice that the Hochschild
coboundary $\delta f$   can be written as $\delta f
=-m\circ f - f\circ m$, where $m: A \otimes A \to A$ is the
multiplication of $A$.  These observations lead to the
associativity equations \eqref{ae}.

To solve these equations starting with $B_0=m$, by
antisymmetrizing we can always assume that $B_1$ is antisymmetric
and hence we can assume  $B_1= \frac{1}{2}\{ \, , \, \}$.   Having
found $B_0, B_1, \cdots B_{n-1}$, we can find a $B_n$ satisfying
\eqref{ae}  if and only if  the cocycle $\sum_{i=1}^{n-1} B_i\circ B_{n-i}$ is a
coboundary, i.e. its class in $H^3 (A, A)$ should vanish.
 The upshot is that the third Hochschild cohomology $H^3 (A, A)$ is the {\it space
of obstructions} for the deformation quantization problem.  In
particular if $H^3(A, A)=0$ then any Poisson bracket on $A$ can be
deformed. Notice, however,  that this is only a sufficient condition and by no means is necessary, as will be shown below.

In the most interesting examples, e.g. for  $A=C^{\infty}(M)$, $H^3(A, A)\neq 0$. To see this
consider the differential graded Lie
algebra   $(C (A, \,A),\,  [ \, ,  \,], \, \delta)$ of continuous
Hochschild cochains on $A$, and the differential graded Lie
algebra, with zero differential,
 $ ( \bigwedge
(TM), \, [ \, ,  \,], \, 0)$ of polyvector fields on $M$. The bracket in
the first is the Gerstenhaber bracket and in the second is the
Schouten bracket of polyvector fields.
By a theorem of Connes (see  the resolution in Lemma 44 in \cite{co85}), we know that the {\it antisymmetrization map}
$$ \alpha : ( C^{\infty}(\bigwedge TM), \, 0)  \to  (C(A, \, A), \, \delta)$$
sending a polyvector field $X_1 \wedge \cdots \wedge X_k$  to the functional $\varphi$ defined by
$$\varphi (f^1, \cdots, f^k) = df^1 (X_1)df^2(X_2) \cdots df^k(X_k)$$
is a quasi-isomorphism of differential graded algebras. In
particular, it induces an isomorphism of graded commutative
algebras
$$ H^k(A, A) \simeq C^{\infty}({\bigwedge}^k TM).$$

The map $\alpha$, however, is
not a morphism of Lie algebras and that is where the real
difficulty of deforming a Poisson structure is hidden. The {\it
formality theorem} of M. Kontsevich \cite{ko}  states that  as a
differential graded Lie algebra, $(C(A, A), \delta, [ \, , \,])$
is formal in the sense that it is quasi-isomorphic to its cohomology. Equivalently, it means that one
 can correct the map $\alpha$, by adding an
infinite number of terms, to a morphism of $L_{\infty}$ algebras.
This shows that the original deformation problem of Poisson
structures can be transferred to $( C^{\infty}(\bigwedge TM), \,
0)$
 where it is  unobstructed since the differential in the latter DGL is zero. We give a couple
of simple examples where deformations   can be explicitly
constructed.

\begin{example} \rm{ The simplest non-trivial Poisson manifold is the dual
$\mathfrak{g}^*$ of a finite dimensional Lie algebra $\mathfrak{g}$. Let $U_h(\mathfrak{g})=T(\mathfrak{g})/I$, where
the ideal $I$ is generated by
$$ x\otimes y-y\otimes x-h[x, y],  \, \, x, y \in
\mathfrak{g}.$$ By the Poincar\'{e}-Birkhoff-Witt theorem, the
antisymmetrization map $\alpha_h: S(\mathfrak{g}) \to
U_h (\mathfrak{g})$ is a linear isomorphism. We can define a $*$-product
on $S (\mathfrak{g})$ by
$$f*_h g= \alpha_h^{-1} (\alpha_h (f) \alpha_h (g))= \sum_{n=0}^{\infty}h^n B_n (f, g).$$
With some work one can show that $B_n$ are bidifferential operators and
hence the formula extends to a $*$-product on $C^{\infty} (
\mathfrak{g}^*)$. }
\end{example}

\begin{example} (Moyal-Weyl quantization)\rm{ Consider the  algebra
generated by $x$ and $y$ with  relation $xy-yx=\frac{h}{i} 1$.
Iterated application of the Leibniz rule gives the formula for
the product
$$ f *_h g  =\sum_{n=0}^{\infty} \frac{1}{n!} (\frac{-ih}{2})^n B_n
(f, g),$$ where $B_0(f, g)=fg, \, \, B_1 (f, g)= \{f, \,\, g\}$
is the standard Poison bracket, and for $n \geq 2$,
$$B_n (f, \, g)=(-1)^n \sum_{k=0}^n (-1)^k \left( \begin{array}{c} n\\k \end{array} \right)(\partial_x^k  \partial_y^{n-k}
\, f) (\partial_x^{n-k}  \partial_y^{k} \, g) .$$ Notice that this
formula makes sense  for $ f, \, g \in C^{\infty} (\mathbb{R}^2)$ and defines a
deformation of this algebra with its standard  Poisson structure. This
can be extended to arbitrary constant Poisson structures
$$\{f, \, \, g\}= \sum \pi^{i j}\partial_i f \, \partial_j  g.$$
The Weyl-Moyal $*$ product is then given by
$$f*_h g = \text{exp} ( -i\frac{h}{2}\sum \pi^{ij} \partial_i \wedge \partial_j) (f, \, \,g).$$ }
\end{example}

\begin{remark} \rm{ As we mentioned before,  deformation quantization has its
origins in Dirac's quantization rule  in  quantum mechanics \cite{dir}. The
original idea was to replace the classical observables (functions)
by quantum observables (operators) in such a way that Poisson
brackets of functions  correspond to commutators of operators. It
was, however,  soon realized that a rigid  interpretation of
Dirac's rule   is
impossible and it must be understood only in an asymptotic
sense. In fact there are some well known  `no go-theorems' that
state that under reasonable non-degeneracy  assumptions this rigid
notion of quantization is not possible (cf. \cite{gust} and
references therein).
 One approach, as in  \cite{bffls1},     formulates  quantization of a classical system  as formal
deformation quantization of  a  Poisson structure manifold.
 We should mention
that the algebraic underpinnings of deformation theory of
(associative and Lie) algebras and the relevance of Hochschild
cohomology  goes back to
Gerstenhaber's  papers \cite{gesc}.  }
\end{remark}

No discussion of deformation quantizations is complete without
discussing Rieffel's deformation quantization \cite{ rie5}. Roughly
speaking,  one demands that formal power series of formal
deformation theory should actually be convergent. More precisely,
let $(M, \{ \, , \, \})$ be a Poisson manifold. A {\it strict
deformation quantization} of the Poisson algebra
$\mathcal{A}=C^{\infty}(M)$ is a family of pre $C^*$-algebra
structures $(*_h, \, \, \| \, \, \|_h)$  on $\mathcal{A}$ for
$h\geq 0$ such that the family forms a continuous field of
$C^*$-algebras on $[0, \, \, \infty)$ (in particular for any $f
\in \mathcal{A}$, \, $h \mapsto \|f\|_h$ is continuous) and for
all $f, g \in \mathcal{A}$,
$$\| \frac{f*_h g- g*_h f}{ih} \|_h\to \{f, \, \, g\}$$
as $h\to 0$. We therefore have a family of $C^*$-algebras $A_h$
obtained by completing $\mathcal{A}$ with respect to the norm $\|
\,\|_h$.

\begin{example} (noncommutative tori) \rm{ In \cite{rie2} it is shown that
the family of noncommutative tori $A_{\theta}$  form a strict
deformation quantization of the Poisson algebra of smooth
functions on the 2-torus. This in fact appears as a special case
of a more general result. Let $\alpha$ be a smooth action of
$\mathbb{R}^n$ on $\mathcal{A}=C^{\infty}(M)$.  Let $X_i$ denote
the infinitesimal generators for this action. Each skew-symmetric
$n \times n$ matrix $J$ defines a Poisson bracket on $\mathcal{A}$
by
$$ \{f, \, g\}=\sum J_{ij} X_i (f) X_j (g).$$
  For each $h\in \mathbb{R}$, define a new
product $*_h$ on $\mathcal{A}$ by
$$ f *_h g= \int_{\mathbb{R}^n \times \mathbb{R}^n} \alpha_{h J u} (f)
\alpha_v (g)e^{2\pi i u\cdot v} du dv.$$
The $*$ structure is by conjugation and is undeformed  (see  \cite{rie2} for the definition of $\| f\|_h$). For
$ \mathcal{A} =C^{\infty} (\mathbb{T}^2)$  with the natural $\mathbb{R}^2$
action one obtains $A_{\theta}$.  }
\end{example}

\begin{remark}. \rm{Does any Poisson  manifold admit a strict deformation quantization?  This question is still open
(even for symplectic manifolds).  In \cite{rie3},  Rieffel shows that
the canonical symplectic structure on the 2-sphere  admits no
$SO(3)$-invariant strict deformation quantization. An intriguing
idea promoted by Connes and Marcolli  in \cite{coma2} is that existence
of a strict deformation  quantization of a Poisson manifold should
be regarded as an {\it integrability condition} for the formal
deformation quantization. There is a clear analogy with the case of
formal and convergent power series solutions of differential
equations around singular points. They  ask for a possible `theory of ambiguity', i.e. a
cohomology theory  that could capture the difference
between the two cases. }
\end{remark}

\section{Sources of noncommutative spaces}
At present we can identify at least four  methods by which noncommutative spaces
are constructed:\\
 i) noncommutative quotients; \\
 ii) algebraic and $C^*$-algebraic deformations;\\
 iii)  Hopf algebras and quantum groups;\\
 iv) cohomological constructions.\\
It should be stressed  that these are by no means mutually
exclusive; there  are intimate relations between these sources and
 sometimes a  noncommutative space can be described by several
methods, as is the case  with  noncommutative tori.  The
majority of examples, by far,  fall into the first category.
 We won't discuss the last idea, advanced by
Connes and Landi \cite{cola} and Connes and  Dubois-Violette
\cite{codu}. Very briefly the idea is that if one writes the
conditions for the Chern character of an idempotent in cyclic
homology to be trivial on the level of chains, then  one obtains
interesting examples of algebras such as noncommutative spheres
and spherical manifolds,  Grassmannians, and Yang-Mills algebras.

\subsection{Noncommutative quotients}

The quotient space $X / \sim $ of a  Hausdorff space may
easily fail to be Hausdorff and may even be an indiscrete
topological space with only two open sets. This happens, for
example, when at least one equivalence class is dense in $X$.
Similarly the quotient of a smooth manifold may become singular
and fail to be smooth.

 The method of {\it noncommutative quotients}  as advanced by Connes in \cite{cob}
 allows one to replace ``bad quotients'' by nice noncommutative
spaces, represented by noncommutative algebras.  In general these noncommutative algebras are defined
as {\it groupoid algebras}. In some cases, like  quotients by group actions,
the noncommutative quotient can be defined as a crossed product
algebra too,  but in general the use of {\it groupoids} seem to be unavoidable.

 An   equivalence relation  is usually obtained from
a much richer structure by forgetting part of this structure. For
example, it  may arise from an action of a group $G$ on $X$
where $x\sim y$ if and only if $gx=y$ for some $g$ in $G$ ({\it
orbit equivalence}). Note that there may be, in general, many $g$
with this property. That is  $x$ may be identifiable with $y$ in
more than one way. Of course when we form the equivalence relation
this extra information is lost. The key idea in dealing with bad
quotients in Connes' theory  is to keep track of this extra
information, by first forming a groupoid.

 Now Connes'  dictum in forming
noncommutative quotients can be summarized as follows:
$$ \mbox{\bf quotient data} \rightsquigarrow \mbox{\bf groupoid} \rightsquigarrow \mbox{\bf groupoid
algebra}, $$
where the noncommutative quotient is defined to be the groupoid algebra itself.

\begin{definition}
A groupoid is a small category in which every morphism is an
isomorphism.
\end{definition}

 The set of objects of a groupoid
$\mathcal{G}$  shall be denoted by $\mathcal{G}^{(0)}$.   Every
morphism has a {\it source},  {\it target} and  an {\it inverse}.
They define maps
$$ s:\mathcal{G} \longrightarrow \mathcal{G}^{(0)}, \quad
t:\mathcal{G} \longrightarrow \mathcal{G}^{(0)}, \quad
i: \mathcal{G} \longrightarrow \mathcal{G}.$$
Composition $\gamma_1 \circ \gamma_2$ of morphisms $\gamma_1 $ and $
\gamma_2$ is only defined if $ s(\gamma_1)=t(\gamma_2)$. Composition
defines a map
$$ \circ:  \mathcal{G}^{(2)}=\{(\gamma_1, \gamma_2); \, \, s(\gamma_1)=
t(\gamma_2)\}\longrightarrow \mathcal{G},$$
which is associative in an obvious sense. \\

\begin{examples} \rm{
i) Every group $G$ defines  a  groupoid  $\mathcal{G}$ with one object * and
$$\text{Hom}_{\mathcal{G}}(*, *)=G.$$
 The composition of morphisms is simply by  group multiplication.\\

\noindent ii) Let $X$ be a $G$-set.
 We define a groupoid
$\mathcal{G}=X\rtimes G$, called the {\it transformation  groupoid} of the action, as
follows. Let $\text{obj} \;\mathcal{G} =X$, and
$$\text{Hom}_{\mathcal{G}} \,(x, y)=\{g\in G; \; gx=y\}.$$
Composition of morphisms is  defined via group multiplication. It is
easily checked that $\mathcal{G}$ is a groupoid. Its set of
morphisms can be identified as
$$\mathcal{G}\simeq X \times G,$$
where the  composition of morphisms is given by
$$(gx, h) \circ (x, g)= (x, hg).$$

\noindent iii) Let $\sim$ denote an equivalence relation on
a set $X$. We define a groupoid $\mathcal{G}$, called the {\it graph of $\sim$}, as
follows. Let $\text{obj} \; \mathcal{G}= X$, and let
$\text{Hom}_{\mathcal{G}} \, (x, y)$ be a one element set if $x\sim y$, and
be the empty set otherwise.

Note that the set of morphisms of $\mathcal{G}$ is identified with the
graph of the relation $\sim$ in the usual sense:
$$\mathcal{G}=\{ (x, y); \; x\sim y \}\subset X\times X.$$}
\end{examples}

The {\it groupoid algebra} of  a  groupoid $\mathcal{G}$ is the algebra of functions
$$\mathbb{C} \mathcal{G} = \{ f: \mathcal{G} \to \mathbb{C}; \;f
\mbox{ has finite support}\},$$
 with finite support on $\mathcal{G}$.
Under  the
{\it convolution product}
$$ (fg) (\gamma)=\sum_{\gamma_1 \circ \gamma_2
=\gamma}f(\gamma_1)g(\gamma_2),$$
and the involution
$$f^* (\gamma)= \overline{f (\gamma^{-1})}$$
it is a $*$-algebra.

\begin{examples} \rm{ i) Let  $\mathcal{G}$ be the groupoid defined by a group $G$.  Then clearly  the groupoid algebra
$ \mathbb{C} \mathcal{G}$ is isomorphic to the group algebra $\mathbb{C}
G$.

\noindent ii) If  $\mathcal{G}= X \rtimes G$ is  a  transformation
groupoid  then we have an algebra isomorphism
$$ \mathbb{C} \mathcal{G} \simeq C(X)\rtimes G.$$

\noindent iiI) Let $\mathcal{G}$ be the
  {\it groupoid of pairs}  on a set of $n$ elements, i.e.
 $$\mathcal{G}=\{(i, j); \; i, j=1, \cdots n\}$$
 with composition given by
 $$(l, k)\circ (j, i)= (l, i) \quad  \mbox{if}\; k=j.$$
 (Composition is not defined otherwise).
We have
 $$\mathbb{C}\mathcal{G}\simeq M_n(\mathbb{C}).$$
 Indeed, it is easily checked that the map
 $$ (i, j)\mapsto E_{i, j},$$
 where $E_{i, j}$ denote the matrix units, is an algebra
 isomorphism. This is an extremely  important example. In fact,  as
 Connes points out in \cite{cob}, the matrices in Heisenberg's matrix
 quantum mechanics \cite{hei}
  were arrived at by a similar procedure. }
\end{examples}

\begin{exercise} Show that the groupoid algebra of a finite groupoid
(finite set of objects and finite set of morphisms) can be
decomposed as a  direct sum of tensor products of group algebras and
matrix algebras.
\end{exercise}

As the above exercise shows, one cannot get very far with just discrete
groupoids and soon one  needs   to work with topological and smooth  groupoids
associated to, say,  continuous actions of topological groups and to
foliations.

A (locally compact) {\it topological groupoid} is a groupoid such
that its set of morphisms $\mathcal{G}$  and set of objects
$\mathcal{G}^{(0)}$ are  (locally compact) topological spaces, and
its  composition, source, target and inverse maps are continuous.
An {\it \'{e}tale groupoid} is a locally compact  groupoid such
that  the fibers of its  target map $ \mathcal{G}^x= t^{-1}(x),
\quad x\in \mathcal{G}^{(0)},$ are discrete.

A {\it smooth  groupoid}, also known as {\it Lie groupoid},    is
a groupoid such that $\mathcal{G}$ and $\mathcal{G}^{(0)}$ are
smooth manifolds, the inclusion $\mathcal{G}^{(0)} \to
\mathcal{G}$ as well as the maps $s, t, i$ and the composition map
$\circ$ are smooth, and $s$ and $t$ are submersions. This last
condition will guarantee that the domain of the composition map
$\mathcal{G}^{(2)}=\{(\gamma_1, \gamma_2); s(\gamma_1)=
t(\gamma_2)\}$ is a smooth manifold.\\

To define the convolution algebra of a topological groupoid and
its $C^*$-completions, we need an analogue of Haar measure for
groupoids.
 A {\it  Haar measure } on a locally compact groupoid   $\mathcal{G}$ is
 a  family of measures $\mu^x$ on each $t$-fiber $\mathcal{G}^x$. The
 family is supposed to be continuous and left invariant in an obvious sense (cf. \cite{ren},   unlike locally compact
 topological groups, groupoids  need not have an invariant Haar measure).
  Given a Haar measure, we can then define, for functions with compact support $f, g \in
 C_c (\mathcal{G})$ their convolution product
\begin{equation} \label{cpg}
 (f * g) (\gamma)=\int_{\gamma_1 \gamma_2=\gamma} f (\gamma_1) g
 (\gamma_2)=\int_{\mathcal{G}^{t (\gamma)}} f ( \gamma_1) g(\gamma_1^{-1} \gamma) d\mu^{t(\gamma)}.
 \end{equation}
 This turns $C_c (\mathcal{G})$ into a $*$-algebra. The involution is
 defined by $f^*(\gamma)=\overline{f(\gamma^{-1}))}$. For each fiber
 $\mathcal{G}^x$, we have an $*$-representation  $ \pi_x$ of
 $C_c (\mathcal{G})$ on the Hilbert space $L^2 (\mathcal{G}^x, \mu^x)$ defined
 by
 $$ (\pi_x  f)(\xi)(\gamma)=\int_{\gamma_1 \gamma_2=\gamma} f (\gamma_1)
 \xi (\gamma_2)=\int_{\mathcal{G}^{t (\gamma)}} f ( \gamma_1)
 \xi(\gamma_1^{-1} \gamma) d\mu^{t(\gamma)}.$$
 We can then define a pre $C^*$-norm on $C_c (\mathcal{G})$ by
 $$\| f \|:
 = \text{sup} \, \{\|\pi_x ( f )\|; \, x \in \mathcal{G}^0\}.$$
  The completion
 of $C_c (\mathcal{G})$ under this norm is  the {\it
 reduced $C^*$-algebra} of the groupoid $\mathcal{G}$ and will be denoted
 by $ C^*_r (\mathcal{G})$.

 There are two special cases that are particularly important and convenient to work with: \'{e}tale and smooth groupoids.
 Notice that
 for an \'{e}tale groupoid each fiber is a discrete set and  with counting measure on each fiber
  we obtain a Haar measure. The convolution product is then given by
 $$ ( f *  g)(\gamma) =\sum_{\gamma_1 \gamma_2=\gamma} f (\gamma_1) g
 (\gamma_2)=\sum_{\mathcal{G}^{t (\gamma)}} f ( \gamma_1) g (\gamma_1^{-1} \gamma).$$
Notice that for each $\gamma$ this is a finite sum since the
support of $f$ is compact and hence contains only finitely many
points of each fiber.

A second interesting case where one can do away with Haar measures
is smooth groupoids. Let $C^{\infty}_c (\mathcal{G}, \,
\Omega^{\frac{1}{2}})$ be the space of sections, with compact
support, of the line bundle of half-densities  on a smooth
groupoid $\mathcal{G}$. Since the product of two half-densities is
a 1-density which has  a well defined integral, the integral \eqref{cpg}  for
$ f, \,  g \in C^{\infty}_c (\mathcal{G}, \, \Omega^{\frac{1}{2}})$ is well
defined and we obtain the smooth convolution algebra $C^{\infty}_c
(\mathcal{G})$.

\begin{examples} \rm{ 1. Let $\Gamma$ be a discrete group acting by
homeomorphisms on a locally compact space $X$. Then the
transformation groupoid $ X \rtimes \Gamma$ is an \'{e}tale groupoid
and  the groupoid algebra recovers the crossed product algebra:
$$C_c(\mathcal{G}) \simeq  C_c(X)\rtimes \Gamma \quad  \text{and} \quad
C^*_r(\mathcal{G}) \simeq C_0(X)\rtimes_r \Gamma.$$
For $X=S^1$ and
$\Gamma =\mathbb{Z}$ acting through  rotation by an angle $\theta$, we
recover the noncommutative torus as a groupoid algebra,
which is one among many incarnations of $A_{\theta}$.  \\

\noindent 2. Let $X$ be a locally compact space with a Borel
probability measure $\mu$ and $\mathcal{G}$ be the groupoid of
pairs on $X$. Then for  $ f, \, g  \in C_c (X \times X)$ the
convolution product \eqref{cpg}  reduces to
$$ (f *  g )(x, z) =\int_X f (x, y) g (y, z) d\mu (y),$$
which is reminiscent of matrix multiplication or  products
of integral operators.  In fact the map  $ T: C_c (X \times X) \to
\mathcal{K} (L^2 (X, \mu))$ sending $f$ to the integral operator
$$ (Tf) (g) (x) = \int_X f (x, y) g (y) d\mu (y)$$
is clearly an
algebra map and can be shown to be 1-1 and onto.

On the other
extreme, 
if $\mathcal{G}$ is the groupoid of the discrete equivalence relation on a locally compact space  $X$, also known as the {\it groupoid
of pairs}, 
then clearly  $C_c(\mathcal{G}) \simeq C_c(X).$}  \\
\end{examples}

\begin{example} (Non-Hausdorff manifolds) \rm{ Let
$$X = S^1 \times 0 \, \cup \, S^1 \times 1$$
be the disjoint union of two copies of the circle. We identify
$(x, 0) \sim (x, 1)$ for all $x \neq 1$ in $S^1$. The quotient
space $X / \sim$ is a non-Hausdorff manifold.  The groupoid of the
equivalence relation  $\sim$
$$\mathcal{G} \, = \, \{ (x, y)\,  \subset X \times X; \, x \sim y \, \}$$
is a smooth \'{e}tale groupoid. Its smooth groupoid
algebra is given by
$$ C^{\infty} (\mathcal{G}) = \, \{ f \in C^{\infty}( S^1, \,  M_2 (\mathbb{C})); \, \, f(1) \, \, \text{is diagonal}\,\}$$ }

\end{example}

There are many interesting examples of  noncommutative quotients
that we did not discuss here but are of much interest in
noncommutative geometry. They include: foliation algebras, the space
of Penrose tilings, the ad\`{e}le space and the space of $\mathbb{Q}$-lattices in number theory. They
can all be defined as groupoid algebras and variations thereof. We
refer to Connes-Marcolli's article in this volume as well as
\cite{coma1, coma2, coma3, cob} for a proper  introduction.

\begin{exercise} Show that the Hecke algebras $\mathcal{H}(\Gamma, \, \Gamma_0)$ defined in Example 2.7 are  groupoid algebras.
\end{exercise}

The following result of M. Rieffel \cite{rie} clarifies the relation between the classical
quotients and noncommutative quotients for group actions:

\begin{theorem} Assume $G$ acts  freely and  properly
on a locally compact Hausdorff space $X$. Then we have a strong Morita
equivalence between the $C^*$-algebras
$C_0(X/G)$ and $C_0(X)\rtimes_r G.$
\end{theorem}

\subsection{Hopf algebras and quantum groups}
Many examples   of noncommutative spaces are  Hopf algebras or
quantum groups. They are either the algebra of coordinates of  a
quantum group, or, dually, the convolution algebra or the
enveloping algebra of a quantum group. In this  section  we shall
make no distinction between Hopf algebras and quantum groups and
use these words interchangeably.   The  theory of Hopf algebras
(as well as Hopf spaces) was born in the paper of H. Hopf  in his
celebrated computation of the rational cohomology of compact
connected Lie groups \cite{ho}. This line of investigation
eventually led to the Cartier-Milnor-Moore theorem \cite{mimo}
characterizing  connected cocommutative Hopf algebras as
enveloping algebras of Lie algebras.

A purely algebraic theory,
with motivations independent from algebraic topology,  was created
by Sweedler in the 1960's.  This line of investigation took a big leap
forward with the work of Drinfeld and Jimbo resulting in
quantizing all classical Lie groups and Lie algebras \cite{dr}.

In a
different direction, immediately after Pontryagin's duality
theorem for locally compact abelian  groups, attempts were made
 to extend it to noncommutative groups. The Tannaka-Krein duality theorem was an important first step.
 Note that the dual, in any sense of the word, of a noncommutative group is necessarily not a group and  one is naturally interested
 in extending the category of groups to a larger category  which is  closed under  duality and hopefully is even equivalent to its second dual. Hopf von Neumann
 algebras
 of G.I.  Kac and Vainerman  achieve this in the measure theoretic world of von Neumann algebras \cite{ensc}.
 The theory of {\it locally compact quantum groups} of Kustermans and Vaes
 \cite{kuva}
 which was developed much later achieves this goal in the category of $C^*$-algebras.
 An important step towards this  program was the theory of compact quantum
 groups of S. L. Woronowicz (cf. \cite{wo} for a survey). We refer to
 \cite{kas, kl, mj, ma88, ma91,  sw} for the general theory of Hopf algebras
 and quantum groups.

 The first serious interaction  between Hopf algebras and noncommutative geometry
 started in earnest in the paper of Connes and Moscovici on transverse
 index theory \cite{como1} (cf. also \cite{como2, como3, como4} for
 further developments). In this paper a
   noncommutative and non-cocommutative Hopf algebra appears as the quantum symmetries of
   the noncommutative space of codimension one foliations. The same Hopf algebra was
   later shown to act on the noncommutative space  of modular Hecke
   algebras \cite{como5}.

To understand the definition of a Hopf algebra, let us see  what
kind of extra structure exists  on the algebra of functions on a
group. For simplicity, let $G$ be a  finite group, though this is by no
means necessary,  and let
$H=C(G)$ be the algebra of functions from  $G \to \mathbb{C}$. The
multiplication $m: G\times G \to G$,  inversion $i: G \to G$, and
unit element $e\in G$,  once dualized,  define  unital algebra maps
$$\Delta: H \to H\otimes H, \quad S: H\to H, \quad  \epsilon: H\to
\mathbb{C},$$
by the formulas
$$\Delta f(x, y)=f(xy), \quad Sf(x)=f(x^{-1}), \quad \varepsilon
(f)=f(e),$$
where we have identified $C(G\times G)$ with $C(G)\otimes C(G)$.  Let
also
$$ m: H \otimes H \to H,  \quad \quad \eta: \mathbb{C} \to H$$
denote the multiplication and unit maps of the algebra $H$.
 The associativity,
inverse, and unit axioms  for  groups are dualized and in fact are easily seen to be equivalent to the following
 {\it coassociativity, antipode, and counit axioms} for $H$:
\begin{eqnarray*}
(\Delta \otimes I) \Delta &=&(I\otimes \Delta) \Delta ,\\
 (\varepsilon \otimes I) \Delta &=&(I\otimes \varepsilon) \Delta
=I,\\
m (S\otimes I)&=&m(I\otimes S)=\eta \varepsilon.
\end{eqnarray*}

\begin{definition} A unital algebra $(H, \,m, \, \eta)$ endowed with unital algebra homomorphisms
$\Delta : H \to H \otimes H,$ \,  $\varepsilon : H \to \mathbb{C}$ and a linear
map $S : H \to H$ satisfying  the above equations  is called a Hopf algebra.
\end {definition}
It can be shown that the antipode $S$ is unique and is an anti-algebra map. It
is also an anti-coalgebra map.

\begin{example} \rm{
 Commutative or  cocommutative Hopf algebras are closely related to
groups and Lie algebras. We give a few examples to indicate this
connection.\\

\noindent i)  Let $\Gamma$ be a discrete group (it need not be finite) and $H=\mathbb{C}\Gamma$
the group algebra of $\Gamma$. Let
$$ \Delta (g)=g\otimes g, \quad S(g)=g^{-1}, \quad
\varepsilon (g)=1,$$ for all $g\in \Gamma$ and linearly extend them to
$H$. Then it is easy to check that $(H, \, \Delta, \, \varepsilon, \,  S)$
is a cocommutative Hopf
algebra. It is  commutative if and only if $\Gamma$ is commutative.\\

\noindent ii) Let $\mathfrak{g}$  be a Lie algebra and
$H=U(\mathfrak{g})$ be  the universal enveloping algebra of
$\mathfrak{g}$. Using the universal property of $U(\mathfrak{g})$
one checks that there are uniquely defined  algebra homomorphisms
$\Delta: U(\mathfrak{g})\to U(\mathfrak{g})\otimes
U(\mathfrak{g})$, $\varepsilon: U(\mathfrak{g})\to \mathbb{C}$ and an
anti-algebra map $S:U(\mathfrak{g})\to U(\mathfrak{g})$,
determined by
$$ \Delta (X)=X\otimes 1 +1\otimes X, \quad \varepsilon (X)=0, \quad \text{and} \; \; S(X)=-X,$$
for all $X \in \mathfrak{g}$.  One then checks easily that
$(U(\mathfrak{g}), \Delta, \varepsilon, S)$ is a cocommutative Hopf
algebra. It is commutative if and only if $\mathfrak{g}$ is an abelian Lie algebra,
in which case $U(\mathfrak{g})=S(\mathfrak{g})$ is the symmetric algebra of $\mathfrak{g}$. \\

\noindent iii) (Compact groups) Let $ G$ be a compact topological group. A continuous
function $f:G \to \mathbb{C}$ is called {\it representable} if the
set of left translations of $f$ by all elements of $G$ forms a finite
dimensional subspace of  $C(G)$.   It is easy to see  that the set of
 representable functions,
 $H = \text{Rep} \,(G)$, is a subalgebra of
 $C(G)$. Let $m: G\times G \to G$ denote the multiplication of $G$
 and $m^*: C(G\times G)\to C(G), \; m^*f (x, y)=f(xy), $ denote its  dual map. One checks  that
  if $f$
 is representable, then
 $$m^*f \in \text{Rep}\, (G)\otimes \text{Rep}\,  (G) \subset C(G\times G).$$
 Let $e$ denote  the identity of $G$.
   The formulas
 $$ \Delta f =m^*f, \quad \varepsilon f=f(e), \quad \text{and} \; \;
 (Sf)(g)=f(g^{-1}),$$
 define a Hopf algebra structure on $Rep (G)$. Alternatively, one can describe $\text{Rep}\,  (G)$ as the linear
 span of {\it matrix coefficients}  of all finite dimensional complex representations of $G$.
 By the {\it Peter-Weyl Theorem,} $\text{Rep} (G)$ is a dense
 subalgebra of $C(G)$. This algebra is finitely
 generated (as an algebra) if and only if $G$ is a Lie group.\\

\noindent iv) (Affine group schemes) The coordinate ring $H=\mathbb{C} [G]$ of an affine algebraic group $G$ is
a commutative Hopf algebra. The maps $\Delta$, $\varepsilon$, and
$S$ are the duals of the multiplication, unit, and inversion  maps
of $G$, respectively. More generally, an {\it affine group
scheme} over a commutative ring $k$ is a commutative Hopf
algebra over $k$. Given such a Hopf algebra $H$, it is easy to see
that for any commutative $k$-algebra $A$, the set $\text{Hom}_{Alg}\,(H,
A)$ is a group under convolution product and $A\mapsto
\text{Hom}_{Alg}\,(H, A)$ is a functor from the category $ComAlg_k$ of
commutative $k$-algebras to the category of groups. Conversely,
any {\it representable functor} $ComAlg_k \to Groups$ is represented by
a, unique up to isomorphism,  commutative $k$-Hopf algebra. Thus
the category  of affine group schemes is equivalent to the
category of representable
functors $ComAlg_k \to Groups.$  }
\end{example}

\begin{example} (compact quantum groups) \rm{ A prototypical example is Woronowicz's $SU_q(2)$, for $0< q\leq 1$. As a $C^*$-algebra
it is the unital $C^*$-algebra generated by $\alpha$ and $\beta$
subject to the relations
$$ \beta \beta^*= \beta^* \beta,  \, \, \, \alpha \beta =q \beta \alpha,
\, \, \, \alpha \beta^* =q \beta^* \alpha, \, \,\,
\alpha \alpha^* + q^2 \beta^*  \beta = \alpha^* \alpha +  \beta^*   \beta =I.$$
Notice that these relations amount to saying that
$$U= \left( \begin{array}{cc} \alpha &  q \beta\\ -\beta^* &
\alpha^*\end{array}\right)$$ is unitary, i.e. $UU^*=U^* U=I.$  Its
coproduct and antipode are  defined by
$$ \Delta  \left( \begin{array}{cc} \alpha & \beta\\ -\beta^* &
\alpha^*\end{array}\right)= \left( \begin{array}{cc} \alpha & \beta\\ -\beta^* &
\alpha^*\end{array}\right) \otimes \left( \begin{array}{cc} \alpha & \beta\\ -\beta^* &
\alpha^*\end{array}\right)$$
$$ S(\alpha)= \alpha^*, \, \, \, S(\beta)= -q^{-1} \beta^*, \, \,\,
S(\beta^*)= -q\beta, \, \, \,
S(\alpha^*)=\alpha.$$
  Notice that the coproduct is only defined on the algebra
  $\mathcal{O}(SU_q (2))$ of matrix elements on the quantum group, and its
  extension to $C(SU_q(2))$ lands in the completed tensor product
  $$ \Delta: C(SU_q(2)) \longrightarrow C(SU_q(2)) \hat{\otimes} C(SU_q(2)).$$
   At $q=1$ we obtain the algebra
  of continuous functions on $SU(2)$. We refer to \cite{kuva} for a  survey  of
  compact and locally compact quantum groups.  }
\end{example}

\begin{example} \rm{ (Connes-Moscovici Hopf algebras) A very important
example for noncommutative geometry and its applications to
transverse geometry and number theory is the family of {\it
Connes-Moscovici Hopf algebras} $\mathcal{H}_n$ for $n\geq 0$\,
\cite{como2, como3, como4}. They are deformations of the group $G= \,
\text{Diff}\, (\mathbb{R}^n)$ of diffeomorphisms of $\mathbb{R}^n$ and can
also be thought of as deformations of the Lie algebra
$\mathfrak{a}_n$ of formal vector fields on $\mathbb{R}^n$. These
algebras appeared for the  first time as quantum symmetries of
transverse
 frame bundles of codimension $n$ foliations. We briefly treat the case $n=1$ here.
 The main features of  $\mathcal{H}_1$ stem from the fact that
 the group  $G= \,
\text{Diff}\, (\mathbb{R}^n)$ has a {\it factorization} of the form
$$G=G_1G_2,$$
where $G_1$ is the subgroup of diffeomorphisms $\varphi$  that satisfy
 $$\varphi (0)=0, \quad \varphi'(0)=1,$$
and $G_2$ is the $ax+b$- group of affine diffeomorphisms. We introduce two Hopf algebras
corresponding to $G_1$ and $G_2$ respectively.   Let $F$
denote the Hopf algebra
of polynomial   functions on the pro-unipotent group
$G_1$. It can also be defined  as the {\it continuous dual}  of the
enveloping algebra of the Lie algebra of $G_1$. It is a commutative Hopf
algebra generated by the Connes-Moscovici coordinate functions
$\delta_n$, $n=1,2,\dots$, defined by
$$\delta_n (\varphi)=\frac{d^n}{dt^n}(\text{log}\,  (\varphi'(t))|_{t=0}.$$
  The second Hopf algebra,
$U$,  is the universal enveloping algebra of the Lie algebra
$\mathfrak{g}_2$ of the $ax+b$-group. It has generators $X$ and $Y$ and one relation $[X,
Y]=X$.

 The factorization $G=G_1 G_2$ defines a {\it matched pair} of Hopf algebras consisting
 of $F$ and $U$. More precisely, The Hopf algebra $F$ has a right  $U$-module
 algebra structure
 defined by
 $$\delta_n (X)=-\delta_{n+1}, \quad \text{ and}\;\;
\delta_n (Y)=-n\delta_n.$$
 The Hopf algebra $U$, on the other hand, has  a left
$F$-comodule coalgebra structure via
$$X\mapsto 1\ot X+ \delta_1\ot X, \quad \text{ and}\;\;
Y\mapsto 1\ot Y.$$
One can check that they are a matched pair of  Hopf
algebras in the sense of G.I. Kac and Majid \cite{mj} and the resulting bicrossed product Hopf algebra
$$ \mathcal{H}_1 = F\bowtie U$$
is
 the Connes-Moscovici Hopf algebra $\mathcal{H}_1$. (See \cite{como2} for a slightly different
approach and fine points of the proof.)

 Thus $\mathcal{H}_1$
is the universal Hopf algebra generated by  $\{X, \, Y, \, \delta_n ; \, \, n =
1, 2, \cdots \}$  with relations
$$[Y, \,X] = X, \; \;[Y, \, \delta_n] = n \delta_n, \; \; [X, \, \delta_n] =
\delta_{n+1},\; \;
[\delta_k, \, \delta_l] = 0,$$
 $$\Delta Y = Y \otimes 1 + 1 \otimes  Y, \quad  \Delta \delta_1 =
  \delta_1 \otimes  1 + 1 \otimes  \delta_1, $$
$$\Delta X = X \otimes  1 + 1 \otimes  X + \delta_1 \otimes  Y,$$
$$ S(Y ) = -Y, \; \;  S(X) = -X + \delta_1Y, \; \; S(\delta_1) =
-\delta_1, $$
for $ n, k, l = 1 ,
2, \cdots.$}
\end{example}

Another recent point of interaction between Hopf algebras and noncommutative
geometry is the work of Connes and Kreimer in renormalization schemes of  quantum field theory.   We refer to
\cite{cokr1, cokr2, cokr3, coma2, coma3, coma4}  and references therein for this fascinating new subject.

\section{Topological $K$-theory}

The topological $K$-theory of spaces and its main theorem, {\it the
Bott periodicity theorem}, can be  extended to noncommutative
Banach algebras. Of all topological invariants of spaces,
$K$-theory has the distinct feature that it is the easiest to
extend to noncommutative spaces. Moreover, on a large class of
$C^*$-algebras  the  theory can be characterized by a few simple
axioms.  In the next section we take up  the  question of
 Chern character in noncommutative geometry. It is to
address this and similar questions that cyclic cohomology and
Connes' Chern character map enter the game.

$K$-theory was first introduced by Grothendieck in 1958 in his extension
of the Riemann-Roch theorem to algebraic varieties.
The isomorphism classes of bounded complexes of coherent sheaves on a variety $X$
form an abelian monoid and the group that they generated was called $K_0(X)$.
Soon after, Atiyah and Hirzebruch realized that in a similar fashion
 complex vector bundles over a compact
space $X$ define a group $K^0 (X)$ and, moreover,  using standard
methods of algebraic topology, one obtains a {\it generalized
cohomology theory} for spaces in this way. Bott's periodicity
theorem for  homotopy groups of stable unitary groups  immediately
implies that  the new functor is  2-periodic. By the mid-1970's it
was clear to operator algebraists that
 topological $K$-theory and Bott periodicity theorem  can be extended to all
Banach algebras. Our references  for this section include \cite{bl1, fi, cob}.

\subsection{The $K_0$ functor }
Since the definition of $K_0 (A)$  depends only on the underlying ring
structure of $A$ and makes
sense for any ring,
 we shall define $K_0(A)$ for any ring $A$.  Let $A$ be a unital noncommutative  ring.   A right
 $A$-module $P$ is called {\it projective} if it is a direct
 summand of a free module, i.e.
 there exists a right
$A$-module $Q$ such that
$$P\oplus Q\simeq  A^I.$$
 Equivalently, $P$ is projective if and only if any short exact sequence of right $A$-modules
$$0\longrightarrow M\longrightarrow N \longrightarrow P\longrightarrow
0$$
splits. Let $\mathcal{P}(A)$ denote the set of isomorphism
 classes of finitely generated projective right $A$-modules. Under the
 operation of direct sum,  $\mathcal{P}(A)$ is an abelian monoid. The group $K_0(A)$ is, by definition, the universal group
 generated by the monoid $\mathcal{P}(A)$. Thus elements of $K_0 (A)$
 can be written as $[P]-[Q]$ for $ P,\, Q \in \mathcal{P}(A)$, \,  with
 $[P]-[Q]= [P']-[Q']$\,  if and only if there  is an  $R  \in
 \mathcal{P}(A)$ such that $ P\oplus Q' \oplus R \simeq P' \oplus Q
 \oplus R$.

 A unital ring
 homomorphism $f: A \to B$ defines a map ({\it base change}) $f_* :
 \mathcal{P}(A) \to \mathcal{P}(B)$ by
 $$f_* (P)= P\otimes_A B$$
 where the left $A$-module structure on $B$ is induced by $f$. This map
 is clearly additive and hence induces an  additive map
 $$ f_* : K_0(A) \to K_0(B).$$
 This shows that $A \to K_0(A)$ is a functor.

 We need to define $K_0$ of
 non-unital rings. Let $A^+$ be the {\it unitization} of a non-unital
 ring $A$. By definition, $A^+=A \oplus \mathbb{Z}$ with multiplication
 $(a, m)(b, n)= (ab+na+mb, \, mn)$ and unit element $(0, 1)$.  A non-unital ring map $f : A \to B$
 clearly induces a unital ring map $ f^+ : A^+ \to B^+$ by $ f^+ (a,n) =(f(a), n)$. The canonical morphism
 $A^+ \to \mathbb{Z}$, sending $(a, n)\to n$, 
 is unital and we define $K_0(A)$ as the kernel of the induced map
 $K_0(A^+) \to K_0(\mathbb{Z})$. If $A$ is already unital then the
 surjection $A^+ \to \mathbb{Z}$ splits and one shows that the two
 definitions coincide.

The first important result about $K_0$ is its {\it half-exactness}
(cf., for example,  \cite{bl1} for a proof): for any exact
sequence of rings
$$0\longrightarrow I\longrightarrow A \longrightarrow A/I \longrightarrow
0,$$
the induced sequence
\begin{equation}\label{hes}
 K_0(I)\longrightarrow K_0 (A) \longrightarrow K_0 (A/I)
 \end{equation}
is {\it exact in the middle}. Simple examples show that exactness at the 
left and right ends can fail and in fact the extent to which they fail is
measured by higher $K$-groups as we define them in the next section.

\begin{remark} When $A$ is commutative, the tensor product $P\otimes_A
Q$ of $A$-modules is well defined and is an $A$-module again. It is finite and projective if $P$ and
$Q$ are finite and projective. 
This operation turns $K_0(A)$
into  a commutative ring. In general, for noncommutative rings
no such multiplicative structure exists on $K_0(A)$.
\end{remark}

 There is an alternative description of $K_0(A)$ in terms of idempotents in
 matrix algebras over $A$ that is often convenient. An idempotent  $e\in
 M_n(A)$ defines a right $A$-module map
 $$e: A^n \longrightarrow A^n$$
 by left multiplication by $e$. Let $P_e= eA^n$ be the image of $e$. The
 relation
 $$ A^n =eA^n \oplus (1-e)A^n$$
  shows that $P_e$ is a finite
 projective right $A$-module. Different idempotents  can define isomorphic
 modules.  This happens, for example,  if $e$ and $f$ are {\it equivalent idempotents} (sometimes called {\it similar})
 in the sense that
 $$ e=ufu^{-1}$$
 for some invertible $u\in GL(n, A)$. Let $M  (A)=\cup M_n(A)$ be the direct limit of  matrix algebras $M_n(A)$
 under the embeddings $M_n(A) \to M_{n+1}(A)$ defined by $a \mapsto \left( \begin{array}{cc} a& 0\\ 0 & 0\end{array}\right)$. Similarly
 let $GL (A)$ be the direct limit of the groups $GL(n, A)$. It acts on
 $M(A)$ by conjugation.

 \begin{definition}Two idempotents $e \in M_k(A)$ and $f\in M_l(A)$ are called {\it stably equivalent} if their images in $M(A)$
 are equivalent under the action of $GL(A)$.
 \end{definition}

 The following is easy to prove and answers our original question:
 \begin{lem} The projective modules $P_e$ and $P_f$ are isomorphic if
 and only if the idempotents $e$ and $f$ are stably equivalent.
 \end{lem}

Let $ \text{Idem} (M(A))/GL(A)$ denote the set of stable equivalence classes of idempotents over $A$.
This is an abelian monoid under the operation
$$ (e, \,f) \mapsto e \oplus f: = \left( \begin{array}{cc} e& 0\\ 0 &
f\end{array}\right).$$

It is clear that any finite projective module is of the type $P_e$ for
some idempotent $e$. In fact writing $P\oplus Q \simeq A^n$, one can let
$e$ be
 the idempotent corresponding to the projection map $(p, q) \mapsto (p, 0).$ These
observations  prove the following lemma:

\begin{lem} For any unital ring $A$, the map $e \mapsto P_e$ defines an isomorphism of monoids
$$ \text{Idem} \,  (M(A))/GL(A) \simeq \mathcal{P} (A).$$
\end{lem}
Given an idempotent $e =(e_{ij}) \in M_n(A)$, its image under a
homomorphism $f : A \to B$ is the idempotent $f_* (e)=
(f(e_{ij}))$. This is our formula for $f_* : K_0(A) \to K_0 (B)$
in the idempotent picture of $K$-theory.

For a Banach algebra $A$,  $K_0(A)$ can be described in terms of
connected components of the space of idempotents of $M (A)$ under
its inductive limit topology (a subset  $ V \subset M (A)$ is
open in the inductive limit topology if and only if
 $V \cap M_n (A)$ is open for all $n$).  It is based on the
following important  observation:  Let $e$ and $f$ be idempotents
in a unital Banach algebra $A$ and assume $\|e- f\|< 1/ \|2e-
1\|$. Then $e\sim  f$. In fact  with
$$v=(2e-1)(2f-1)+1$$
and $u=\frac{1}{2} v$, we have $ueu^{-1}=f$. To see that $u$ is
invertible note that $\|u-1\|<1$.  One consequence of this fact is
that if $e$ and $f$ are in the same  path component of the space
of idempotents in $A$, then they are equivalent. As a result we
have, for any Banach algebra $A$, an isomorphism of monoids
$$ \mathcal{P} (A) \simeq \pi_0 (\text{Idem}\, (M (A))),$$
where $\pi_0$ is the functor  of path components.

For $C^*$-algebras, instead of idempotents it suffices  to consider only the {\it projections}. A
projection is a self-adjoint idempotent ($p^2=p=p^*$). The reason is that
every idempotent in a $C^*$-algebra is similar  to a projection \cite{bl1}: let  $e$ be an idempotent
and set $z=1+ (e-e^*)(e^*-e)$. Then $z$ is invertible and positive  and one shows that $p=ee^*z^{-1}$ is a projection and is similar
to $e$.
\begin{exercise} Show that the set of projections of a $C^*$-algebra is
homotopy equivalent (in fact a retraction) of the set of idempotents.
\end{exercise}
Let $\text{Proj} (M (A))$ denote the space of projections in $M
(A)$. We
 have established  isomorphisms of monoids
$$ \mathcal{P}(A) \simeq \pi_0 (\text{Idem}\, (M (A))) \simeq \pi_0 (\text{Proj}\,
(M (A))).$$

From the above homotopic interpretation of $K_0$ for Banach
algebras, its {\it homotopy invariance}  and {\it continuity}
easily follows. Let $f, g : A \to B$ be continuous homomorphisms
between Banach algebras. They are  called {\it homotopic} if there exists a continuous homomorphism
$$ F: A \to C ([0, 1], \,B)$$
such that $f= e_0  F$ and $g= e_1 F$, where $e_0,  \, e_1 : C( [0, 1], B)
\to B$ are the evaluations at 0 and 1 maps. Now by our definition of $K_0$
via $\pi_0$,  it is clear that   $e_{0*} =e_{1*} : K_0 ( C( [0, \, 1], \, B)) \to K_0
(B)$ and hence
$$ f_* =g_* : K_0(A) \to K_0 (B), $$
which shows that $K_0$ is homotopy invariant.

In a similar way one can
also show that $K_0$ preserves direct limits of Banach algebras:
if  $A = \underset{\longrightarrow}{\text{Lim}}  \, ( A_i, \,
f_{ij})$ is an inductive limit of Banach algebras then $ K_0 (A) =
\underset{\longrightarrow}{\text{Lim}} \, ( K_0 (A_i), \, f_{ij
*}).$ This property is referred to as {\it continuity} of $K_0$.

In addition to its homotopy invariance and continuity, we collect a couple of other properties of $K_0$ which hold for all rings: \\

$\bullet$ {\it Morita Invariance}: if $A$ and $B$ are  Morita
equivalent unital rings then $ K_0(A)\simeq K_0(B).$ This is clear
since Morita equivalent rings, by definition,  have equivalent
categories of modules and the equivalence can be  shown to
preserve the categories of finite projective modules.
Therefore $\mathcal{P}(A) \simeq \mathcal{P}(B)$. \\

$\bullet$ {\it Additivity}:  $K_0 (A \oplus B) \simeq K_0(A) \oplus K_0 (B)$
for unital rings $A$ and $B$. This is a consequence of $ \mathcal{P}(A \oplus B) \simeq \mathcal{P}(A) \oplus \mathcal{P}(B),$
which is easy to prove.

 \begin{example} (commutative algebras) \rm{ For $A=C_0(X)$ we have
\begin{equation} \label{sst}
K_0(C_0(X))\simeq K^0(X),
\end{equation}
where $K^0$ is the topological $K$-theory of spaces. The reason
for this is the {\it Swan theorem} \cite{swa} (cf. also Serre \cite{ser} for the corresponding result in the context
of affine varieties),  according to
which for any compact Hausdorff space $X$ the category of finite
projective  $C(X)$-modules is equivalent to the category of
complex vector bundles on $X$.  The equivalence is via  the {\it
global section functor}.
 Given a vector bundle $p: E\rightarrow
X$,  let
$$P= \Gamma (E)=\{s: X\rightarrow E;\;  ps= \text{id}_X \}$$
be the set of all continuous global sections  of $E$. It is clear
that under fiberwise scalar multiplication and addition, $P$ is a
$C(X)$-module. Using the local triviality of $E$  and a partition
of unity one shows that there is a vector bundle $F$ on $X$ such
that $E\oplus F$ is a trivial bundle, or, equivalently,
  $$P \oplus Q\simeq A^n,$$ where $Q$ is the module of global sections
  of $F$. This shows that
  $P$ is finite and projective. It is not difficult to show that all
  finite projective modules are obtained in this way and $\Gamma$ is an equivalence of categories
  (see exercise below).
   Now the rest of the proof of \eqref{sst} is clear
  since $K^0(X)$ is, by definition, the universal group  defined by the
  monoid of complex vector bundles on $X$. }
  \end{example}

\begin{exercise}
 Given  a finite projective $C(X)$-module $P$, let $Q$ be a $C(X)$-module
such that $P \oplus Q\simeq  A^n,$ for some integer $n$. Let $e:A^n\rightarrow A^n$ be the right $A$-linear
projection map    $(p, q)\mapsto (p, 0)$. It is
obviously an
idempotent in $M_n(C(X))$. Define the vector   bundle $E$ to be the
image
of $e$:
$$E=\{(x, v); \; e(x)v=v, \mbox{for all}\,\,  x \in X, \,v\in \mathbb{C}^n\}\subset X\times \mathbb{C}^n.$$
Now it is easily seen  that $\Gamma (E)\simeq P$. With some more work it is shown that the functor $\Gamma$
is full and faithful and hence defines  an equivalence of categories \cite{swa}.
\end{exercise}

Motivated by  the Serre-Swan theorem, one usually thinks of finite
projective modules over  noncommutative algebras as {\it
noncommutative vector bundles}.
\begin{example} \rm{
Here is a nice example of a projection in $M_2 (C_0(\mathbb{R}^2)^+)$. Let
$$p = \frac{1}{1+|z|^2} \left( \begin{array}{cc} |z|^2 & z\\ \bar{z}& 1\end{array}\right)$$
 It does not define an element of $K_0 (C_0(\mathbb{R}^2))$ since
 $p(\infty)=\left( \begin{array}{cc} 1& 0\\ 0 & 0\end{array}\right) \neq 0. $
 The difference
 $$\beta = p-\left( \begin{array}{cc} 1& 0\\ 0 & 0 \end{array}\right)$$
 is however the generator of $K_0 (C_0(\mathbb{R}^2))\simeq \mathbb{Z}$.
 This is a consequence of the Bott periodicity theorem that we recall later in this section.
 $\beta$ is called the {\it Bott generator} and
 $p$ the {\it Bott projection}.  Now we have $C(S^2)= C_0(\mathbb{R}^2)^+$. Let $[1]$ denote the class of the
 trivial line bundle on $S^2$. It follows that  $[1]$ and $\beta$ form a basis
 for $K_0 (C(S^2)) \simeq K^0(S^2)  \simeq \mathbb{Z} \oplus  \mathbb{Z}$.}
 \end{example}

\begin{examples} \rm{ i)  $K_0 (\mathbb{C}) \simeq \mathbb{Z}$. In fact any finite projective module over $\mathbb{C}$ is simply
a finite dimensional complex vector space whose isomorphism class
 is determined by its dimension. This shows that $\mathcal{P} (\mathbb{C})
\simeq \mathbb{N}$,  from which  our claim follows.  \\
ii) By Morita
invariance, we then have $K_0 (M_n (\mathbb{C})) \simeq
\mathbb{Z}$. \\
iii) The algebra of compact operators $\mathcal{K} =
\mathcal{K} (H)$ on a separable Hilbert space is the direct limit
of matrix algebras. Using the continuity of $K_0$ we conclude that
$ K_0 (\mathcal{K}) \simeq \mathbb{Z}$. \\
iv) On the other hand,  $ K_0
(\mathcal{L})=0$ where $\mathcal{L} =\mathcal{L} (H)$ is the
algebra of bounded operators on an infinite dimensional Hilbert
space. To prove this let $e \in M_n ( \mathcal{L})= \mathcal{L}
(H^n)$ be an idempotent.
 The   idempotents $e\oplus I$ and  $0 \oplus I$ are equivalent  in
 $M_{2n}(\mathcal{L}) =\mathcal{L} (H^{2n})$ since both have infinite dimensional range. This shows that
 $[e]=0$.
Notice that for   commutative algebras $A =C(X)$, $K_0(A)$ always contains a copy of $\mathbb{Z}$.
 }
  \end{examples}

  Consider the exact sequence of $C^*$-algebras
  $$ 0 \longrightarrow \mathcal{K} \longrightarrow \mathcal{L} (H)
  \longrightarrow \mathcal{C} \longrightarrow 0,$$
  where $\mathcal{C}: = \mathcal{L} (H)/ \mathcal{K}$ is called  the {\it Calkin algebra}.
  From the above example we see that the corresponding $K_0$ sequence \eqref{hes}
  fails to be exact on the  left.

\begin{example} \label{trkt} (The trace map) \rm{ Let $A$ be a unital algebra,
$V$ be a  vector space  and  $\tau : A \to V$  be a trace on $A$. Then
$\tau$ induces an additive map
$$ \tau : K_0(A) \to V $$
as follows. Given an idempotent $e=(e_{ij})\in M_k(A)$, let
\begin{equation} \label{trk}
\tau ([e])=\sum_{i=1}^k \tau (e_{ii}).
\end{equation}
Using the trace property of $\tau$, one has $\tau (ueu^{-1})=\tau
(e)$. The additivity property $\tau (e \oplus f) = \tau (e) + \tau
(f)$ is clear.  This shows that \eqref{trk} is a well defined map
on $K$-theory.

Alternatively, given a finite
 projective $A$-module $P$, let
 $$\tau ([P])=\text{Tr}\,  (\text{id}_P)$$
 where
 $$\text{Tr}: \, \text{End}_A \, (P)\simeq P^* \otimes_A P   \to V, $$
the
 {\it Hattori-Stallings trace},
  is the natural extension of $\tau$
  defined by
$$ \text{Tr}\,  (f\otimes \xi)= \tau (f(\xi))$$
for all $f\in P^*=\text{Hom}_A (P, \, A)$ and $\xi \in P$. This is the
simplest example
 of a  pairing between cyclic cohomology and $K$-theory (Connes' Chern character),  to be defined in
 full generality
   later in these notes.
Notice that  if $\tau ([e])\neq 0$ then we can conclude that $[e]\neq 0$
in $K_0(A)$. This is often very useful in applications. }
 \end{example}

\begin{example} (Hopf line bundle on quantum spheres) \rm{Let $ 0< q \leq 1$
be a real number. The algebra $C (S_q^2)$ of functions on the {\it standard
 Podle\'{s}  quantum sphere} $S^2_q$ is, by definition,  the
unital $C^*$-algebra generated   by elements $a$ and $b$ with
relations
\[
a a^{\ast}+q^{-4}b^2=1,\; a^{\ast}a+b^2=1,\; ab=q^{-2}ba,\; a^{\ast}b=q^2ba^{\ast}.
\]
The quantum analogue of the Dirac (or Hopf) monopole line bundle over
$S^2$ is given by the following projection
in $M_2(C(S^2_q))$ \cite{bm}:
$$\mathbf{e}_q=\frac{1}{2}\left[\begin{array}{cc} 1+q^{-2}b & qa \\ q^{-1} a^{\ast} & 1-b \end{array}\right].$$
It can be directly checked that $\mathbf{e}_q^2=\mathbf{e}_q =\mathbf{e}_q^*.$  For $q=1$,
$C (S_1^2)=C (S^2)$ and the corresponding projection defines  the Hopf line
bundle on $S^2$. We refer to the article of  Landi  and   van Suijlekom in this volume \cite{lan1} for a survey of noncommutative
bundles and instantons
in general. }
\end{example}

\begin{example} ($K_0(A_{\theta})$) \rm{ We shall see later in this
Section, using the 
Pimsner-Voiculescu exact sequence, that
$$K_0(A_{\theta})  \simeq \mathbb{Z}\oplus \mathbb{Z}.$$
One generator is the class of the trivial idempotent $[1]$. Notice that $[1]\neq 0$
 because $\tau (1)=1 \neq 0$, where $\tau : A_{\theta} \to \mathbb{C}$ is the canonical trace.
  When $\theta$ is irrational a
second generator is given  by the {\it Powers-Rieffel projection}
$p \in A_{\theta}$.
   The projection $p$ is of the form
\begin{equation}\label{prp}
p=U_2^* \, g(U_1)^* +f(U_1) +g(U_1)\,U_2,
\end{equation}
where $f, g \in C^{\infty} (S^1)$. By $g(U_1)$ we mean of course
$\sum \hat{g}_n U_1^n$ where $\hat{g}_n$ are the Fourier
coefficients of $g$. To fulfill the projection condition
$p^2=p=p^*$,   $f$ and $g$ must satisfy certain
relations (cf. \cite{cob})  one of which implies that $\int_0^1 f(t) dt = \theta$.
 There are many such solutions
 but their corresponding
projections are all homotopic and hence define the same class in
$K_0(A_{\theta})$. Now in \eqref{prp},   the only contribution  to
the trace $\tau (p)$ comes  from the constant term of the middle
term and hence
$$ \tau (p)= \int_0^1 f(t) dt =\theta.$$
It follows that the range of the trace map $ \tau: K_0 (A_{\theta}) \to
\mathbb{C}$ is in fact  the subgroup $\mathbb{Z} +\theta \mathbb{Z}
\subset \mathbb{R}$.}
 \end{example}

 \begin{example}  (Relation with Fredholm operators)
 \rm{ The space of Fredholm operators, under the  norm topology,
 is  a {\it classifying space} for the $K$-theory of spaces.  Let $[X, \, \,
 \mathcal{F}] $ denote the set of homotopy classes of continuous maps from a compact space
 $X $ to the space of Fredholm operators $\mathcal{F}$ on an infinite dimensional Hilbert space.
 Such continuous maps should be thought of as  families  of Fredholm
 operators parameterized by $X$. By  a theorem of Atiyah and J\"{a}nich (cf. \cite{atse}
  for a proof and a generalization) there
 exists a well defined {\it index map}  $\text{ index}\,: [X, \, \,
 \mathcal{F}] \to  K^0 (X)$ which
  induces
 an isomorphism
 \begin{equation}\label{ajt}
  \text{index} : [X, \, \mathcal{F}]\simeq K^0 (X).
  \end{equation}
Its definition is as follows.
 Given  a Fredholm family $T: X \to \mathcal{F}$,
 if $ \text{dim} \, \text{Ker} (\, T_x)$ and    $ \text{dim} \,
 \text{Coker} (\, T_x)$ are locally constant functions of $x$, then  the family of finite dimensional subspaces
 $ \text{Ker} (\, T_x)$ and $  \text{Coker} (\, T_x)$,  $x \in X$,   define vector
 bundles denoted
 $\text{Ker}\, (T)$ and  $\text{Coker}\, (T)$  on $X$,  and their difference
$$\text{index} \,(T): = \text{Ker}\, (T) -\text{Coker}\, (T)$$
 is the $K$-theory class associated to $T$. For a  general family the dimensions of the subspaces
 $ \text{Ker} (\, T_x)$ and $  \text{Coker} (\, T_x)$
 may be discontinuous,
 but one shows that it can always be continuously deformed to a family where these dimensions are
  continuous .

 The isomorphism  \eqref{ajt} is  fundamental.  For example, the index of a family of elliptic operators which are fiberwise
 elliptic, by this result, is an element of the $K$-theory of the base manifold (and not an integer).
 In  noncommutative geometry, for example in transverse index theory on foliated manifolds,
 the parameterizing space $X$
 is highly   singular and is replaced by a noncommutative  algebra $A$.  The above
 analytic index map \eqref{ajt}, with values in $K_0(A)$,
 still  can be defined
 and its identification is one of the major problems of  noncommutative index theory \cite{cob}. }
 \end{example}

\subsection{ The higher $K$-functors}
Starting with $K_1$, algebraic and topological $K$-theory begin to differ from each other.  In this section we shall
first briefly indicate the definition of algebraic $K_1$ of rings and then define, for Banach algebras,  a
sequence of functors $K_n$ for $n \geq 1$.

For a unital ring $A$,   let $GL (A)$ be the {\it direct
limit} of groups $GL (n,  A)$ of invertible $n\times n$ matrices over $A$  where   the direct system
$GL (n,  A) \hookrightarrow GL(n+1, A)$ is
 defined by $x \mapsto \left( \begin{array}{cc} x&0\\0 &1
\end{array}\right).$ The algebraic $K_1$ of $A$ is defined as the
abelianization of $GL(A):$
$$ K_1^{\text{alg}}(A):= GL (A)/[GL(A), \, GL(A)],$$
where $ [ \, ,  \, ]$ denotes the commutator subgroup.

Applied to $A= C(X)$, this definition does not reproduce the
topological  $K^1 (X)$. For example for $A =\mathbb{C} =C(\text{pt})$ we
have $K_1^{\text{alg}} (\mathbb{C}) \simeq \mathbb{C^{\times}}$ where the
isomorphism is induced  by the determinant map
$$\text{det}: \,GL (
\mathbb{C}) \to  \mathbb{C^{\times}},$$
  while $K^1(\text{pt})=0$.
  It turns out that, to obtain the right result, one should
 divide $GL(A)$  by a bigger subgroup, i.e. by the {\it  closure} of its commutator subgroup. This works for all Banach algebras
 and will give the right definition of topological $K_1$.  A better approach however is to define the higher $K$ groups in terms
 of $K_0$ and      the {\it suspension functor}.

The {\it suspension} of a Banach algebra $A$ is the Banach
algebra
$$SA = C_0 (\mathbb{R}, \, \,A)$$
 of continuous functions
from $\mathbb{R}$  to $A$ vanishing at infinity.  Notice that for
$A= C(X)$, $SA$ is isomorphic to the algebra of continuous
functions on $X \times [0, 1]$ vanishing on $X \times \{0, 1\}$.
It follows that  $ SA^+ \simeq  C(\Sigma X)$, where  $\Sigma X $
is  the {\it suspension} of $X$   obtained by collapsing $X \times
\{0, 1\}$ to a point in $X \times [0, 1]$.

\begin{definition}
The higher topological  $K$ groups of a Banach algebra $A$ are  the $K_0$
groups of the iterated suspensions of $A$:
$$ K_n(A) =K_0 (S^nA), \quad \quad  n\geq \,1.$$
\end{definition}

This is a bit too abstract. It is better to think of  higher $K$
groups of $A$ as higher homotopy groups of $GL(A)$. To do this  we
need the following  Lemma.
 Let $GL^{\circ} (n,A)$ denote
the connected component of the identity in $GL(n, A)$.

 \begin{lem} i) Let  $f: A \to B$ be a surjective unital homomorphism of unital Banach
algebras.   Then $f: GL^{\circ} (1, A) \to GL^{\circ} (1, B)$ is
surjective.\\
ii) For any $u\in GL(n, A)$,\, \, $\text{diag}\,  (u, u^{-1}) \in GL^{\circ} (2n,
\,A).$
\end{lem}
To prove  the first statement notice that the group generated by the
exponentials $e^y$, $y\in B$, coincides with $GL^{\circ} (1, B)$.
Now since $f$ is surjective  we have $e^y =e^{f(x)}= f (e^x)$
which implies  that any  product of exponentials is in the image
of $f$. To prove the second statement we can use the path
$$ z_t= \text{diag}\,
(u, \,u^{-1}) u_t  \, \text{diag}\,  (u, \,u^{-1}) \, u_t^{-1},$$ where
$$ u_t =\left( \begin{array}{cc} cos \frac{\pi}{2} t & -sin
\frac{\pi}{2} t\\ sin \frac{\pi}{2} t& cos \frac{\pi}{2} t\end{array}
\right),$$
connecting $\text{diag} (u, \,u^{-1})$ to $ \text{diag} (1, \, 1)$.

We can now show  that
\begin{equation}\label{kho}
K_1(A) \simeq \pi_0 (GL(A)),
\end{equation}
 the group of
connected components of $GL(A)$. To see this  let $u\in GL(n, A)$.
Then, by  the above lemma, there is a path  $\alpha_t$
 in $GL (2n, \, A)$ connecting  diag $(u, \, u^{-1})\in GL^0
(2n, \, A)$ to $I_{2n}$. Let $p_n = \text{diag} (I_n, \,0)$. Then
$e_t=\alpha_tp_n\alpha_t^{-1}$ is an idempotent in $ (SA)^+$ and  the map
$ [u] \mapsto  \, [e_t]-[p_n]$ implements the isomorphism in
\eqref{kho}.

Now since $\pi_n (GL (A)) \simeq \pi_{n-1}(GL (SA))$, using  \eqref{kho}
we obtain
\begin{equation} \label{kho1}
K_n (A) \simeq \pi_{n-1} (GL (A)).
\end{equation}

\begin{example} \rm{ Let $A=\mathbb{C}$. Then by \eqref{kho1}, we have
$$K_n(\mathbb{C})\simeq \pi_{n-1} (GL (\mathbb{C}))\simeq \pi_{n-1}
(U(\mathbb{C})).$$
By the Bott  periodicity theorem, the homotopy groups of the  stable unitary
groups $U(\mathbb{C})$ are periodic, i.e.   for all $n$ we have
$$\pi_n (U(\mathbb{C}))\simeq \pi_{n+2} (U(\mathbb{C})),$$
and hence
$$K_{n+2}(\mathbb{C})\simeq K_{n}(\mathbb{C}).$$
This is the simplest  instance of the general Bott periodicity theorem for the $K$-theory of Banach
algebras to be discussed in the next section.}
\end{example}

\begin{example} \rm{ For any Banach algebra $A$, we have  a surjection
$$ K_1^{alg} (A) \twoheadrightarrow K_1 (A).$$
Using \eqref{kho}, this follows from
$$ GL^{\circ} (A) = \overline{[GL(A), \, GL(A)]},$$
which we leave as an exercise.
 For $A=\mathbb{C}$, we have $K_1^{\text{alg}}(\mathbb{C})\simeq \mathbb{C}^{\times}$
with the isomorphism given by the determinant map $\text{det}: GL (\mathbb{C})
\to \mathbb{C}^{\times}$, while $K_1 (\mathbb{C})=0 $ (see the next example). }
\end{example}

\begin{example} \rm{
 (i) Since  $GL (n, \mathbb{C})$ is connected for all $n$, we have
  $K_1(\mathbb{C}) =0$.  Similarly, using polar decomposition, one shows
  that
  for any von Neumann algebra $A$,
   $GL (n, A)$ is connected for all $n$  and hence
 $K_1(A)=0$.\\
(ii) By Morita invariance we have $K_1(M_n(\mathbb{C})) =0$.\\
(iii) Since the algebra $\mathcal{K}$ of compact operators is the direct limit of finite matrices,
by continuity we have
$K_1(\mathcal{K})=0.$}
\end{example}

\begin{exercise} \rm {Starting from  the definitions, show that for $i=0, 1$
 $$ K_i(C(S^1))\simeq \mathbb{Z}.$$
  Under these isomorphisms a projection $e: S^1 \to M_n(\mathbb{C})$ is sent to $\text{tr} (e) \in \mathbb{Z}$ and
  an invertible $u: S^1 \to
  GL(n, \mathbb{C})$ is sent to the winding number of $ \text{det}\,
  (u(z))$.}
  \end{exercise}

The suspension functor is {\it exact} in the sense that for any exact sequence
of Banach algebras  $ 0 \to I \to A \to A /I \to 0$ the  sequence
$$0\longrightarrow SI\longrightarrow SA \longrightarrow S (A/I) \longrightarrow
0$$
is exact too. Coupled with the half exactness of $K_0$, we  conclude
that the sequences
$$ K_n(I)\longrightarrow K_n (A) \longrightarrow K_n (A/I) $$
are  exact in the middle for all $n\geq 0$.

One can  splice these half exact  sequences into a long exact sequence
\begin{equation} \label{ktes} \cdots   \to K_{n}(I) \to K_{n}(A) \to K_{n}(A/I) \to
\cdots   \to K_0(I) \to K_0(A) \to K_0(A/I)
\end{equation}
To do  this  it suffices to show that there exists a {\it
connecting homomorphism}
$$\partial : K_1(A) \longrightarrow K_0(I)$$
which renders the sequence
\begin{equation}\label{kes}
 K_{1}(I) \to K_1(A) \to K_{1}(A/I)  \overset{\partial}{\to} K_{0}(I) \to K_{0}(A) \to K_{0}(A/I)
 \end{equation}
exact.  It is sometimes called the {\it (generalized) index map}  since   for  the  Calkin extension
\begin{equation}\label{caex}
0\longrightarrow \mathcal{K}  \longrightarrow  \mathcal{L} (H)
\longrightarrow \mathcal{C} \longrightarrow 0
\end{equation}
 it coincides with the index of Fredholm operators
$$ \partial = \text{index}: K_1 (\mathcal{C}) \longrightarrow K_0 (\mathcal{K})
\simeq
\mathbb{Z}.$$

Let $u\in GL_n(A/I)$ and let $w\in GL_{2n}(A)$ be a
lift of $\text{diag}\, (u, u^{-1})$. Define
$$ \partial  ([u])=[wp_nw^{-1}]-[p_n] \in K_0(I),$$
where the projection $p_n = \text{diag} (I_n, \, 0).$
It can be shown that this map is well defined and \eqref{kes} is exact.

\begin{example} \rm{Let $A$ be a unital $C^*$-algebra. Using polar
decomposition one shows that any invertible in $GL (n, A)$ is
homotopic to a unitary. Given  such a unitary $u$, we can find a
partial isometry $v$ in $M_{2n} (A)$ lifting $ \text{diag} (u,
1)$. Now the unitary
$$ w= \left( \begin{array}{cc} v & 1-vv^*\\1-v^*v & v^* \end{array}
\right)$$ lifts $ \text{diag} \, ((u, 1), (u^{-1}, 1))$  and hence
$$\partial ([u]) =[wp_{2n}w^{-1}]-[p_{2n}]= [1-v^*v]-[1-vv^*].$$
 For
the Calkin extension \eqref{caex}  this maps sends an invertible in the
Calkin algebra $\mathcal{C}$, i.e. a Fredholm  operator, to its
Fredholm index in $K_0( \mathcal{K}) \simeq \mathbb{Z}$.}
\end{example}

\begin{remark} A  more conceptual way to get the long exact sequence \eqref{ktes} would be to derive it from the homotopy
exact sequence of a fibration.
\end{remark}

\subsection{ Bott periodicity theorem}

Homotopy invariance, Morita invariance, additivity and the  exact sequence \eqref{ktes} are  essential features of  topological
$K$-theory. The deepest result of $K$-theory, however, at least in the commutative case,   is
  the Bott periodicity
theorem.  It states that there is a natural isomorphism  between
$K_0$ and $K_2$. The isomorphism is given by
   the {\it Bott map}
$$\beta : K_0(A) \to K_1(SA). $$

Since $K_1 (SA) \simeq \pi_1 (GL (A))$ is the homotopy group of the
stable general linear group of $A$,  $\beta$ should somehow turn
an idempotent in $M (A)$ into a loop of invertibles in $GL (A) $.
We assume $A$ is unital (the general case easily follows).
Given an idempotent $e\in M_n(A)$, define a map $u_e : S^1 \to GL
(n, A)$ by $u_e(z)=ze+(1-z)e.$ It defines a loop in $GL (A)$ based
at 1, whose homotopy class is an element of $\pi_1 (GL (A)) \simeq
K_1 (SA)$. Now the {\it Bott map} $ \beta : K_0 (A) \longrightarrow
K_1 (SA)$ is defined by
 $$ \beta ([e]-[f])= u_e u_f^{-1}.$$
 Notice that, since $u_e$ is a group homomorphism the additivity of $\beta$ follows.

\begin{theorem} (Bott periodicity theorem) For a complex Banach
algebra $A$ the Bott map
$$ \beta : K_0(A) \to K_2(A)$$
is a natural isomorphism.
\end{theorem}

It  follows that  for all $n\geq 0$,
$$K_n(A)\simeq K_{n+2} (A)$$
and the long exact sequence \eqref{ktes} reduces to a periodic 6-term exact
sequence
$$\begin{CD}
K_0(I) @ >i_* >> K_0(A) @ >\pi_*>> K_0(A/I) \\
@A\partial AA   @.  @V V\partial V \\
K_1(A/I) @ <<\pi_*<  K_1(A) @ <<i_* <  K_1(I)
\end{CD}$$

\begin{example} \rm{ For $A= \mathbb{C}$ we already knew that $K_0
(\mathbb{C}) \simeq  \mathbb{Z}$ and $K_1
(\mathbb{C}) \simeq  0$. Using Bott periodicity we obtain $K_{2n} (\mathbb{C}) \simeq \mathbb{Z}$, and
$K_{2n+1} (\mathbb{C}) \simeq 0$. This
  is a non-trivial result and
 in fact of the same magnitude of difficulty as Bott periodicity. Since the spheres $S^n$ are iterated suspensions of a point,
 we obtain $ K^0 (S^{2n}) \simeq \mathbb{Z} \oplus \mathbb{Z}$,   $ K^1
 (S^{2n}) \simeq 0$, and
 $ K^0 (S^{2n+1}) \simeq K^1 (S^{2n+1})\simeq \mathbb{Z}.$  }
\end{example}

\begin{example} ( The Toeplitz extension) \rm{
The  Toeplitz extension
$$ 0\longrightarrow \mathcal{K} \longrightarrow
\mathcal{T}\overset{\sigma}{\longrightarrow} C(S^1)
\longrightarrow 0 $$
was introduced in  Section 3.1. Now the index map
$\partial : K_1 (C (S^1)) \to K_0( \mathcal{K})$  in the 6-term
exact sequence is an isomorphism (this is, more or less,
equivalent to the Gohberg-Krein index theorem for Toeplitz operators).
Since  $K_0 (\mathcal{K}) = \mathbb{Z}$,  $K_0 (\mathcal{K}) =0$,
and  $K_i ( C(S^1))= \mathbb{Z}$ for $i=0, \, 1$, from   the above 6-term exact sequence  we
deduce that $K_0(\mathcal{T}) \simeq K_1(\mathcal{T}) =\mathbb{Z}$.
 }
\end{example}

\subsection{Further results}
In the commutative case the properties of  homotopy invariance,
long exact sequence,  and Bott periodicity, suffice to
give a good understanding of the $K$-theory of spaces such as CW
complexes. Noncommutative spaces are of course much richer and
more complicated. Here we give
 two  further results  that
 have no counterpart in the commutative  case. Proofs of both can be found in \cite{bl1} or in the original
articles cited below.

\begin{theorem} (Connes' Thom isomorphism \cite{ac82}). If  $\alpha:
\mathbb{R} \to Aut (A)$ is a continuous one-parameter group of
automorphisms of a $C^*$-algebra $A$,  then
$$K_i(A\rtimes_{\alpha} \mathbb{R})\simeq K_{1-i}(A), \quad \quad i=0,
1.$$
\end{theorem}
In particular this result shows that  the $K$-theory of
$A\rtimes_{\alpha} \mathbb{R}$ is independent of the action $\alpha$.

The dimension shift is reminiscent of the dimension shift in the 
classical Thom isomorphism theorem relating the $K$-theory with
compact support of the total space of a vector bundle with the
$K$-theory of its base. Note that if the action is trivial then
the theorem reduces to the Bott periodicity theorem. In fact in this
case we have
$$ A\rtimes_{\alpha} \mathbb{R} \simeq A \otimes C_0(\mathbb{R})\simeq SA.$$

The second result we would like to highlight in this section is
the 6-term exact sequence of Pimsner and Voiculescu:
\begin{theorem} (Pimsner-Voiculescu exact sequence \cite{pv}). For any automorphism  $\alpha \in Aut( A)$
of a  $C^*$-algebra $A$ there is a 6-term exact
sequence\\

$$\begin{CD}
K_0(A) @ >1-\alpha_* >> K_0(A) @ >i_*>> K_0(A\rtimes_{\alpha}\mathbb{Z}) \\
@AAA   @.  @VVV\\
K_1(A\rtimes_{\alpha}\mathbb{Z}) @ <<i_*<  K_1(A) @ <<1-\alpha_* <  K_1(A)
\end{CD}$$\\
\end{theorem}

\begin{example} (K-theory of noncommutative tori) \rm{A beautiful
application of this result is to the $K$-theory of the
noncommutative torus. We have $A_{\theta} =C(S^1)\rtimes_{\alpha}
\mathbb{Z},$ where the automorphism $\alpha$ is through rotation by $ 2
\pi \theta$. But $\alpha$ is homotopic to the identity through
rotations by $2 \pi \theta t$, $t\in [0, \, 1]$. By homotopy
invariance of $K_i$, we obtain
$\alpha_{*}=1$. Using $K_i (C(S^1))\simeq \mathbb{Z}$ for $i=0, 1$
and a simple diagram chase we conclude from the Pimsner-Voiculescu
exact sequence that
$$K_i(A_{\theta})\simeq \mathbb{Z}\oplus \mathbb{Z} \quad \text{for}
\,\, \, i=0, 1.$$ Thus it seems that $K$-theory by itself cannot
distinguish the isomorphism class of $A_{\theta}$ for different
$\theta$. There is however an extra piece of structure in $K_0
(A)$, for $A$ a $C^*$-algebra,
 that can be used in this regard. Notice that $K_0(A)$ is an {\it ordered group}  with its positive cone
 defined by projections in $M (A)$. Equipped with this extra structure one can then show that $A_{\theta_1} \simeq A_{\theta_2}$
 iff $\theta_1 = \theta_2$ or $\theta_1=1-\theta_2$ (cf. \cite{bl1} and references therein).  }
\end{example}

\subsection{Twisted $K$-theory}
 Twisted $K$-theory has been around for quite some time (cf. Donovon-Karoubi \cite{doka}, Rosenberg \cite{ro1}). The
  recent
 surge of interest in the subject has to do with  both mathematics  and high energy physics. In mathematics, a recent
 result of Freed, Hopkins and Teleman shows that the twisted equivariant $K$-theory of a  compact Lie group is isomorphic
 to the Verlinde algebra (fusion algebra)  of the  group \cite{fr, fht}. The latter algebra  is  the
 algebra  of projective representations of the loop group of the group at a  fixed {\it level}.
 In some semiclassical limits of string theory over a background spacetime $X$,
 the strengths of $B$-fields are elements of $H^3 (X, \mathbb{Z})$. When this $B$-field is non-trivial the topological
 charges of D-branes are interpreted as elements of twisted $K$-theory with  respect to the twisting defined
 by $B$ (cf. \cite{bcmms, boma} for a mathematical perspective).   A recent comprehensive study of twisted $K$-theory  can be found
 in Atiyah and
 Segal's article \cite{atse}.

The twisting coefficients (local systems) of twisted $K$-theory
are cohomology  classes in $H^3 (X, \mathbb{Z})$. There are at
least two approaches to the subject. One can either extend the
definition of $K$-theory through Fredholm operators and the
relation \eqref{ajt} to include twistings as in \cite{atse}, or
one can define the twisted $K$-theory   as the $K$-theory of a
noncommutative algebra as is done in \cite{ro1}. We shall briefly
describe  this latter definition.

   Let $X$ be  a locally compact, Hausdorff,  and second countable space. We recall the classification of locally trivial bundles of
algebras with  fibers isomorphic to the  algebra $\mathcal{K} =\mathcal{K}(H)$ of
compact operators on an infinite dimensional Hilbert space, and with structure group  $ \text{Aut} \, (\mathcal{K})$.
As we saw in Section 2.4  there is a one  to one correspondence between
isomorphism classes of such bundles and $H^3 (X, \, \mathbb{Z})$. Given
 such a  bundle of algebras $\mathcal{A}$, its Dixmier-Douady invariant
 $$\delta (\mathcal{A}) \in H^3 (X, \mathbb{Z})$$ is a complete
 isomorphism invariant of such bundles.

Now given a pair $(X, \, \delta)$ as above, the {\it twisted
K-theory} of $X$ can be defined as the $K$-theory of  the $C^*$-
algebra $A =\Gamma (X, \,  \mathcal{A})$ of continuous  sections
of $\mathcal{A}$ vanishing at infinity:
$$ K^i_{\delta} (X):= K_i (A).$$

There is also an equivariant version of twisted $K$-theory,
denoted  by $K^{\delta}_{G} (X)$,   that is  specially important
in view of the recent work \cite{fr}. The coefficients for this
theory are elements of the equivariant cohomology $H^3_G (X, \,
\mathbb{Z})$, where $G$ is  a compact Lie group  acting  on a
space $X$. We refer to \cite{atse} for its definition. For
simplicity, let $G$ be a compact, connected,  simply connected,
and simple  Lie group. Then the central extensions of the loop
group $LG$ of $G$ are characterized by a positive integer $k$,
called the level \cite{prse}. For each positive integer $k$,  the
positive energy representations of this central extension, up to
equivalence, constitute  a finite set and  we denote by   $V_k (G)$  the
free abelian group generated by this set. There is a commutative
ring structure on this set, corresponding to tensor product of
representations.

  Now let  $G$ act on itself by conjugation ($G=X$). Then the    equivariant
cohomology $ H^3 (G, \, \mathbb{Z})$ is  a  free group of rank one
whose elements we shall denote by integers. The theorem of
Freed-Hopkins-Teleman  states that, at each level $k$, the fusion
ring $V_k (G)$
 is isomorphic to  the twisted
equivariant $K$-theory of $G$:
$$K^n_G (G) \simeq V_k (G),$$
where the integer $n$ can be explicitly defined in terms of $k\geq 0$
and $G$ \cite{fr, fht}.

\begin{example} \rm{ Let $G =SU (2) \simeq S^3$. Then $H^3 (G, \, \mathbb{Z})
\simeq \mathbb{Z}$. For each integer $n$ representing a class in
$H^3 (G, \, \mathbb{Z})$  there is a bundle $\mathcal{A}_n$ of
algebras of compact operators $\mathcal{K}$ over $S^3$ obtained by
gluing the trivial bundles $S^3_+ \times \mathcal{K}$  and $S^3_{-}
\times \mathcal{K}$ on the upper and lower hemispheres
respectively. The gluing is defined by  a map $S^2 \to \text{Aut} \,
(\mathcal{K})$ of degree $n$.  Let $A_n$ denote the $C^*$-algebra
of continuous sections of $\mathcal{A}_n$. We have then, by
definition,
$$  K^{i, n} (S^3)= K_i (A_n).$$

The representation ring of $G= SU (2)$ is a polynomial algebra
whose generator is the fundamental 2-dimensional representation of
$G$. The
 twisted equivariant $K$-theory of $SU (2)$   can be explicitly computed and
 shown to be a quotient of the representation  ring (cf.
 \cite{fr, fht}).}
\end{example}

\subsection{K-homology}
In \cite{at} and a little later and independently in \cite{bdf}, Atiyah
and Brown, Douglas and Fillmore proposed theories dual to
topological $K$-theory,  using techniques of functional analysis and
operator algebras. The cycles for Atiyah's theory are {\it abstract
elliptic operators} $ ( H, \, F)$ over $C(X)$ where $H=H^+ \oplus H^-$
is a $\mathbb{Z}_2$-graded Hilbert space, $ \pi: C(X) \to
\mathcal{L} (H)$ is an even representation of $C(X)$, and $F : H \to H$ is an even bounded operator
with $F^2-I \in \mathcal{K} (H)$.
 This data must satisfy the condition
$$[F, \, \pi (a)] \in \mathcal{K} (H).$$
We see that an abstract elliptic operator is the same as a
Fredholm module over $C(X)$ as in definition \eqref{psabl} except
that  now instead of $F^2=I$ we have the above condition. The two
definitions are however essentially equivalent. In particular the
formula
$$ \langle (H, \, F)\, , \, [e] \rangle = \text{index} \, F_e^+$$
where the Fredholm operator $F_e^+$ is defined in \eqref{fop} defines a
pairing between the $K$-theory of $X$ and abstract elliptic operators on $X$.

Let $\text{Ext}  (A)$ denote the set of isomorphism classes   of
extensions of a $C^*$-algebra $A$ of the form
$$0 \longrightarrow  \mathcal{K} \longrightarrow E \longrightarrow A
\longrightarrow 0.$$
There is a natural operation of addition of extensions which turn
$\text{Ext}  (A)$ into an abelian   monoid. It can be shown that if
$A$ is a nuclear $C^*$-algebra,
for example if $A=C(X)$ is commutative, then $\text{Ext} (A)$ is actually
a  group.  There is a pairing
$$ K_1 (A) \times  \text{Ext}  (A)  \to \mathbb{Z}$$
which is defined as follows.  Let $\delta$ be the connecting homomorphism in the 6-term exact sequence
of an extension $\mathcal{E}$ representing an element of $\text{Ext} \, (A)$, and let $[u] \in \,K_1\, (A)$. Then
$$\langle [u] \, , \, [\mathcal{E}]\rangle: = \delta  ( [u] )
\in K_0 \,  (\mathcal{K}) \simeq \mathbb{Z}.$$

\begin{example} \rm{ A simple example of a   non-trivial extension is the
Toeplitz extension
$$ 0 \longrightarrow  \mathcal{K} \longrightarrow \mathcal{T} \longrightarrow
C(S^1)\longrightarrow 0 $$ from Section 3.1. It can be shown
that   the class of this extension generates the $K$-homology
group $ \text{Ext}  (C (S^1))$ \cite{bdf}. A more elaborate
example  is the  pseudodifferential extension of the algebra of
functions on the cosphere bundle:
$$0\longrightarrow \mathcal{K}(L^2 (M)) \longrightarrow \Psi^0 (M)
\overset{\sigma}{\longrightarrow} C (S^* M)\longrightarrow 0.$$
briefly discussed in Section 3.1. For a comprehensive introduction
to $K$-homology see Higson and Roe's \cite{hr}.}
\end{example}

For noncommutative spaces, the above approach to $K$-homology works best when
the corresponding $C^*$-algebra is nuclear. For arbitrary
$C^*$-algebras, Kasparov's $KK$-theory provides a unified approach to both $K$-theory
and $K$-homology in  a bivariant theory (cf. \cite{kasp} and
\cite{bl1}).

\section{ Cyclic Cohomology}
Let $A$ be an algebra and $e, f \in M_n(A)$ be two idempotents.
How can we show that  $[e]\neq [f]$ in $K_0(A)$? Here is a simple
device that is often helpful in this regard. As we saw in Example
\eqref{trkt}, any trace $\tau : A \longrightarrow \mathbb{C}$
induces an additive map
$$\tau: K_0 (A)\longrightarrow \mathbb{C},$$
via the formula
$$\tau ([e]):=\sum_{i=1}^n e_{ii}.$$
Now  if $\tau  ([e]) \neq  \tau ([f])$, then, of course,  $[e] \neq [f]$.

 \begin{exercise} Let $A = C(X)$, where $X$ is compact and connected,  and let $\tau
(f)= f(x)$, for some fixed $x \in X$.   Show that if $E$ is a vector
bundle on $X$ and $e\in M_n(C(X))$ an idempotent defining $E$,
then $\tau ([e])= \text{dim}\, (E_x)$ where $E_x$ is the fibre of $E$ over
$x$. Thus $\tau$ simply measures the rank of the vector bundle.
   \end{exercise}

 The topological information hidden in an idempotent is much more
subtle than just its `rank' and in fact traces can only capture
zero-dimensional information. To know more about idempotents and
$K$-theory we need higher dimensional analogues of traces. They
are called {\it cyclic cocycles}  and their study is the subject
of {\it cyclic cohomology}. As we shall see in this section,
cyclic cohomology is the right noncommutative analogue of de Rham
homology of currents on smooth manifolds.

Cyclic cohomology was first discovered  by Alain Connes \cite{ac81, ac83, co85}. Let us first recall   a remarkable
subcomplex of the Hochschild complex called  the {\it Connes complex}  that was introduced by him for the definition of cyclic cohomology.
For an algebra $\A$ let
$$C^n(\A)= \text{Hom} \, ({\A}^{\otimes (n+1)}, \, \mathbb{C}), \quad n=0, 1, \cdots,$$
denote the space of $(n+1)$-linear functionals on $\A$ with values in
$\mathbb{C}$.
 The Hochschild differential $b: C^n(\A) \to C^{n+1}(\A)$ is defined by
\begin{eqnarray*}
(b \varphi)(a_0, \cdots ,a_{n+1})&=&\sum_{i=0}^n(-1)^i \varphi (a_0, \cdots,a_ia_{i+1},\cdots ,a_{n+1})\\
& &+(-1)^{n+1}\varphi (a_{n+1}a_0, a_1,\cdots ,a_{n}).
\end{eqnarray*}
One checks that $b^2=0.$ The cohomology of the complex
$(C (\A), \, b)$ is, by definition, the Hochschild cohomology of $\A$ (with
coefficients in the $A$-bimodule $A^*$)  and will be denoted by
$HH^n(\A)$. An $n$-cochain $ \varphi \in C^n (\A)$ is called {\it cyclic} if
$$\varphi (a_n, a_0, \cdots ,a_{n-1})=(-1)^n \varphi (a_0, a_1, \cdots, a_n)$$
for all $a_0, \cdots ,a_n$ in $\A$. Though it is not obvious at
all, one can check that cyclic cochains form a subcomplex
$$ (
C_{\lambda} (\A), \, b) \subset (C (\A), \, b)$$ of the Hochschild
complex.

We shall  refer to $( C_{\lambda} (\A), \, b)$ as the {\it
Connes complex} of $\A$. Its cohomology,  by definition,  is the
cyclic cohomology of $\A$ and will be denoted by $HC^n (\A)$. We
start our introduction to cyclic cohomology by some concrete
examples of cyclic cocycles \`{a} la Connes.

\subsection{Cyclic cocycles}
We give an example of a  cyclic cocycle. Let $M$ be a closed (i.e.
compact without boundary),  smooth, oriented, $n$-manifold. For
$f^0, \cdots, f^n \in \mathcal{A}=C^{\infty} (M)$, let
$$\varphi (f^0, \cdots , f^n)=\int_M f^0 df^1 \cdots df^n.$$
The $(n+1)$-linear cochain
$$ \varphi :  \mathcal{A}\times \cdots \times \mathcal{A} \to \mathbb{C}$$
has three properties: it is continuous with respect to the natural
Fr\'{e}chet space topology on $\mathcal{A}$; it is a Hochschild
cocycle;  and it  is
 cyclic.  For the cocycle condition, notice that
\begin{eqnarray*}
(b \varphi)(f^0, \cdots ,f^{n+1})&=&\sum_{i=0}^n(-1)^i \int_M
f^0 df^1\cdots d(f^if^{i+1}) \cdots  df^{n+1}\\
& &+(-1)^{n+1} \int_M f^{n+1}f^0 df^1 \cdots df^{n}\\
&=&0,
\end{eqnarray*}
where we only used   the  Leibniz rule for the de Rham
differential $d$ and the graded commutativity of the algebra
$(\Omega M, \, d)$ of differential forms on $M$. The  cyclic
property of $\varphi$
$$ \varphi (f^n, f^0, \cdots, f^{n-1})= (-1)^n \varphi (f^0,  \cdots,
f^{n})$$ is more interesting. In fact since
$$ \int_M ( f^ndf^0 \cdots df^{n-1}-(-1)^n f^0 d f^1 \cdots
df^{n})=\int_M d (f^n f^0df^1 \cdots df^{n-1}),$$
we see that the cyclic property of $\varphi$ follows from Stokes' formula
$$ \int_M d\omega =0,$$
which is valid for any $(n-1)$-form $\omega$ on a closed manifold  $M$.

A remarkable property of cyclic cocycles is that, unlike de Rham
cocycles which make sense only over commutative algebras, they can
be defined  over any noncommutative algebra and the resulting
cohomology theory is the right generalization of de Rham homology
of currents on a smooth manifold. Before developing cyclic
cohomology any further we give one more example. In the above
situation it is clear that if $V \subset M$ is a closed
$p$-dimensional oriented submanifold then the formula
$$\varphi (f^0, \cdots ,f^p) = \int_V f^0 df^1 \cdots df^p$$
defines a cyclic $p$-cocycle on $\A$. We can replace $V$ by {\it
closed  currents} on $M$ and obtain
 more cyclic cocycles.

Recall that a $p$-dimensional {\it current} $C$ on $M$ is a continuous linear functional
$C : \Omega^p M \to \mathbb{C}$ on the space
of $p$-forms on $M$. We write $\langle C, \, \omega \rangle $ instead of $C (\omega)$. For example a
zero dimensional current on $M$ is just a distribution on $M$. The differential of a current is
defined by  $ \langle dC, \,  \omega \rangle= \langle C, \, d \omega \rangle $ and in this way one  obtains the
complex of currents on $M$
whose homology is the
{\it de Rham homology}  of $M$.
\begin{exercise} Let $C$ be a $p$-dimensional current on $M$. Show that
the $(p+1)$-linear functional
$$ \varphi_C (f^0, \cdots,  f^p) =\langle  C \, , \, f^0 df^1 \cdots df^p \rangle$$
is a Hochschild cocycle on $\mathcal{A}$. Show that if $C$ is closed
then
$\varphi_C$ is a cyclic $p$-cocycle on $\A$.
\end{exercise}

Let $\mathcal{A}$ be an  algebra. Define the operators
$$   b': C^n(\A) \to C^{n+1} (\A), \quad \text{and} \quad  \lambda : C^n(\A) \to C^n(\A),$$
by
\begin{eqnarray*}
(b' \varphi )(a_0, \cdots ,a_{n+1})&=&\sum_{i=0}^n(-1)^i \varphi (a_0,
\cdots,a_ia_{i+1},\cdots ,a_{n+1}),\\
(\lambda \varphi )(a_0, \cdots ,a_n)&=&(-1)^n \varphi (a_n, a_0, \cdots ,a_{n-1}),\\
 \end{eqnarray*}
By a direct computation one checks that
\begin{equation} \label{blf}
(1-\lambda)b=b'(1-\lambda), \quad \quad b'^2=0.
\end{equation}
Notice that a cochain $\varphi \in C^n$ is cyclic if and only if
$(1-\lambda) \varphi =0.$ Using \eqref{blf} we obtain

\begin{lem} The space of cyclic cochains is invariant under  $b$, i.e. for all $n$,
$$b \,C^n_{\lambda} (\A) \subset C^{n+1}_{\lambda} (\A).$$
\end{lem}

We therefore have a subcomplex of the Hochschild complex, called the
{\it Connes complex} of $\A$:
\begin{equation} \label{ccom}
 C_{\lambda}^0 (\A)\overset{b}{\longrightarrow} C_{\lambda}^1 (\A)
\overset{b}{\longrightarrow}C^2_{\lambda}(\A) \overset{b}{\longrightarrow}\cdots
\end{equation}

The cohomology of this complex is called the {\it cyclic cohomology} of
$\A$ and will be denoted by $HC^n (\A)$, $n=0, 1, 2,\cdots.$
A cocycle for cyclic cohomology is called a {\it cyclic cocycle}. It
satisfies the two conditions:
$$ (1-\lambda )\varphi =0, \quad {\text and} \quad b \varphi =0.$$

\begin{examples}  \rm{ 1. Clearly $HC^0 (A)=HH^0(A) $ is the space of traces on $A$. In particular if $A$ is commutative then
 $HC^0(\A) \simeq A^*$ is the linear dual of $A$. \\

\noindent 2. Let $\mathcal{A} =C^{\infty} (M, M_n(\mathbb{C}))$ be the space of smooth matrix valued functions
on a closed smooth oriented manifold $M$. For any closed de Rham $p$-current on
$M$
$$\varphi_C (f^0, \cdots ,f^p)= \langle  C \, ,  \, Tr (f^0 df^1 \cdots df^p)\rangle$$
is a cyclic $p$-cocycle on $\mathcal{A}$. \\

\noindent 3. Let $\delta : \mathcal{A} \to \mathcal{A}$ be a
derivation and $\tau: \mathcal{A} \to \mathbb{C}$ an invariant
trace, i.e. $\tau (\delta (a))=0$ for all $a\in \mathcal{A}$. Then
one checks that
\begin{equation} \label{nctc1}
\varphi (a_0, a_1) =\tau( a_0 \delta(a_1))
\end{equation}
is a cyclic 1-cocycle on $\mathcal{A}$. This example can be generalized.
Let $\delta_1$ and $\delta_2$ be  a pair of {\it commuting} derivations which leave a trace $\tau$ invariant. Then
\begin{equation} \label{nctc2}
\varphi (a_0, a_1, a_2)=\tau (a_0 (\delta_1 (a_1) \delta_2 (a_2)-\delta_2
(a_1) \delta_1 (a_2)))
\end{equation}
is a cyclic 2-cocycle on $\mathcal{A}$.

Here is a concrete example
with $\A= \A_{\theta}$ a  smooth noncommutative torus.  Let
$\delta_1, \, \delta_2: \A_{\theta} \to \A_{\theta} $ be the unique
derivations defined by
$$\delta_1 (U)=U, \quad \delta_1 (V)=0; \quad
\delta_2 (U)=0, \quad \delta_2 (V)=V. $$
They commute with each other and
preserve  the standard trace $\tau$ on $\A_{\theta}$. The
resulting cyclic 1-cocycles  $\varphi_1 (a_0, a_1) =\tau( a_0
\delta_1(a_1))$ and  $\varphi_1' (a_0, a_1) =\tau( a_0
\delta_2(a_1))$  form a basis for the periodic cyclic cohomology
$HP^1 (\A_{\theta})$.  Similarly, the corresponding cocycle \eqref{nctc2} together with $\tau$ form a basis for
$HP^0 (\A_{\theta})$. }
\end{examples}

Consider the short exact sequence of complexes
$$ 0 \to C_{\lambda}\to C \to C /C_{\lambda} \to 0$$
Its associated long exact sequence is
\begin{equation} \label{ses1} \cdots \longrightarrow HC^n (\A) \longrightarrow HH^n (\A)
\longrightarrow H^n (C/C_{\lambda})
\longrightarrow  HC^{n+1}(\A) \longrightarrow \cdots
\end{equation}
We need to identify the cohomology groups $H^n(C/C_{\lambda})$. To this
end, consider  the short exact sequence
\begin{equation}\label{ses2} 0 \longrightarrow  C/C_{\lambda} \overset{1-\lambda}{\longrightarrow} (C, b')
\overset{N}{\longrightarrow}
C_{\lambda} \longrightarrow
 0 ,
 \end{equation}
 where the  operator $N$ is defined by
 $$ N=1+\lambda +\lambda^2 +\cdots +\lambda^n : C^n \longrightarrow
 C^n.$$  The relations
 $$ N (1-\lambda)=(1-\lambda) N =0, \quad \text{and} \quad bN=Nb'$$
 can be verified and they show that $1-\lambda$ and $N$ are morphisms of
 complexes in \eqref{ses2}.
\begin{exercise} Show that \eqref{ses2} is exact (the interesting part is to show
that $\text{Ker}\,  N \subset \text{Im}\, (1-\lambda)$).
\end{exercise}

 Now, assuming $A$ is unital,  the middle complex $(C, \, b')$  in  \eqref{ses2} can be shown to be exact with  contracting
 homotopy
 $ s: C^n \to C^{n-1}$
 given by
 $$ (s \varphi) (a_0, \cdots, a_{n-1})=(-1)^n \varphi (a_0, \cdots, a_{n-1}, 1).
 $$ It follows that
 $H^n(C/C_{\lambda})\simeq HC^{n-1} (\A)$. Using this in \eqref{ses1},  we obtain {\it Connes'
   long exact} sequence relating Hochschild  and cyclic cohomology:
  \begin{equation}\label{ibs}
   \cdots \longrightarrow  HC^n(\A) \overset{I}{\longrightarrow}HH^n(\A) \overset{B}{\longrightarrow}
  HC^{n-1}(\A) \overset{S}{\longrightarrow} HC^{n+1}(\A) \longrightarrow
  \cdots
  \end{equation}

The operators $B$ and $S$ can be made more explicit  by finding the connecting homomorphisms in
the above  long exact sequences. Remarkably, there is a formula for {\it Connes' boundary operator} $B$ on the level
of cochains  given by
$$ B =Ns (1-\lambda)= NB_0,$$
where $B_0: C^n \to C^{n-1}$ is defined by
$$ B_0 \varphi (a_0, \cdots, a_{n-1})=\varphi (1, a_0, \cdots, a_{n-1})
-(-1)^n \varphi (a_0, \cdots, a_{n-1}, 1).$$
The operator  $S: HC^n(\A) \to HC^{n+2}(\A)$  is called the {\it periodicity
operator} and is in fact related to Bott periodicity.   The {\it periodic
cyclic cohomology} of $\A$ is defined as the direct limit under the
operator $S$ of cyclic cohomology groups:
$$ HP^i (\A) = \underset{\longrightarrow}{\text{Lim}} \, HC^{2n+i} (\A), \quad \quad i=0, \, 1$$

A typical application of \eqref{ibs}  is to extract information
about cyclic cohomology from Hochschild cohomology.  For example,
assume $f: \mathcal{A} \to \mathcal{B}$ is an algebra homomorphism
such that $f^* : HH^n(\mathcal{B}) \to HH^n (\mathcal{A})$ is an
isomorphism for all $n\geq 0$. Then, using the five lemma, we
conclude that   $f^* : HC^n(\mathcal{B}) \to HC^n (\mathcal{A})$
is an isomorphism for all $n$. In particular from Morita
invariance of Hochschild cohomology one obtains  the Morita
invariance of cyclic cohomology.

\subsection{Connes' spectral sequence}

The cyclic complex \eqref{ccom} and the long exact sequence
\eqref{ibs}, as useful as  they are, are not powerful enough for
computations. A much deeper relation between Hochschild and cyclic
cohomology groups is encoded in Connes' spectral sequence that we
recall now. This spectral sequence resembles in many ways the
Hodge to de Rham spectral sequence for complex manifolds.  About
this connection we shall say nothing in these notes but see \cite{kal}
where a conjecture of Kontsevich and Soibelman  about the degeneration of this
spectral sequence is proved.

 Let $\mathcal{A}$ be a unital algebra.  Connes' $(b, B)$- bicomplex of
 $\mathcal{A}$ is the bicomplex

$$\begin{CD}
\vdots @.\vdots @.\vdots \\
C^2(\mathcal{A})@>B>> C^1(\mathcal{A})@>B>> C^0(\mathcal{A})\\
@AbAA @AbAA\\
C^1(\mathcal{A})@>B>>C^0(\mathcal{A})\\
@AbAA\\
C^0(\mathcal{A})
\end{CD}
$$
Of the three  relations
$$ b^2=0, \quad bB+Bb=0, \quad B^2=0,$$
only the middle relation is not obvious. But this follows from the 
relations $b's+sb'=1, \, (1-\lambda) b=b' (1-\lambda)$ and $Nb'
=bN$,  already used in this section.

\begin{theorem} (Connes \cite{co85}) The map
$ \varphi \mapsto  (0, \cdots, 0,\varphi)$ is a quasi-isomorphism
of complexes
$$(C_{\lambda} (\A),  \, b) \to (Tot \mathcal{B} \, (\A),  \, b+B)$$
\end{theorem}
This is a consequence of the  vanishing  of the $E^2$ term  of the
second spectral sequence (filtration by columns) of $\mathcal{B}
(A)$. To prove this consider the short exact sequence of
$b$-complexes
$$ 0 \longrightarrow \text{Im}\, B \longrightarrow \text{Ker}\,  B
\longrightarrow  \text{Ker}\,
B/ \text{Im}\, B \longrightarrow 0$$
By a hard lemma of Connes (\cite{co85},  Lemma 41 ), the induced map
$$ H_b (\text{Im} B) \longrightarrow  H_b (\text{Ker}\, B)$$
is an isomorphism. It follows that $H_b (\text{Ker}\, B / \text{Im}\, B)$ vanish. To
take care of the first column one appeals to the fact that
$$ \text{Im}\, B \simeq Ker (1-\lambda)$$
is the space of cyclic cochains.

\subsection{Topological algebras}
There is no difficulty in defining {\it continuous} analogues of
 Hochschild and cyclic cohomology  groups for Banach algebras. One simply replaces bimodules
 by Banach bimodules (where  the left and right module actions are
 bounded operators) and
 cochains
 by continuous cochains. Since the multiplication of a Banach algebra is a bounded map,  all
 operators including the Hochschild boundary and the cyclic operator
 extend to this continuous setting. The resulting Hochschild and cyclic
 theory for Banach and  $C^*$-algebras,  however, is hardly useful and tends  to vanish in
 many interesting
 examples. This is hardly surprising since the definition of any Hochschild and cyclic
 cocycle of dimension bigger than zero involves differentiating the
 elements of the algebra in one way or another. This is in sharp
 contrast with topological $K$-theory where the right   setting is the setting of Banach
 or $C^*$-algebras.

   \begin{exercise} Let $X$ be a compact Hausdorff space. Show that any derivation
   $\delta : C(X) \longrightarrow
   C(X)$ is identically zero.  (hint: first show that if
   $f=g^2$  and $g(x)=0$ for some $x \in X$, then $\delta (f)(x)=0$.)
   \end{exercise}

   \begin{remark} By results of Connes and Haagerup (cf. \cite{cob} and references therein),
   we know that a $C^*$-algebra
   is {\it amenable} if and only if it is {\it nuclear}. Amenability refers to the property that
   for all $n\geq 1$,
   $$ H^n_{cont}(A, M^*)=0,$$
   for any Banach dual bimodule $M^*$.
     In particular, by  using Connes' long exact
   sequence, we find that, for any nuclear $C^*$-algebra $A$,
    $$ HC^{2n}_{cont}(A)=A^*, \quad  \text{and} \quad
    HC^{2n+1}_{cont}(A)=0,$$
    for all $n\geq 0.$
    \end{remark}

The right class of topological algebras for cyclic cohomology
turns out to be the class of {\it locally convex algebras} \cite{co85}. An
algebra $\mathcal{A}$ equipped with a locally convex topology is
called a locally convex algebra if its multiplication map
$\mathcal{A} \times \mathcal{A} \to \mathcal{A}$ is jointly
continuous. Basic examples of locally convex algebras include the
algebra $\A = C^{\infty} (M)$ of smooth functions on a closed
manifold and the smooth noncommutative tori $ \A_{\theta}$ and
their higher dimensional analogues. The topology  of
$C^{\infty}(M)$ is defined by the sequence of seminorms
$$\|f\|_n = \text{sup} \, | \partial^{\alpha} \, f|;  \quad |\alpha | \leq n,$$
where the supremum is over a fixed, finite,  coordinate cover for
$M$ (see the exercise below for the topology of $\A_{\theta}$).

Given locally convex topological vector spaces $V_1$ and $V_2$,
their {\it projective tensor product} is a locally convex space
$V_1 \hat{\otimes}V_2 $ together with a  universal jointly
continuous bilinear map $V_1 \otimes V_2 \to V_1 \hat{\otimes}V_2
$ \cite{gro}. It follows from the universal property that for any
locally convex space $W$, we have a natural isomorphism between
continuous bilinear maps $V_1 \times V_2 \to W$ and continuous
linear maps $V_1 \hat{\otimes}V_2 \to W$. One  of the nice
properties of the projective tensor product is that for smooth
compact manifolds  $M$ and $N$, the natural map
  $$C^{\infty}(M)\hat{\otimes} C^{\infty}(N)\to C^{\infty}(M \times N)$$
  is an isomorphism.

A {\it topological left $\mathcal{A}$-module} is a locally convex
topological vector space $\mathcal{M}$  endowed with a continuous
left $\mathcal{A}$-module action $\mathcal{A} \times \mathcal{M}
\to \mathcal{M}$. A {\it topological free} left
$\mathcal{A}$-module  is a module of the type
$\mathcal{M}=\mathcal{A}\hat{\otimes} V$ where $V$ is a locally
convex space.  A {\it projective} module is a module which is a
direct summand in a free module.

Given a locally convex algebra $\A$, let
$$C^n_{\text{cont}} (\A) = \text{Hom}_{\text{cont}} (\A^{\hat{\otimes}n}, \,
\mathbb{C}) $$ be the space of continuous $(n+1)$-linear
functionals on $A$ and let $C^n_{\text{cont}, \, \lambda} (\A)$ denote
the space of continuous cyclic cochains on $\A$. All the algebraic
definitions and results of sections 6.1 and 6.2 extend to this
topological setting. In particular   one defines topological
Hochschild and cyclic cohomology groups of a locally convex
algebra. The right class of topological projective and free
resolutions are those resolutions that admit a continuous linear
splitting. This extra condition is needed when one  wants to prove
comparison theorems for resolutions. We won't go into details here
since this is very well explained in Connes' original article
\cite{co85}.

\begin{exercise} The sequence of norms
$$p_k(a)= \text{Sup}\,  \{ (1+|n|+|m|)^k |a_{mn}|\}$$
defines a locally convex topology on the smooth noncommutative
torus $\mathcal{A}_{\theta}$. Show that the multiplication of
$\mathcal{A}_{\theta}$ is continuous in this topology.
\end{exercise}

\subsection{The deformation complex}
What we called the Hochschild cohomology of $\mathcal{A}$ and
denoted by $HH^n( \mathcal{A})$ is in fact the Hochschild
cohomology of $\mathcal{A}$ with coefficients in the
$\mathcal{A}$-bimodule $\mathcal{A}^*$. In general, given  an
$\mathcal{A}$-bimodule $\mathcal{M}$, the Hochschild complex of
$\mathcal{A}$ with coefficients in the bimodule $\mathcal{M}$ is
the complex
$$ C^0 (\mathcal{A}, \,\mathcal{M}) \overset{\delta}{ \longrightarrow} C^1 (\mathcal{A}, \, \mathcal{M})
\overset{\delta}{ \longrightarrow} C^2 (\mathcal{A}, \, \mathcal{M}) \longrightarrow
\cdots $$ where $C^0 (\mathcal{A}, \, \mathcal{M})=\mathcal{M}$
and $C^n (\mathcal{A}, \, \mathcal{M})=\text{Hom}_{\mathbb{C}}\,
(A^{\otimes n}, \, \mathcal{M})$ is the space of $n$-linear
functionals on $\mathcal{A}$ with values in $\mathcal{M}$. The
differential $\delta $ is given by
\begin{eqnarray*}
(\delta \varphi)(a_1, \cdots, a_{n+1})&=& a_1 \varphi (a_2, \cdots, a_{n+1})+\sum_{i=1}^n (-1)^{i+1}
\varphi (a_1, \cdots, a_ia_{i+1},\cdots , a_{n+1})\\
& & + (-1)^{n+1} \varphi (a_1, \cdots, a_{n})a_{n+1}.
\end{eqnarray*}

Two special cases are particularly important. For  $\mathcal{M}=\mathcal{A}^*$, the linear dual of $\A$
  with the bimodule action
$$ (afb)(c)=f(bca)$$
for all $a, b, c$ in $\mathcal{A}$ and $f\in \mathcal{A}^*$, we
obtain the Hochschild groups $H^n (\A, \A^*) = HH^n (A)$. This is
important in cyclic cohomology since as we saw it enters into a
long exact sequence with cyclic groups. The second important case
is when $\mathcal{M} =\mathcal{A}$ with bimodule structure given
by  left and right multiplication. The resulting complex $ (C (\A,
\, \A), \, \delta)$ is called the {\it deformation complex} of
$\A$. It is the complex that underlies the deformation theory of
associative algebras as studied by Gerstenhaber \cite{ger}.

There is a  much deeper structure hidden in the deformation complex
$C(\mathcal{A}, \, \mathcal{A}), \, \delta)$
 than meets the eye and we will only barely scratch  the
 surface. The first piece of structure is the cup product. The {\it cup product}
 $ \cup : C^p \times C^q \to C^{p+q}$ is defined by
 $$ (f \cup g ) (a^1,  \cdots, a^{p+q})= f(a^1, \cdots, a^p) g (a^{p+1},\cdots , a^{p+q}).$$
 Notice that $\cup$ is associative and one  checks that this product
 is compatible with the differential
 $\delta$ and hence induces an associative graded product on
 $ H(\A, \, \A) =\oplus H^n (\A, \, \A)$. What is not so obvious however is that
 this product is graded commutative for any $\A$ \cite{ger}.

 The second piece of structure on $  (C(\mathcal{A}, \, \mathcal{A}),\,
 \delta)$ is a graded Lie bracket. It is based on
the Gerstenhaber circle product $ \circ : C^p \times C^q
\longrightarrow C^{p+q -1}$  defined by
 $$ (f\circ g) (a_1, \cdots, a_{p+q-1})=\sum_{i=1}^{p-1}
 (-1)^{|g|(|f|+i-1)} f (a^1, \cdots, g(a^i, \cdots, a^{i+p}), \cdots,
 a^{p+q-1}).$$
 Notice that $\circ$ is not an associative product. Nevertheless one can
 show that \cite{ger} the corresponding graded bracket $[ \, , \, ] : C^p \times C^q \to
 C^{p+q-1}$
 $$ [f, \, g] = f\circ g -(-1)^{(p-1)(q-1)} g\circ f$$
  defines a graded Lie algebra structure on deformation cohomology
  $ H (\A,\,  \A)$. Notice  that
  the Lie algebra grading  is  now shifted by one.

  What is most interesting is that the cup product and the Lie algebra structure
  are compatible in the sense that $[ \, , \, ]$ is a graded derivation
  for the 
  cup product;  or in short $ (H (\A, \, \A), \, \cup, \, [ \,  , \, ])$
  is a graded Poisson algebra.

  The {\it fine structure} of the Hochschild cochain complex
  $(C(\mathcal{A}, \, \mathcal{A}), \, \delta)$, e.g. the existence of higher
  order products and homotopies between them is the subject of  many
  studies in recent years \cite{ko,ko1,ko2}.
   While it is relatively easy to write down these higher order products in the form of
  a brace algebra structure on the Hochschild complex,  relating
  them to  known geometric structures such as moduli of curves is quite hard.

\begin{remark} The graded Poisson algebra structure on  deformation
cohomology $ H (\A, \, \A)$  poses a natural question:  is $ H (\A, \,
\A)$ the semiclassical limit of a quantum cohomology theory for
algebras?
\end{remark}

 \begin{example} \label{hhofm}\rm{ Let $\mathcal{A} = C^{\infty}(M)$, where $M$ is a compact $n$-dimensional
 manifold. In \cite{co85} Connes
gives a projective resolution of the topological  left $\mathcal{A}\hat{\otimes}\mathcal{A}$-module
$\mathcal{A}$,
\begin{equation}\label{cres}
 \mathcal{A}\leftarrow \mathcal{M}_0 \leftarrow \mathcal{M}_1 \cdots
\leftarrow \mathcal{M}_n \leftarrow 0
\end{equation}
 where $\mathcal{M}_i$ is
the space of  smooth sections of the vector bundle $p_1^*
(\bigwedge^i TM)$  and $p_1 : M \times M \to M$ is the
projection on the second factor.
After applying the $ \text{Hom}_{\A \hat{\otimes}\A} (-, \, \A)$ functor to  \eqref{cres},
one obtains a complex with
zero differentials, which shows that
$$ H^p (\A, \, \A) \simeq C^{\infty} ({\bigwedge}^p TM), \quad p=0, 1, \cdots.$$
The latter is the space of polyvector fields on $M$.}
\end{example}

\subsection{ Cyclic homology}

Cyclic cohomology is a contravariant functor on the category of
algebras. There is a dual covariant theory   called {\it cyclic
homology} that we introduce now. The relation between the two is
similar to the relation between currents and differential forms on
manifolds.

For each $n\geq 0$,  let
$C_n(\mathcal{A})= \mathcal{A}^{\otimes
(n+1)}$.  Define the operators
\begin{eqnarray*}
b: C_n(\mathcal{A}) &\longrightarrow& C_{n-1}(\mathcal{A})\\
b': C_n(\mathcal{A}) &\longrightarrow&
 C_{n-1}(\mathcal{A})\\
  \lambda: C_n(\mathcal{A}) &\longrightarrow& C_{n}(\mathcal{A})\\
    s: C_n(\mathcal{A})&\longrightarrow&  C_{n+1}(\mathcal{A})\\
 N:  C_n(\mathcal{A}) &\longrightarrow& C_{n}(\mathcal{A})\\
  B: C_n(\mathcal{A})& \longrightarrow & C_{n+1}(\mathcal{A})
\end{eqnarray*}
by
\begin{eqnarray*}
b(a_0\otimes \cdots \otimes a_n) & =& \sum_{i=0}^{n-1} (-1)^i( a_0\otimes
\cdots \otimes a_i a_{i+1} \otimes \cdots \otimes a_n)\\
& +& (-1)^n(a_na_0\otimes a_1\cdots \otimes a_{n-1})\\
b'(a_0\otimes \cdots \otimes a_n) & =& \sum_{i=0}^{n-1} (-1)^i ( a_0\otimes
\cdots \otimes a_i a_{i+1} \otimes \cdots \otimes a_n)\\
\lambda (a_0\otimes \cdots \otimes a_n) & =& (-1)^n( a_n\otimes
a_{0}
\cdots  \otimes a_{n-1})\\
s(a_0\otimes \cdots \otimes a_n)&=&(-1)^n (a_0\otimes \cdots \otimes
a_n\otimes 1)\\
N &=& 1+\lambda +\lambda^2 +\cdots + \lambda^n\\
B &=& (1-\lambda) s N
\end{eqnarray*}
They satisfy the relations
\begin{eqnarray*}
 b^2=0, \quad  \quad b'^2=0, \quad  \quad (1-\lambda
 )b'=b(1-\lambda)\\
b'N= Nb, \quad  \quad B^2=0, \quad \quad  bB+Bb=0
\end{eqnarray*}

The complex  $(C_{\bullet}(\mathcal{A}), \, b)$ is  the
Hochschild complex of $\mathcal{A}$ with coefficients in the $\mathcal{A}$-bimodule
$\mathcal{A}$.   The complex
$$ C_n^{\lambda}(\mathcal{A}):= C_n(\mathcal{A})/\text{Im} (1-\lambda)$$
is called the
{\it Connes complex} of $\mathcal{A}$ for cyclic homology. Its homology, denoted by
$HC_n (\mathcal{A}), \, \, n=0, 1, \cdots$,  is
called the {\it cyclic homology} of $\mathcal{A}$. It is clear that the space of
cyclic cochains is the linear dual of the space of cyclic chains
$$C_{\lambda}^n (\A)\simeq \text{Hom}\, (C_n^{\lambda}(\A), \, \mathbb{C})$$
and
$$HC^n (\A)\simeq HC_n (\A)^*. $$

Similar to cyclic cohomology, there is a long exact sequence relating Hochschild and cyclic homologies, and also there is
a spectral sequence from Hochschild to cyclic homology. In particular
cyclic homology can be computed using  the following bicomplex.

$$\begin{CD}
\vdots @.\vdots @.\vdots \\
A^{\otimes 3}@<B<< A^{\otimes 2}@<B<< A \\
@VVbV @VVbV\\
A^{\otimes 2}@<B<<A\\
@VVbV\\
A
\end{CD}
$$

\begin{example} (Hochschild-Kostant-Rosenberg and Connes theorems) \rm{Let
$$\mathcal{A}\overset{d}{\longrightarrow}\Omega^1 \mathcal{A}
 \overset{d}{\longrightarrow}\Omega^2
 \mathcal{A}\overset{d}{\longrightarrow}\cdots$$
 denote  the de Rham complex of
   a commutative unital algebra $\mathcal{A}$. By definition  $d:
 \mathcal{A}\to \Omega^1 \mathcal{A}$ is a {\it universal derivation} into a  symmetric
 $\mathcal{A}$-bimodule
 and $\Omega^n \A: =\wedge^n_{\A} \Omega^1 \A$ is the
 $k$-th exterior power of $\Omega^1 \mathcal{A}$ over $\mathcal{A}$. One
usually defines $\Omega^1 \A$,  the {\it module of K\"{a}hler differentials}, as $I/I^2$,  where $I$ is the kernel of the
multiplication map $A\otimes A \rightarrow A$.
  $d$ is then defined
by
$$d(a)= a\otimes 1-1\otimes
a \quad \text{mod} (I^2).$$
The universal derivation $d$ has a unique  extension to a
graded derivation of degree one on $\Omega \mathcal{A}$, denoted by $d$.

  The {\it antisymmetrization map}
 $$\varepsilon_n:  \Omega^n \mathcal{A} \longrightarrow
 \mathcal{A}^{\otimes (n+1)}, \quad n=0, 1, 2, \cdots,$$
 is defined by
 $$\varepsilon_n ( a_0 da_1\wedge  \cdots \wedge da_n)=\sum_{\sigma \in S_n} sgn (\sigma)
 a_0\otimes a_{\sigma (1)}\otimes \cdots \otimes a_{\sigma (n)},$$
 where $S_n$ is  the symmetric group of order  $n$ .
We also have a  map
 $$\mu_n : \mathcal{A}^{\otimes n} \longrightarrow \Omega^n \mathcal{A}, \quad n=0, 1, \cdots$$
 $$\mu_n (a_0 \otimes a_1\otimes \cdots \otimes a_n)= a_0 da_1 \wedge \cdots \wedge da_n.$$
 One checks that the resulting  maps
 $$(  \Omega  \mathcal{A}, \, 0) \to ( C (\A), \, b), \quad
 \text{and}, \quad    ( C (\A), \, b)   \to (  \Omega  \mathcal{A}, \, 0)$$
                are morphisms of complexes, i.e.
 $$ b \circ \varepsilon_n =0,  \quad \mu_n \circ b =0.$$
 Moreover, one can easily check that
 $$\mu_n \circ \varepsilon_n =n! \, \text{Id}_n.$$
 It follows that, for any commutative  algebra $\A$,   the antisymmetrization map induces  an inclusion
 $$ \varepsilon_n: \Omega^n \mathcal{A} \hookrightarrow HH_n(\mathcal{A}),$$
 for all $n$.

 The  Hochschild-Kostant-Rosenberg  theorem \cite{hkr} states
 that if $\mathcal{A}$ is a {\it regular  algebra}, e.g.    the algebra of regular functions
 on a smooth affine variety, then $\varepsilon_n$ defines an
 algebra isomorphism
$$\varepsilon_n:  \Omega^n \mathcal{A} \simeq HH_n(\mathcal{A})$$
 between  Hochschild homology of $\A$ and  the algebra of  differential forms
  on $\A$.

To compute the   cyclic homology of $\A$, we first show that under the map $\mu$ the
operator $B$ corresponds to the de Rham differential $d$. More
precisely, for each integer $n\geq 0$ we have a commutative diagram:
$$\begin{CD}
C_{n}(\A)@> \mu >> \Omega^{n} \A \\
@VVBV @VVdV\\
C_{n+1}(\A)@>\mu >> \Omega^{n+1} \A \\
\end{CD}
$$\\

We have
\begin{eqnarray*}
\mu B (f_0\otimes \cdots \otimes f_n)&=&\mu \sum_{i=0}^n (-1)^{ni}(1\otimes f_i \otimes
\cdots \otimes f_{i-1}-(-1)^n f_i \otimes
\cdots f_{i-1} \otimes 1)\\
&= & \frac{1}{(n+1)!}\sum_{i=0}^n (-1)^{ni}df_i \cdots df_{i-1}\\
&= &\frac{1}{(n+1)!}(n+1)df_0 \cdots df_n\\
&=& d \mu (f_0\otimes \cdots \otimes f_n).
\end{eqnarray*}

It follows that $\mu$ defines a morphism of bicomplexes
$$ \mathcal{B}(\A) \longrightarrow \Omega (\A),$$
where $\Omega (\A)$ is the bicomplex
$$\begin{CD}
\vdots @.\vdots @.\vdots \\
\Omega^2 \A @<d<< \Omega^1 \A@<d<< \Omega^0 \A \\
@VV0V @VV0V\\
\Omega^1 \A@<d<<\Omega^0 \A\\
@VV0V\\
\Omega^0 \A
\end{CD}
$$

 Since $\mu$ induces isomorphisms on row homologies, it induces isomorphisms on
 total homologies as well. Thus we have \cite{co85, lod}:
 $$HC_n( \A)\simeq \Omega^n \A/\text{Im}\, d \oplus H^{n-2}_{dR}(\A) \oplus \cdots \oplus H^k_{dR}
 (\A),$$
 where $k=0$ if $n$ is even and $k=1$ if $n$ is odd.

Using the same map $\mu$ acting between the corresponding periodic
complexes, one concludes that the periodic cyclic homology of $\A$
is given by
$$HP_k (\A)\simeq \bigoplus_i H^{2i+k}_{dR}(\A), \quad k=0, 1.$$

By a completely similar method one can compute the {\it continuous cyclic homology} of
the algebra $\A = C^{\infty} (M)$ of smooth
functions on a smooth closed manifold $M$.  Here by continuous cyclic homology we
mean the homology of the cyclic complex
where instead of algebraic tensor products $ \A \otimes \cdots \otimes
\A$, one uses the topological {\it projective tensor product} $  \A
\hat{\otimes} \cdots \hat{\otimes}
\A$.  The  continuous Hochschild homology of $\A$ can be computed using
Connes' topological resolution for $\A$ as an $\A$-bimodule as in  Example \eqref{hhofm}. The result is
$$ HH_n^{cont} (C^{\infty} (M)) \simeq \Omega^n M$$
with  isomorphism induced by the  map
$$ f_0 \otimes f_1 \otimes \cdots \otimes f_n \mapsto f_0 df_1 \cdots
df_n. $$ The rest of the computation of continuous cyclic homology
follows the same pattern as in the case of regular algebras above. The end
result is \cite{co85}:
$$HC_n^{\text{cont}} ( C^{\infty} (M))\simeq \Omega^n M/\text{Im}\, d \oplus H^{n-2}_{dR}( M) \oplus \cdots \oplus H^k_{dR}
 (M),$$
 and
 $$HP_k^{cont} (C^{\infty} (M))\simeq \bigoplus_i H^{2i+k}_{dR}(M ), \quad k=0, 1.$$

 The cyclic (co)homology of (topological) algebras is computed in many cases. We refer to \cite{co85} for smooth  noncommutative
 tori, to \cite{bur} for group algebras, and to \cite{ft, nis} for
 crossed product algebras. In \cite{ak} the spectral sequence for group
 crossed products has been extended to Hopf algebra crossed products. We
 refer to \cite{cu, cq1, cq3} 
 for an alternative approach to cyclic (co)homoloy due to Cuntz and
 Quillen.

}
\end{example}

\subsection{ Connes-Chern character}
 We can now indicate Connes'
generalization of the pairing  between
$K_0(\mathcal{A})$ and traces on $\mathcal{A}$ (Example \eqref{trkt})  to a full fledged
pairing between $K$-theory and cyclic cohomology:
\begin{eqnarray}\label{p1}
 K_0 (\mathcal{A}) \times HC^{2n}(\mathcal{A}) &\longrightarrow &
\mathbb{C},\\
K_1^{\text{alg}}(\mathcal{A}) \times HC^{2n+1}(\mathcal{A}) &\longrightarrow &
\mathbb{C},
\end{eqnarray}
These maps are defined for all $n\geq 0$ and are compatible with the $S$-operation on
cyclic cohomology. In the dual setting  of cyclic homology, these
pairings translate into noncommutative {\it Connes-Chern characters}
\begin{eqnarray*}
 \text{Ch}_0^{2n}: K_0 (\mathcal{A})  &\longrightarrow & HC_{2n}(\mathcal{A}),\\
\text{Ch}_1^{2n+1}: K_1^{\text{alg}} (\mathcal{A})  &\longrightarrow & HC_{2n+1}(\mathcal{A}),
\end{eqnarray*}
compatible with the $S$-operation on cyclic homology. As a
consequence of compatibility with the periodicity operator $S$,
we obtain maps
\begin{eqnarray*}
 \text{Ch}_0: K_0 (\mathcal{A})  &\longrightarrow & HP_{0}(\mathcal{A}),\\
\text{Ch}_1: K_1^{\text{alg}} (\mathcal{A})  &\longrightarrow & HP_1(\mathcal{A}).
\end{eqnarray*}
For $\mathcal{A}= C^{\infty}(M)$, these maps reduce to the classical
Chern character as defined via the connection and curvature formalism
of Chern-Weil theory \cite{ms}.

The  definition of these pairings rest on the following three facts \cite{co85,cob, lod}:\\

\noindent 1) For any $k\geq 1,$ the map
$ \varphi \mapsto  \varphi_k$   from  $C^n (\mathcal{A}) \to C^n (M_k(\mathcal{A}))$ defined by
$$ \varphi_k (m_0\otimes a_0, m_1 \otimes a_1, \cdots , m_n\otimes
a_n)=\text{Tr}\, (m_0m_1 \cdots m_n)\varphi ( a_0,  a_1, \cdots ,
a_n)$$
commutes with the operators $b$ and  $\lambda$. It follows that if $\varphi$ is
 a cyclic cocycle on $\A$, then $\varphi_k$ is  a cyclic  cocycle on
 $M_k(\mathcal{A})$. \\

\noindent  2) Inner automorphisms act  by the identity on Hochschild
and hence on
 cyclic cohomology.  \\

 \noindent 3) (normalization) The inclusion of {\it normalized cochains} $C_{\lambda}^{\text{norm}}(\mathcal{A})
 \to  C_{\lambda}(\mathcal{A})$ is a quasi-isomorphism in dimensions
 $n\geq 1$.
 A  cyclic cochain $\varphi$ is called {\it normalized} if
 $\varphi (a_0, \cdots, a_n)=0$ if $a_i=1$ for some $i$. \\

Now let $e \in
M_k(A)$ be an idempotent and $[\varphi] \in HC^{2n}(\A)$. The pairing \eqref{p1} is defined by  the  bilinear map
$$\langle [\varphi], \, [e]\rangle = {\varphi}_k  (e, \cdots,e).$$
Let us first check that the value of the pairing depends only on the cyclic cohomology class of
$\varphi$. It suffices to assume $k=1$ (why?). Let
$\varphi =b\psi$ with $\psi \in C_{\lambda}^{2n-1}(A)$.  Then we have
\begin{eqnarray*}
\varphi (e, \cdots,e)&=&b\psi (e, \cdots, e)\\
&=&\psi (ee, e,
\cdots,e)
-\psi (e, ee, \cdots, e)+ \cdots +(-1)^{2n} \psi(ee, e,
\cdots, e)\\
&=&\psi (e, \cdots,e)\\
&=&0,
\end{eqnarray*}
where the last relation follows from the cyclic property of
$\psi$. The pairing is clearly invariant under the inclusion
$M_k(\mathcal{A}) \to M_{k+1}(\mathcal{A})$.

It remains to show that  the value of $ \langle [\varphi]\, , \, [e]\rangle $, for fixed
$\varphi$, only depends on the class of $[e] \in K_0(\A)$. It suffices to check that
 for $u\in GL_k( \A)$,
$\langle [\varphi], \, [e] \rangle = \langle [\varphi], \, [ueu^{-1}]\rangle $.  But this is
exactly fact 2) above.

\begin{exercise} Let $e\in \mathcal{A}$ be an idempotent. Show that
$$Ch_0^{2k} (e):= (-1)^k \frac{(2k)!}{k!} \text{tr}\,  ((e-\frac{1}{2})\otimes e^{\otimes (2k+1)}), \quad k=0, \cdots,
n,$$ defines a cycle in the $(b, B)$-bicomplex of $\mathcal{A}$.
This is the formula for the Connes-Chern character in the
bicomplex picture of cyclic homology.

Dually, given   a cocycle  $\varphi =(\varphi_0, \varphi_2,
\cdots, \varphi_2n)$ in the $(b, B)$-bicomplex,  its pairing with
an idempotent $e\in M_k (\mathcal{A})$ is given by
$$ \langle  [\varphi], \, \, [e] \rangle = \sum_{k=1}^{n} (-1)^k \frac{k!}{(2k)!}\varphi_{2k}(e-\frac{1}{2}, e, \cdots, e).$$
\end{exercise}

Given a normalized cyclic cocycle  $\varphi \in HC^{2n+1} (\A)$ and an
invertible $u \in GL (k, A)$, let
$$ \langle [\varphi], \, [u] \rangle = \varphi_k (u, u^{-1}, \cdots, u, u^{-1}).$$
It can be shown (cf. \cite{co85}, Part II) that the above formula defines a pairing between $K_1^{\text{alg}}$
and $HC^{2n+1}$.

\begin{exercise}
Given an invertible $u \in GL (n, \A)$, show that
$$ Ch_1^{2k+1}: =(-1)^k k! \, Tr (u^{-1} \otimes u)^{\otimes 2k}$$
defines a cycle in the {\it normalized} $(b, B)$-bicomplex of
$\mathcal{A}$. This is the formula for the Connes-Chern character in the $(b, B)$-bicomplex picture
of cyclic homology. Dually, given a normalized  $(b, B)$-cocycle
$\varphi =(\varphi_1, \cdots,
\varphi_{2n+1})$, the formula
$$ \langle [\varphi], \, \, [u] \rangle  = \sum_{k=1}^n (-1)^k \varphi_{2k+1} (u, u^{-1},
\cdots, u, u^{-1})$$
defines the pairing between $K_1^{\text{alg}}$ and $HC^{2n+1}$.
\end{exercise}

It often happens that an element of $ K_0 (\mathcal{A})$ is represented by
a finite projective module and not  by an explicit idempotent.  It is then important
to have a formalism that would give the value of its pairing with cyclic cocycles.  This is
based on a noncommutative version of Chern-Weil theory developed by Connes in \cite{ac80, co85} that  we sketch
next.

Let $\mathcal{E}$ be a finite projective right $\mathcal{A}$-module,  $(\Omega, \, d)$  a
differential
calculus on  $\mathcal{A}$ and let $\int : \Omega^{2n} \to \mathbb{C}$ be a
closed graded trace representing a cyclic cocycle $\varphi$ on
$\mathcal{A}$. Thanks to its  projectivity,   $\mathcal{E}$ admits a
{\it connection}, i.e. a degree one map
$$ \nabla : \mathcal{E}\otimes_{\A} \Omega \to \mathcal{E}\otimes_{\mathcal{A}} \Omega$$
which satisfies the  graded Leibniz rule
$$\nabla (\xi \omega) = \nabla (\xi) \omega + (-1)^{\text{deg}\, \xi}\xi
d\omega$$ with respect to the right $\Omega$-module structure on $\mathcal{E}\otimes_{\mathcal{A}} \Omega$. The
{\it curvature}  of $\nabla$ is the operator    $\nabla^2 $,  which
can be
 easily checked to be $\Omega$-linear,
$$\nabla^2 \in \text{End}_{\Omega} (E \otimes_{\mathcal{A}} \Omega )= \text{End}_{\mathcal{A}}  (\mathcal{E})\otimes \Omega.$$
Now since $\mathcal{E}$ is finite projective over $\A$ it follows
that $\mathcal{E}\otimes_{\A} \Omega$ is finite projective over
$\Omega$ and therefore the trace $\int : \Omega \to \mathbb{C}$
extends to a trace, denoted again by $\int$,  on $
\text{End}_{\mathcal{A}}  (\mathcal{E})\otimes \Omega$ (cf.
formula \eqref{trk}). The following result of Connes relates the
value of the pairing as defined above to its value  computed
through the Chern-Weil formalism:
$$ \langle [\mathcal{E}], \, [\varphi] \rangle  = \frac{1}{n!} \int \nabla^n.$$

\begin{example}
\rm{Let $\mathcal{S}(\mathbb{R})$ denote the Schwartz space of rapidly
decreasing functions on the real line. The operators
$u, v : \mathcal{S}(\mathbb{R}) \to \mathcal{S}(\mathbb{R})$ defined
by
$$ (u f)(x)= f(x-\theta), \quad \quad (v f)(x)= e^{2\pi i \theta}
f(x)$$ satisfy the relation $u v =e^{2\pi i \theta}vu$
and hence turn $\mathcal{S}(\mathbb{R})$ into   a right
$\mathcal{A}_{\theta}$-module via the maps $U \mapsto u, \,
V\mapsto v$. We denote this module by $\mathcal{E}_{0, 1}$. It is the
simplest  of a series of modules $\mathcal{E}_{p, q}$ defined by Connes in \cite{ac80}.
It turns out that $\mathcal{E}_{0, 1}$  is finite  projective, and
for the canonical trace $\tau $ on $\mathcal{E}_{0, 1}$ we have
$$ \langle \tau, \,  \mathcal{E}_{0, 1} \rangle = -\theta.$$
Using the two derivations $\delta_1, \delta_2$ as a basis for
``invariant vector fields" one can define a differential calculus
$\Omega^0 \oplus \Omega^1 \oplus \Omega^2$ on
$\mathcal{A}_{\theta}$ with $\Omega^i =\bigwedge^i\{de_1, de_2\}
\otimes \mathcal{A}_{\theta}$. A connection on $\mathcal{E}$ with
respect to this calculus  is a pair of operators
$\nabla_1, \, \,  \nabla_2: \mathcal{E} \to \mathcal{E}$
satisfying
$$ \nabla_j (\xi a)= (\nabla_j) (\xi a) + \xi \delta_j (a)$$
for all $\xi \in \mathcal{E}_{0, 1}$ and $a \in \mathcal{A}_{\theta}$
and $j=1, 2$. One can check that the following formula defines a
connection on  $\mathcal{E}_{0, 1}$ \cite{ac80, cob}:
$$\nabla_1 (\xi) (s)= -\frac{s}{\theta}\,  \xi (s), \quad \quad \nabla_2
(\xi) (s)=  \frac{d \xi}{ds} (s).$$
The curvature of this connection is constant and is given by $\nabla^2 =[\nabla_1, \nabla_2]= \frac{1}{\theta} I \in
\text{End}_{\mathcal{A}}(\mathcal{E}_{0, 1})$. }

\end{example}

\begin{remark} Chern-Weil theory is a theory of characteristic classes
for smooth  principal $G$-bundles, where $G$ is a Lie group. The
above theory  for  noncommutative vector bundles should be
generalized to noncommutative analogues of principal bundles. A
good point to start would be  the theory of Hopf-Galois extensions.
\end{remark}

In the remainder of this section we shall briefly introduce  Connes' Chern
character in $K$-homology (cf. \cite{co85, cob} for a full account). In
fact one of the main reasons for  introducing
cyclic cohomology was to define a Chern character in $K$-homology
\cite{ac81}. In the even case, let Let $(H, \, F)$ be an even $p$-summable
Fredholm module over an algebra $\A$ as in definition \eqref{fm}. For each even integer $n \geq p-1$, define an $n$-cochain
$\varphi_n$ on $\A$ by
$$\varphi_n (a_0, \cdots, a_n)= \, \text{Trace}\, (\varepsilon \, a_0 [F,\,
a_1]\cdots  [F, \, a_n]).$$ The $p$-summability condition on $(H,
\, F)$ ensures that the above product of commutators is in fact a
trace class operator and $\varphi_n$ is finite for all $n \geq
p-1$ (this is obvious if $n> p-1$; in general one has to
manipulate the commutators a bit to prove this).
\begin{exercise} Show that $\varphi_n$ is a cyclic cocycle. Also show that if $n$ is odd then $\varphi_n =0$.
\end{exercise}

Although this definition depends on $n$, it can be shown that the cyclic cocycles $\varphi_n$ are related to each other by
the periodicity operator $S$,
$$ S \varphi_n =\varphi_{n+2},$$
and therefore define an even periodic cyclic cohomology class.
This is Connes' Chern character  of a $K$-homology class in the
even case. For applications it is important to have an index
formula  which computes the value of the pairing $ \langle  [\varphi_n],
\, [e] \rangle $ as the index of a Fredholm operator similar to Example
2.11.   We refer to \cite{co85, cob} for this.

Let  $A$ be a  nuclear $C^*$-algebra.  In the odd case, smooth $p$-summable elements of $K$-homology group $K^1 (A)$
can be represented either by
 smooth Brown-Douglas-Fillmore extensions
$$ 0\longrightarrow \mathcal{L}^p \longrightarrow \mathcal{E}
\longrightarrow \A \longrightarrow 0$$ or by Kasparov modules.
  We refer to \cite{co85} for the definition  of  the Connes-Chern character in the odd case. We have already
met one example of this though in the case of smooth Toeplitz
extension
$$ 0\longrightarrow \mathcal{K}^{\infty} \longrightarrow
\mathcal{T}^{\infty}
\longrightarrow C^{\infty} (S^1) \longrightarrow 0.$$
The map $f \to T_f$, sending a function on the circle to the corresponding Toeplitz operator,  is a section for
the symbol map $\sigma$. The extension is $p$-summable for all $ p\geq
1$. Its  Connes character is  represented by  the cyclic 1-cocycle on  $C^{\infty} (S^1)$ defined by
$$ \varphi_1 (f, \, g) = \,\text{Tr} \, ([T_f, \, T_g]).$$

\subsection{Cyclic modules}
  Cyclic cohomology of  algebras was
first defined by Connes  through explicit
complexes or bicomplexes \cite{ac81, co85}. Soon after   he introduced the notion of  {\it cyclic
 module} and defined its cyclic cohomology \cite{ac83}. Later developments
 proved that this extension was of great significance. Apart from earlier applications, here we have the very
 recent
 work \cite{cocoma}  in mind where the abelian category of
 cyclic modules plays the role of
 the  category of motives in noncommutative geometry.
   Another  recent example is the cyclic
 cohomology of Hopf algebras \cite{como2, como3, hkrs1, hkrs2},  which cannot be defined as the cyclic cohomology of
an algebra or a coalgebra but only   as the cyclic cohomology of a cyclic
module naturally attached to the given Hopf algebra.

The original motivation of \cite{ac83} was
 to define  cyclic cohomology of algebras as a derived functor. Since the category
 of algebras and algebra homomorphisms is not even an additive category (for the simple
 reason  that the sum of two algebra homomorphisms is not an algebra homomorphism in general),
 the standard
(abelian) homological algebra is not applicable. In Connes'
approach,  the category $\Lambda_k$  of
 cyclic $k$-modules appears as an  ``abelianization" of the category of
$k$-algebras.
  Cyclic cohomology is then shown to be the derived
functor of the  functor of traces,  as we explain in this section.

The {\it simplicial category} $\Delta$ is the category whose
objects are totally ordered sets
 $$[n]=\{0< 1< \cdots < n\},$$
 for $n=0, 1, 2, \cdots$. A morphism $f: [n] \to [m]$ is an order preserving, i.e.
 monotone non-decreasing, map
 $f: \{0, 1, \cdots, n\}\to \{0, 1, \cdots, m \}. $ Of particular
 interest among the morphisms of $\Delta$ are  {\it faces} $\delta_i$
 and {\it degeneracies} $\sigma_i$,
 $$ \delta_i: [n-1]\to [n], \quad \sigma_i: [n]\to [n-1], \quad
 \quad i= 1, 2, \cdots $$
 By definition $\delta_i$ is the unique injective morphism missing $i$ and
 $\sigma_i$ is the unique surjective  morphism identifying $i$ with $i+1$. It can be checked that they satisfy
 the following {\it simplicial identities}:

\begin{eqnarray*}
 \delta_j  \delta_i &=& \delta_{j-1}  \delta_i \quad
 \text{if} \quad i <j,\\
 \sigma_i  \sigma_i &=& \sigma_i \sigma_i \quad
 \text{if} \quad i <j,\\
\sigma_i \delta_i &=&
 \begin{cases}
\sigma_{j-1} \delta_i  &\text{$i<j$}\\
\text{id}      &\text{$i=j$ or $i=j+1$}\\
\sigma_j \delta_{i-1} &\text{$i>j+1$}.
\end{cases}
\end{eqnarray*}
Every morphism of $\Delta$ can be uniquely decomposed as a product of faces followed by a product of degeneracies.

The cyclic category $\Lambda$ has the same set of objects as $\Delta$
and in fact contains $\Delta$ as a subcategory. An, unfortunately
unintuitive,
definition of its morphisms  is as follows (see \cite{cob} for a more intuitive definition in terms of homotopy
classes of maps from $S^1 \to S^1$). Morphisms of   $\Lambda$ are generated by simplicial morphisms  $\delta_i$,
$\sigma_i $ as above and   $\tau_n : [n] \to [n]$ for $n \geq 0$. They are  subject to the above simplicial
as well as  the following
extra relations:

\begin{eqnarray*}
\tau_n\delta_{i}&=&\delta_{i-1} \tau_{n-1} \hspace{43 pt} 1\le i\le
n\\
\tau_n \delta_0 &=& \delta_{n}
\\ \tau_n \sigma_i &=& \sigma _{i-1} \tau_{n-1}  \hspace{43 pt} 1\le i\le n\\
\tau_n \sigma_0 &=& \sigma_n \tau_{n+1}^2 \\
\tau_n^{n+1} &=& \mbox{id}  .
\end{eqnarray*}

 A {\it cyclic object} in a category $\mathcal{C}$ is
a  functor $\Lambda^{\text{op}}  \rightarrow \mathcal{C}$. A {\it cocyclic object} in
 $\mathcal{C}$ is a  functor $\Lambda \rightarrow \mathcal{C}$.
For any commutative  ring $k$, we denote the category of cyclic
  $k$-modules by $\Lambda_k $. A morphism of cyclic $k$-modules is a natural transformation between the corresponding
  functors. Equivalently, a morphism  $f: X \to Y$ consists of a sequence
  of $k$-linear maps $f_n: X_n \to Y_n$ compatible with the
  face, degeneracy, and cyclic operators. It is clear that
  $\Lambda_k$ is an abelian category. The kernel and cokernel of
  a morphism  $f$ is defined pointwise: $(\text{Ker}\, f)_n= \text{Ker}\, f_n: X_n \to Y_n$ and $(\text{Coker}\, f)_n=
  \text{Coker}\,  f_n: X_n \to Y_n$.  More generally,  if $\mathcal{A}$ is any abelian category then the category
  $\Lambda \mathcal{A}$  of cyclic objects in $\mathcal{A}$ is itself
  an abelian category.

Let $Alg_k$ denote the category of unital $k$-algebras and unital
algebra homomorphisms. There is a functor
$$ \natural : Alg_k \longrightarrow \Lambda_k,$$
defined as follows. To an algebra $A$, we associate the
cyclic module
$A^ \natural$
 defined by $A_n^\natural =A^{\otimes(n+1)},  n\geq 0$, with  face,  degeneracy
and cyclic operators  given  by
\begin{eqnarray*}
\delta_i(a_0 \otimes a_1\otimes \dots \otimes a_n)&=&a_0 \otimes \dots \otimes
a_{i}a_{i+1}\otimes \dots \otimes a_n\\
\delta_n(a_0 \otimes a_1\otimes\dots \otimes a_n)&=&a_na_0 \otimes a_1 \otimes
\dots  \otimes a_{n-1}\\
\sigma_i(a_0 \otimes a_1\otimes\dots \otimes a_n)&=&a_0 \otimes\dots \otimes
a_i \otimes  1 \otimes \dots\otimes
 a_n\\
\tau_n(a_0 \otimes a_1\otimes\dots \otimes a_n)&=& a_n \otimes a_0\dots \otimes a_{n-1}.
\end{eqnarray*}
A unital algebra map $f: A \to B$ induces a morphism of cyclic modules
$f^{\natural}: A^ \natural \to B^ \natural$ by $f^\natural (a_0\otimes
\cdots \otimes a_n)= f(a_0)\otimes
\cdots \otimes f(a_n).$

\begin{example} We have
$$\text{Hom}_{\Lambda_k} \,(A^ \natural, \, \, k^ \natural) \simeq T(A),$$
where $T(A)$ is the space of traces from $A \to k$. Under this
isomorphism  a trace $\tau$ is sent to the cyclic map $(f_n)_{n\geq 0}$,  where
$$ f_n(a_0\otimes a_1\otimes \cdots \otimes a_n)=\tau (a_0a_1\cdots
a_n), \quad \, n \geq 0.$$
\end{example}

Now we can state the following fundamental theorem of Connes
\cite{ac83}:
\begin{theorem} For any unital  $k$-algebra $A$, there is a    canonical
isomorphism
$$ HC^n(A) \simeq Ext^n_{\Lambda_k} (A^ \natural, \, k^ \natural), \quad \quad \text{for all} \, \, n\geq 0.$$
\end{theorem}
Now the above Example and Theorem, combined together, say that cyclic cohomology is the derived functor of the functor
of traces $A \to T(A)$ where the word derived functor is understood to
mean as above.

Motivated by the above theorem, one defines the cyclic cohomology and
homology of any cyclic module $M$  by
$$ HC^n (M):= \text{Ext}^n_{\Lambda_k} (M,\, k^ \natural),$$
and
$$HC_n(M):=\text{Tor}_n^{\Lambda_k} (M, \,k^\natural),$$
One can use the injective resolution used to prove  the above Theorem  to show that
 these Ext and Tor groups can be computed by explicit complexes and bicomplexes, similar to the situation with algebras. For example
one has the  following first quadrant  bicomplex, called
 the {\it cyclic bicomplex}  of  $M$
$$\begin{CD}
\vdots @.\vdots @.\vdots @.\\
M_2@<1-\lambda<< M_2 @<N<< M_2@<1-\lambda<< \dots  \\
@VV bV @VV-b'V @VV bV  \\
M_1@<1-\lambda<< M_1 @<N<< M_1@<1-\lambda<< \dots \\
@VV bV @VV-b'V @VV bV  \\
M_0@<1-\lambda<< M_0 @<N<< M_0@<1-\lambda<< \dots
\end{CD}
$$
whose total homology is naturally isomorphic to cyclic homology.
Here the operator $\lambda : M_n \to M_n$ is defined by $\lambda
=(-1)^n\tau_n$, while
$$
 b = \sum _{i=0}^n(-1)^i \delta _i, \quad \quad
 b' = \sum _{i=0}^{n-1} (-1)^i \delta _i,$$
 and
$N= \sum_{i=0}^n \lambda ^i$.
Using the simplicial and cyclic relations, one can check that $b^2=b'^2=0$,
$b(1-\lambda)=(1-\lambda)b'$
 and $b'N=Nb'$. These relations amount to saying that the above is a
 bicomplex.

The $(b, B)$-bicomplex of a cyclic module is the bicomplex
$$\begin{CD}
\vdots @.\vdots @.\vdots \\
M_2@<B<< M_1@<B<< M_0 \\
@VVbV @VVbV\\
M_1@<B<<M_0\\
@VVbV\\
M_0
\end{CD}
$$
whose total homology is again isomorphic to the cyclic homology of
$M$ (this time we have to assume that  $k$ is a field of
characteristic 0). Here  $B : M_n \rightarrow M_{n+1}$ is Connes'
boundary operator defined  by $B=(1-\lambda) sN$,
where $s= (-1)^n \sigma _n. $

 A remarkable property of the  cyclic category $\Lambda$, not shared  by the simplicial category,
  is its {\it self-duality} in the sense that there is a
natural isomorphism  of categories $\Lambda \simeq \Lambda^{\text{op}}$.   Roughly  speaking, Connes'  duality
functor $\Lambda^{\text{op}} \longrightarrow
\Lambda $ acts as the identity on objects
of $\Lambda$ and exchanges face and degeneracy operators while sending
the cyclic operator to its inverse.
Thus to a cyclic (resp. cocyclic) module one can
associate a cocyclic (resp. cyclic) module by applying Connes' duality
isomorphism.
In the next section we shall see  examples of cyclic modules in   Hopf cyclic
(co)homology that are dual to each other in the above sense.

\subsection{Hopf cyclic cohomology}

In their fundamental work on  index theory of transversally elliptic operators \cite{como2},  Connes and Moscovici
 developed  a new cohomology theory for Hopf algebras based on ideas in cyclic cohomology. This theory
 can be regarded  as the right  noncommutative  analogue    of both
 group and Lie algebra homology, although this was not the original motivation behind it. Instead,
      the main reason
   was to obtain a
  {\it noncommutative characteristic map}
\begin{equation} \label{hcm}
  \chi_{\tau}: {HC}_{(\delta,\sigma)}^\ast (H) \longrightarrow  HC^\ast (A),
  \end{equation}
  for an action of a  Hopf
  algebra $H$ on an  algebra $A$ endowed with an ``invariant trace'' $\tau:
  A \to \mathbb{C}$. Here, the pair $(\delta, \sigma)$ consists of a grouplike element $\sigma \in H$ and
   a character $\delta: H\to \mathbb{C}$ satisfying certain
   compatibility conditions to be discussed later in this section. While in this section we confine ourselves to Hopf algebras,  we
   refer to the recent surveys \cite{como6} and \cite{kr3} for later developments
   in the subject inspired by \cite{como2}.

   The  characteristic map \eqref{hcm} is induced,  on the level of
   cochains, by a map
   $$ \chi_{\tau}: H^{\otimes n} \longrightarrow C^n(A)$$
   defined by
   \begin{equation}\label {cm}
   \chi_{\tau} (h_1 \otimes \cdots \otimes h_n)(a^0, \cdots, a^n)=\tau (a_0 h_1(a^1)\cdots h_n(a^n)).
   \end{equation}
   Maps like this  have quite a history in cyclic cohomology, going back to \cite{ac80}. Notice, for example,  that
   the fundamental cyclic cocycles \eqref{nctc1}, \eqref{nctc2}  on the noncommutative
   2-torus  are of this form, where  $H$ is the enveloping algebra of a two-dimensional abelian Lie algebra
   acting by a pair of commuting derivations $\delta_1$ and $\delta_2$
   on $\mathcal{A}_{\theta}$.

   \begin{exercise} Let $\delta_1, \cdots, \delta_p$ be a commuting family of derivations on an algebra $\A$ and let $\tau : \A
   \to \mathbb{C}$ be an invariant trace, i.e. $\tau \, (\delta_i (a))=0$ for all $a\in \A$ and $i=1, \cdots, p$. Show that
   $$\varphi (a^0, \cdots, a^p)=\sum_{\sigma \in S_n} (-1)^{\sigma} \tau (a^0
   \delta_{\sigma_1}(a^1)\cdots \delta_{\sigma_p}(a^p))$$
   is a cyclic $p$-cocycle on $\A$.
   \end{exercise}
 Applied to higher dimensional  noncommutative tori, these cocycles give
 a  basis for its periodic cyclic cohomology.

 For applications to
transverse geometry and number theory \cite{como4, como5, como6}, it is important to
formulate a notion of `invariant trace' under the presence of a modular
pair.  Let $A$ be an $H$-module algebra, $\delta $ a character of
$H$, and $\sigma \in H$ a grouplike element.  A linear map
$\tau :A \rightarrow \mathbb{C}$ is called $\delta$-{\it invariant} if for
all $h \in H$ and $a\in A$,
$$ \tau (h(a))=\delta (h) \tau (a).$$
$\tau$ is called a $\sigma$-{\it trace} if for all $a, b$ in $A$,
$$ \tau (ab)=\tau (b \sigma (a)).$$
 For $a, b \in A,$ let
$$\langle a, \, b \rangle :=\tau (ab).$$
 Then the $\delta$-invariance property of $\tau$ is equivalent to
 the {\it integration by parts formula}:
\begin{equation} \label{ibp}
\langle h(a), \, b \rangle = \langle a, \,  \widetilde{S}_{\delta}(h)(b)\rangle,
\end{equation}
where the  $\delta$-{\it twisted antipode}
$\widetilde{S}_{\delta}: H \rightarrow H$ is defined by $\widetilde{S}_{\delta}=\delta *
S$. That is,
$$\widetilde{S}_{\delta}(h)= \delta (h^{(1)})S(h^{(2)}).$$
Loosely speaking, this   amounts to saying that  the formal adjoint of the differential operator $h$ is
$\widetilde{S}_{\delta}(h)$.
 Following \cite{como2, como3}, we
say  $(\delta, \, \sigma)$ is a {\it modular pair} if $\delta
(\sigma)=1$, and a {\it modular pair in involution} if
$$  \widetilde{S}_{\delta}^2 (h)=\sigma
h\sigma^{-1},$$
for all $h$ in $H$.

\begin{examples} \rm{ 1. For any commutative or cocommutative Hopf algebra  we have $S^2=1$. It follows that
$(\varepsilon, 1)$ is a
modular pair in involution. \\

\noindent 2. The original non-trivial example of a modular pair in involution
is the pair $(\delta, \, 1)$ for the 
Connes-Moscovici  Hopf algebra $\mathcal{H}_1$. Let $\delta $
denote the unique extension of the modular character
$$\delta : \mathfrak{g}_{\text{aff}} \to \mathbb{R};  \quad \delta(X)=1, \quad  \delta(Y)=0,$$
to a character $\delta : U (\mathfrak{g}_{\text{aff}})\to \mathbb{C}.$
There is a unique extension of $\delta $ to a character, denoted
by the same symbol  $\delta : \mathcal{H}_1 \to \mathbb{C}.$
Indeed the relations $[Y, \delta_n]=n\delta_n$ show that we must
have $\delta (\delta_n)=0$, for $n=1, 2, \cdots.$ One can then
check that these relations are compatible with the algebra
structure of $\mathcal{H}_1.$

Now the   algebra
$A_{\Gamma}=C^{\infty}_0(F^+(M))\rtimes \Gamma$ from Section 4.2
 admits a $\delta$-invariant trace $\tau : A_{\Gamma} \to
\mathbb{C}$ under its canonical $\mathcal{H}_1$ action. It is  given by
\cite{como2}:
$$\tau (fU^{\ast}_{\varphi})=\int_{F^+(M)}f(y, y_1)\frac{dy dy_1}{y_1^2},
\quad  \text{if} \; \varphi =1,$$
and $\tau (fU^{\ast}_{\varphi})=0$, otherwise.\\

\noindent 3. Let $H=A(SL_q(2))$ denote the Hopf algebra of functions on quantum $SL(2)$.
 As an algebra it is generated by  $x,\; u,\; v,\; y,$  subject to the  relations
$$ ux=qxu, \;\; vx=qxv, \;\; yu=quy,\;\;yv=qvy,$$
$$uv=vu,\;\; xy-q^{-1}uv=yx-quv=1. $$
The coproduct, counit and antipode of  $H$ are defined by
$$\Delta (x)=x \otimes x+u \otimes v,\;\;\;\Delta (u)=x \otimes u+u \otimes y, $$
$$\Delta (v)=v \otimes x+y \otimes v,\;\;\;\Delta (y)=v \otimes u+y \otimes y, $$
$$\epsilon (x)=\epsilon (y)=1,\;\;\;\epsilon (u)=\epsilon (v)=0,$$
$$S(x)=y,\;\; S(y)=x,\;\;S(u)=-qu,\;\;S(v)=-q^{-1}v. $$
Define a character  $\delta : H \to  \mathbb{C} $ by
$$\delta (x)=q,\; \; \delta(u)=0,\; \; \delta (v)=0,\; \; \delta (y)=q^{-1}.$$
One checks that  $\widetilde {S}_{\delta}^2= \text{id}$. This shows that $(\delta, \, 1)$ is a modular pair in involution for
$H$. This example and its Hopf cyclic cohomology are studied in
\cite{kr1}.

More generally, it is shown in \cite{como3} that {\it coribbon Hopf algebras} and compact
quantum groups are endowed with canonical modular pairs in involution of the form $(\delta, \, 1)$
and,
dually,  ribbon Hopf algebras have canonical modular pairs in involution of the type
$(1, \sigma)$.
\\

\noindent 4. It is shown in \cite{hkrs1}  that modular pairs in
involution are in fact
 one-dimensional examples of {\it  stable anti-Yetter-Drinfeld modules} over Hopf algebras introduced there.
  These modules are
noncommutative coefficient systems for the general Hopf cyclic
cohomology theory developed in \cite{hkrs2}. }
\end{examples}

  Now let $(H, \, \delta, \, \sigma)$ be a Hopf algebra endowed with a modular
 pair in involution. In \cite{como2}
Connes and Moscovici attach  a cocyclic module  $H_{(\delta,
\sigma)}^\natural$ to this data
 as follows.
Let
$$H_{(\delta,\sigma)}^{\natural,0}= \mathbb{C} , \quad \text{ and} \quad
H_{(\delta,\sigma)}^{\natural,n}=H^{\otimes n},  \quad  \text{for}\, \, \, \,   n\geq 1.$$
Its  face, degeneracy and cyclic operators $\delta_i$,  $\sigma_i$,  and $\tau_n$ are
defined by
\begin{eqnarray*}
\delta_0(h_1 \otimes \dots \otimes h_n)&=& 1  \otimes h_1 \otimes\dots
 \otimes h_n \\
\delta_i(h_1 \otimes\dots \otimes h_n)&=& h_1 \otimes \dots \otimes\Delta (h_i)
\otimes\dots \otimes h_n
  \;\;\text{for}\;\;1\leq i \leq n \\
\delta _{n+1}(h_1 \otimes\dots \otimes h_n)&=& h_1 \otimes\dots \otimes h_n
\otimes \sigma \\
\sigma_i(h_1 \otimes\dots \otimes h_n)&=& h_1 \otimes\dots \otimes\epsilon
(h_{i+1})\otimes\dots \otimes h_n
 \;\;\text{for}\;\;0 \leq i \leq n \\
\tau_n (h_1 \otimes\dots \otimes h_n)&=&\Delta^{n-1}\widetilde{S}(h_1)\cdot(h_2 \otimes\dots \otimes h_n \otimes \sigma).
\end{eqnarray*}

The cyclic cohomology of the cocyclic module
$H_{(\delta,\sigma)}^\natural$ is called the Hopf cyclic cohomology of
the triple $(H, \,  \delta, \, \sigma)$ and will be denoted by
 $HC^n_{(\delta,\sigma)}( H)$.

 \begin{examples} \rm{
1. For  $H=\mathcal{H}_n$,   the  Connes-Moscovici Hopf algebra, we have
\cite{como2} (cf. also \cite{mora})
$$HP^n_{(\delta,1)}(\mathcal{H}_n)\simeq \bigoplus_{i=n\; (\text{mod}\;2)} H^i (\mathfrak{a}_n,\mathbb{C})$$
where $\mathfrak{a}_n$ is the Lie algebra of formal vector fields on
$\mathbb{R}^n$.\\

\noindent 2. For  $H=U(\mathfrak{g})$  the enveloping algebra of a Lie algebra  $\mathfrak{g}$, we
have \cite{como2}
$$HP^n_{(\delta,1)}(H)\cong \bigoplus_{i=n\; (\text{mod}\; 2)} H_i(\mathfrak{g},\mathbb{C}_\delta)$$

\noindent 3. For  $H=\mathbb{C}\lbrack G \rbrack $ the coordinate
ring of a nilpotent affine algebraic  group $G$, we have
\cite{como2}
$$HP^n_{(\epsilon,1)}(H)\cong \bigoplus_{i=n\;(\text{mod}\; 2)} H^i(\mathfrak{g},\mathbb{C}),$$
where $\mathfrak{g}=Lie(G)$.\\

\noindent 4. If $H$ admits a normalized left Haar integral, then \cite{cr}
$$HP^1_{(\delta,\sigma)}(H)=0,\qquad HP^0_{(\delta,\sigma)}(H)=\mathbb{C}.$$
Recall that a linear map $\int : H \rightarrow \mathbb{C}$ is called a normalized left  Haar integral
if for all $h\in H$,   $\int h=\int(h^{(1)})h^{(2)}$ and
$\int 1=1$. It is known that a Hopf algebra defined over a field admits
a normalized left Haar integral if and only if it is cosemisimple
\cite{sw}.
Compact quantum groups
 and group  algebras are known to admit a normalized Haar
 integral in the above sense. In the latter case $\int : \mathbb{C} G \rightarrow k$
sending $g \mapsto 0$ for all $g\neq e$  and $e \mapsto 1$ is a Haar
integral. Note that $G$ need not  be finite. In this regard, we should also mention
that there are interesting examples of finite-dimensional
non-cosemisimple Hopf algebras defined as quantum groups at roots of
unity. Nothing is known about the cyclic
(co)homology of these Hopf algebras.

\noindent 5. If $H=U_q(\text{sl}_2 )$ is the quantum universal algebra of $\text{sl}_2$, we have~\cite{cr}
$$HP^0_{(\epsilon,\sigma)}(H)=0,\quad  HP^1_{(\epsilon,\sigma)}(H)=\mathbb{C} \oplus \mathbb{C}.$$
\noindent 6. Let $H$ be a commutative Hopf algebra. The periodic cyclic cohomology of the cocyclic module
        $\mathcal{H}^\natural_{(\epsilon,1)}$  can be  computed in terms of the Hochschild homology of the coalgebra
         $H$ with trivial
        coefficients.

  \begin{prop}(\cite{kr1})
  Let $H$ be a commutative Hopf algebra. Its periodic cyclic
  cohomology in the sense
  of Connes-Moscovici is given by
  $$HP_{(\epsilon,1)}^n(H)=\bigoplus_{i=n\;(\text{mod}\;2)}H^i
  (H, \mathbb{C}).$$
    \end{prop}
    For example, if $ H=\mathbb{C} \lbrack G\rbrack$ is the algebra of
    regular functions on
     an affine algebraic group $G$, the coalgebra complex of
     $ H$ is isomorphic
     to the group cohomology complex of $G$ where  instead
     of regular cochains one uses regular functions
     $G\times G\times \dots \times G \rightarrow k\mathbb{C}$. Denote this cohomology by
      $H^i(G, \mathbb{C})$. It follows that
   $$HP_{(\epsilon,1)}^n( \mathbb{C} \lbrack G\rbrack)=\bigoplus_{i=n\;(\text{mod}\;2)}H^i(G, \mathbb{C}).$$
   As is remarked in \cite{kr1},  when  the Lie algebra Lie$(G)=\mathfrak{g}$
   is nilpotent, it  follows from
   Van Est's theorem
   that $H^i(G, \mathbb{C})\simeq H^i(\mathfrak{g}, \mathbb{C})$. This gives an alternative proof of
    Proposition 4 and
   Remark 5 in~\cite{como2}.}
   \end{examples}

   Now given $(H, \, \delta, \, \sigma)$,  a Hopf algebra endowed with a modular pair in involution as above,
   let $\A$ be an algebra with an $H$-action and let $\tau : \A \to \mathbb{C}$ be a
$\delta$-invariant $\sigma$-trace. Then one can check, using in
particular the integration by parts formula  \eqref{ibp},   that
the characteristic map  \eqref{cm} is a morphism of cocyclic
modules. It follows that we have  a well defined map
\begin{eqnarray*}
 \chi_{\tau} : HC_{(\delta,\sigma)}^n(H)\rightarrow HC^n(\A)
 \end{eqnarray*}
\begin{example} \rm{ Let $\mathcal{A}$ be an  $n$-dimensional smooth noncommutative torus
with canonical commutating derivations $\delta_1, \cdots,
\delta_n$ defined by  $\delta_i (U_j)= \delta_{ij} U_i$. They
define an action of $H=U(\frak{g})$ on $\A$ where $\frak{g}$ is
the abelian $n$-dimensional Lie algebra. The canonical trace
$\tau$ is invariant under this action. The  characteristic
 map $\chi_{\tau}$
combined with the antisymmetrization map
$${\bigwedge}^k \frak{g} \to U(\frak{g})^{\otimes k}$$
defines a map
$${\bigwedge}^k \frak{g} \to HC^k (\A).$$
This is the map of  Exercise 6.8. }
\end{example}

In the rest of   this  section  we recall  a dual cyclic theory
for Hopf
   algebras which was defined and studied in
   \cite{kr1}.
    This theory is needed, for
   example, when one studies
    coactions of Hopf algebras and  quantum
   groups on noncommutative spaces. Notice that for compact quantum
   groups coactions are more natural.
   Also, as we mentioned before, for cosemisimple Hopf algebras,  i.e.
   Hopf algebras endowed with a  normalized   Haar integral, the Hopf
    cyclic cohomology   is trivial in positive dimensions,  but
    the dual theory is non-trivial. There is a clear analogy with
    continuous group cohomology here.

Let  $(\delta, \sigma) $ be a modular pair on $H$ such that
$\widehat{S}^2=\text{id}_H$, where $\widehat{S} (h):= h^{(2)}\sigma S(h^{(1)})$.  We define a cyclic module
$\widetilde{ H}_{\natural}^{(\delta,\sigma)} $  by
$$  \widetilde{H}_{\natural, 0}^{(\delta,\sigma)}= \mathbb{C}, \quad \quad
\widetilde{H}_{\natural, n}^{(\delta,\sigma)}  = H^{\otimes n}, \,
 \, \, \,  n> 0. $$
   Its face,  degeneracy, and cyclic operators are defined  by
\begin{eqnarray*}
 {\delta}_0 (h_1 \otimes h_2 \otimes  ... \otimes h_n) &=&\epsilon (h_1)h_2 \otimes h_3
  \otimes ... \otimes h_n \\
 {\delta}_i (h_1 \otimes h_2 \otimes
 ... \otimes h_n) &=& h_1\otimes h_2 \otimes ...\otimes
h_i h_{i+1}\otimes ... \otimes h_n   \\
{\delta}_n(h_1 \otimes h_2 \otimes  ... \otimes h_n) &=& \delta (h_n)h_1\otimes
 h_2 \otimes ... \otimes h_{n-1} \\
 {\sigma}_0(h_1 \otimes h_2 \otimes  ... \otimes h_n) &=&1 \otimes h_1  \otimes
 ... \otimes h_n \\
\hspace{2cm}{\sigma}_i(h_1 \otimes h_2 \otimes ... \otimes h_n) &=&h_1 \otimes
h_2 ...\otimes h_i \otimes 1 \otimes
 h_{i+1} ... \otimes h_n     \\
 {\sigma }_n(h_1 \otimes h_2 \otimes  ... \otimes h_n) &=& h_1 \otimes h_2
 \otimes ...\otimes 1, \\
\tau _n(h_1 \otimes h_2 \otimes  ... \otimes h_n) &=&
 \delta(h_n^{(2)})\sigma S (h_1^{(1)} h_2^{(1)} ...h_{n-1}^{(1)}h_n^{(1)})
\otimes h_1^{(2)} \otimes ...\otimes h_{n-1}^{(2)}.
\end{eqnarray*}
We denote the cyclic homology of this cyclic module by $\widetilde{HC}^{(\delta ,\sigma)}_\bullet
(H).$

\begin{remark}It is not difficult to check that $(\delta\circ
S^{-1},\sigma^{-1})$, is a modular pair in involution if and only
if $(\delta,\sigma)$ is a modular pair with
$\widehat{S}^2=id_{H}$. In other words $(\delta,\sigma)$
is a modular pair in involution in the
 sense of Connes and Moscovici \cite{como2} if and only if $(\delta\circ S,\sigma^{-1})$ is a modular pair in
 involution in the sense of
 \cite{kr1}.
 \end{remark}

Now let $A$ be an $H$-comodule algebra.  A linear map, $\tau
:A\rightarrow \mathbb{C}$ is called a
$\delta$-trace  if
\begin{equation*}
\tau(ab)=   \tau(b^{(0)}a)\delta (b^{(1)}) \hspace{2cm}\forall a ,
b\in  A.
\end{equation*}
It is called $\sigma $-invariant if for  all $a  \in  A$,
$$\tau(a^{(0)})a^{(1)}=\tau(a)\sigma.$$
Now consider the map $\chi_{\tau} :A_\natural \rightarrow \widetilde{H}^
{(\delta,\sigma)}_\natural$ defined by
$$\chi_{\tau}(a_0\otimes a_1\otimes\cdots \otimes a_n)=\tau(a_0a_1^{(0)}\cdots a_n^{(0)})
a_1^{(1)}\otimes
  a_2^{(1)}\otimes\cdots \otimes a_n^{(1)}.$$
It is proved in \cite{kr1} that $\chi_{\tau}$ is a morphism of
cyclic modules. This looks rather uninspiring, but once dualized
we obtain a characteristic map for coactions
$$\chi_{\tau}^*: \widetilde{HC}_{(\delta ,\sigma)}^n (H)\longrightarrow
HC^n(A),$$
which can be useful as will be shown  in Example \eqref{dcmg} below.

Next we state a theorem which computes  the Hopf cyclic homology of cocommutative
Hopf algebras in terms of the Hochschild homology of the underlying
algebra:
   \begin {theorem}(\cite{kr1})\label{t1}
If $H$ is a cocommutative Hopf algebra, then
$$ \widetilde{HC}^{({\delta,1})}_n(H) = \bigoplus _{i\geq 0} H_{n-2i}
(H, \,   \mathbb{C}_\delta),$$  where $\mathbb{C}_\delta$ is the one-dimensional module defined by $\delta$.
\end{theorem}

\begin{example}\rm{
One knows that for any Lie algebra  $\mathfrak{g}$,
$$H_n(U(\mathfrak{g}), \, \mathbb{C}_\delta) \simeq  H_n^{Lie}(\mathfrak{g}, \, \mathbb{C}_\delta).$$
 So by Theorem \textnormal{\ref{t1}} we have
$$\widetilde{HC}^{(\delta,1)}_n(U(\mathfrak{g})) \simeq \bigoplus _{i\geq 0}H_i^{Lie}(\mathfrak{g},  \mathbb{C}_\delta). $$
}
\end{example}

\begin{example} \label{dcmg}
\rm{Let $H=\mathbb{C}\Gamma$ be the group algebra of a  discrete group  $\Gamma$. Then from Theorem \textnormal{\ref{t1}} we have
\begin{center}
$\widetilde{HC}_n^{(\epsilon,1)}(\mathbb{C}\Gamma)\simeq\bigoplus _{i\ge 0}H_{n-2i}(\Gamma, \mathbb{C}),$\\
and
$\widetilde{HP}_n^{(\epsilon,1)}(\mathbb{C}\Gamma)\simeq\bigoplus _{i=n\;(\text{mod}\; 2)}H_{i}(\Gamma, \mathbb{C}).$
\end{center}

Now any Hopf algebra
$H$ is a comodule algebra over itself
 via the coproduct map $H\longrightarrow
H\otimes H$. The map
 $\tau :\mathbb{C}\Gamma \rightarrow \mathbb{C}$  defined by
 $$\tau (g)=\left\{\begin{matrix}
  1 \quad\qquad g=e\\
0 \quad\qquad g \neq e \end{matrix}\right.$$
is a $\delta$-invariant $\sigma$-trace for $\delta=\epsilon$, $\sigma=1$.
The dual characteristic map
$$\chi_{\tau}^{\ast}: \widetilde{HC}^n_{(\epsilon,1)}(\mathbb{C}\Gamma)\rightarrow
{HC}^n(\mathbb{C}\Gamma)$$
 combined with the inclusion
$H^n(\Gamma, \mathbb{C})\hookrightarrow \widetilde{HC}^n_{(\epsilon, 1)}
(\mathbb{C}\Gamma)$ gives us a map
 $$H^n(\Gamma,\mathbb{C})\rightarrow HC^n(\mathbb{C}\Gamma)$$
 from group cohomology to cyclic cohomology. The image of a normalized group $n$-cocycle $\varphi (g_1, \cdots, g_n)$ under
 this map is the cyclic $n$-cocycle $\hat{\varphi}$ defined by
$$\hat{\varphi}(g_0, \cdots, g_n)=\left\{\begin{matrix}
  \varphi (g_1, \cdots, g_n) \quad\qquad g_0g_1\cdots g_n=e\\
0  \qquad \qquad\qquad g_0g_1\cdots g_n \neq e \end{matrix}\right. $$
Thus the characteristic map for Hopf cyclic homology reduces to a well known map in noncommutative geometry \cite{cob}. This
should be compared with Example 6.5. }
\end{example}

It would be very interesting to compute the  Hopf cyclic homology
$\widetilde{HC}_n $ of
quantum groups. We cite one of the very few results known in this direction.
Let  $H=A(SL_q(2, \mathbb{C}))$ be the Hopf algebra of quantum $SL_2$.
 As an algebra it is generated by  $a, b, c, d,$ with  relations
$$ ba=qab, \;\; ca=qac, \;\; db=qbd,\;\;dc=qcd,$$
$$bc=cb,\;\; ad-q^{-1}bc=da-qbc=1. $$
The coproduct, counit and antipode of   $H$ are defined by
$$\Delta (a)=a\otimes a+b\otimes c,\;\;\;\Delta (b)=a\otimes b+b\otimes c $$
$$\Delta (c)=c\otimes a+d\otimes c,\;\;\;\Delta (d)=c\otimes b+d\otimes d $$
$$\epsilon (a)=\epsilon (d)=1,\;\;\;\epsilon (b)=\epsilon (c)=0,$$
$$S(a)=d,\;\; S(d)=a,\;\;S(b)=-qb,\;\;S(c)=-q^{-1}c. $$
  We  define  a modular pair $(\sigma, \delta)$ by
$$\delta (a)=q,\; \; \delta(b)=0,\; \; \delta (c)=0,\; \; \delta (d)=q^{-1},$$
 $\sigma =1$. Then we have $\widetilde S_{(1, \delta)}^2=\text{id}$.

\begin{theorem}(\cite{kr1})
For $q$ not a root of unity,   one has
 $$\widetilde{ HC}_1(A(SL_q(2, \mathbb{C})))=\mathbb{C}\oplus
 \mathbb{C},
 \quad   \quad
\widetilde{HC}_n(A(SL_q(2,  \mathbb{C})))=0,  \quad   n\neq 1,$$
and
$$\widetilde{HP}_0(A(SL_q(2, \mathbb{C})))
=\widetilde{HP}_1(A(SL_q(2, \mathbb{C})))=0.$$
\end{theorem}

In \cite{hkrs1, hkrs2}, following the lead of \cite{ak, kr1, kr2},  Hajac-Khalkhali-Rangipour-\\Sommerh\"auser define
 a full fledged Hopf cyclic cohomology theory for algebras or coalgebras endowed with actions
or coactions of a Hopf algebra. This extends the pioneering work
of Connes and Moscovici in two different directions. It  allows
 coefficients for the theory  and instead of Hopf algebras one
now works with algebras or coalgebras with a Hopf action. It turns out that  the
periodicity condition $\tau_n^{n+1}=\text{id}$ for the cyclic operator
puts very stringent conditions on the type of   coefficients that are
allowable
and the correct class of Hopf  modules turned to be
 the class of {\it stable anti-Yetter-Drinfeld modules} over a
 Hopf algebra. It also sheds light on Connes-Moscovici's modular pairs in involution by interpreting them
 as  one-dimensional stable anti-Yetter-Drinfeld modules.

 The category of anti-Yetter-Drinfeld modules over a Hopf algebra $H$ is a
 twisting, or variant rather,   of the category of Yetter-Drinfeld $H$-modules. Notice that this latter category is widely
 studied primarily
 because of its connections with quantum group theory and with invariants of knots and low-dimensional
 topology \cite{ye}.
 Technically it is obtained from the latter  by replacing the antipode $S$ by $S^{-1}$ although this
 connection is hardly illuminating. We refer to \cite{kr3} and
 references therein for  a survey of
 anti-Yetter-Drinfeld modules and Hopf cyclic cohomology.

\begin{remark} The picture that is emerging is quite intriguing and seems to be at the crossroads of three different areas: von Neumann
algebras,  quantum groups and low-dimensional topology, and cyclic
cohomology.  We have Connes-Moscovici's modular pairs in
involution which was suggested by type III factors and
non-unimodular Lie groups and turned out to be examples of
anti-Yetter-Drinfeld modules. As we saw above the latter category
is of fundamental importance in Hopf cyclic cohomology. One
obviously needs to understand these connections much better.
\end{remark}

\end{document}